\documentclass[lecture,11pt,a4paper]{pcms-l}

\makeatletter

\newif\iffirstlecture\firstlecturefalse

\newcommand{\lectureseries}{\firstlecturetrue
              \secdef\@lectureseries\@slectureseries} 

\newcommand{\@lectureseries}[2][default]{\chapter*{#2}%
              \gdef\thelectureseries{#1}} 

\newcommand{\@slectureseries}[1]{\chapter*{#1}}

\renewcommand{\auth}{\secdef\@auth\@sauth}

\newcommand{\@auth}[2][default]{\vspace{-1pc}{\raggedleft
        \Large\bf\noindent
        #2\endgraf
        \vspace*{2pc}
        }
        \def\@author{#1}%
}

\newcommand{\@sauth}[1]{\vspace{-1pc}{\raggedleft
        \Large\bf\noindent
        #1\endgraf
        \vspace*{2pc}
        }
        \def\@author{#1}%
}

\def\lecture#1{\global\Lecturetrue\global\Monographfalse
\iffirstlecture\else\chapter*{}\fi\firstlecturefalse
  \global\let\sectionmark\@gobble 
  \addtocounter{lecture}1\relax
  \refstepcounter{chapter}%
\gdef\thelecturename{#1\unskip}
  {\Large\bfseries
   \raggedleft
   \@xp\uppercase\@xp{\thelecturelabel} {\LARGE\thelecturenum}\\
   \vspace*{3pt}%
   \thelecturename
   \endgraf}%
  \let\@secnumber=\thelecturenum
  \@xp\lecturemark\@xp{\thelecturename}%
  \addcontentsline{toc}{chapter}{%
    \thelecturelabel\ \thelecturenum.\ \thelecturename}%
  \vspace{34\p@}\noindent}
  
\def\lecture{\global\Lecturetrue\global\Monographfalse
\iffirstlecture\else\chapter*{}\fi%
  \global\let\sectionmark\@gobble 
\secdef\@lecture\@slecture}

\def\@lecture[#1]#2{%
  \addtocounter{lecture}1\relax
  \refstepcounter{chapter}%
\gdef\thelecturename{#1\unskip}\firstlecturefalse
  {\Large\bfseries
   \raggedleft
   \@xp\uppercase\@xp{\thelecturelabel} {\LARGE\thelecturenum}\\
   \vspace*{3pt}%
    #2\unskip
   \endgraf}%
  \let\@secnumber=\thelecturenum
  \@xp\lecturemark\@xp{\thelecturename}%
  \addcontentsline{toc}{chapter}{%
    \thelecturelabel\ \thelecturenum.\ #2}%
  \vspace{34\p@}\noindent}
  
\def\slecturerunhead#1#2#3{%
    \let\@tempa\chaptername
    \uppercasenonmath{\@tempa}%
    \def\@tempb{#3\unskip}%
    \uppercasenonmath{\@tempb}%
    {\normalfont\@tempb}
    }
\def\slecturemark{
    \@secmark\markright\slecturerunhead\chaptername}%

\def\@slecture#1{%
\iffirstlecture
\gdef\thelecturename{#1\unskip}\firstlecturefalse
  {\Large\bfseries
\noindent\thelecturename
   \endgraf}%
  \let\@secnumber=\thelecturenum
  \@xp\slecturemark\@xp{\thelecturename}%
  \addcontentsline{toc}{chapter}{%
    \thelecturename}%
 \vspace{-6\p@}\noindent
\else
\gdef\thelecturename{#1\unskip}\firstlecturefalse
  {\Large\bfseries
   \raggedleft
   \@xp\uppercase\@xp{\thelecturename}
   \endgraf}%
  \let\@secnumber=\thelecturenum
  \@xp\slecturemark\@xp{\thelecturename}%
  \addcontentsline{toc}{chapter}{%
    \thelecturename}%
  \vspace{34\p@}\noindent
\fi}

\ifLecture
  \def\chapterrunhead#1#2#3{%
    \let\@tempa\@author
    \uppercasenonmath{\@tempa}%
    \uppercasenonmath{\thelectureseries}%
    \textmd{\@tempa, \thelectureseries}
    }
  \def\lecturerunhead#1#2#3{%
    \let\@tempa\chaptername
    \uppercasenonmath{\@tempa}%
    \def\@tempb{#3\unskip}%
    \uppercasenonmath{\@tempb}%
    \textmd{\@tempa\ #2. \@tempb}
    }
\else
  \let\chapterrunhead\partrunhead
\fi

\newif\ifBibliographyIsASection\BibliographyIsASectionfalse

  \def\bibliomark{
    \@secmark\markright\bibliorunhead\chaptername}%

  \def\bibliorunhead#1#2#3{%
    \let\@tempa\chaptername
    \uppercasenonmath{\@tempa}%
    \def\@tempb{#3\unskip}%
    \uppercasenonmath{\@tempb}%
    \textmd{\@tempb}
    }

\def\thebibliography#1{%
  \ifBibliographyIsASection
    \section*\refname
    \if@backmatter
      \markboth{\refname}{\refname}%
    \fi
  \else
\chapter*{}
  {\Large\bfseries
   \raggedleft
   \@xp\uppercase\@xp{\bibname} \\
   \endgraf}%
  \let\@secnumber=\thelecturenum
  \@xp\bibliomark\@xp{\bibname}%
  \addcontentsline{toc}{chapter}{%
    \bibname}%
  \vspace{34\p@}\noindent
  \fi
  \normalsize\labelsep .5em\relax
  \list{\@arabic\c@enumi.}{\settowidth\labelwidth{\@biblabel{#1}}%
  \leftmargin\labelwidth
  \advance\leftmargin\labelsep
	\usecounter{enumi}}\sloppy
  \clubpenalty9999 \widowpenalty\clubpenalty  \sfcode`\.\@m}

  \def\indexmark{
    \@secmark\markright\indexrunhead\chaptername}%

  \def\indexrunhead#1#2#3{%
    \let\@tempa\chaptername
    \uppercasenonmath{\@tempa}%
    \def\@tempb{#3\unskip}%
    \uppercasenonmath{\@tempb}%
    \textmd{\@tempb}
    }

\ifLecture
\def\theindex{\cleardoublepage
\@restonecoltrue\if@twocolumn\@restonecolfalse\fi
\columnseprule \z@ \columnsep 35pt
\def\indexchap{\@startsection
		{chapter}{1}{\z@}{8pc}{34pt}%
		{\raggedleft
		\Large\bfseries}}%
 \twocolumn[\indexchap[{\indexname}]{\@xp\uppercase\@xp{\indexname}}]
  \@xp\indexmark\@xp{\indexname}%
	\thispagestyle{plain}\let\item\@idxitem\parindent\z@
	 \footnotesize\parskip\z@ plus .3pt\relax\let\item\@idxitem}
\fi

\def\@makefntext{\noindent\@makefnmark}

\def\setaddress{%
  {\let\@makefnmark\relax  \let\@thefnmark\relax
        \nobreak
        \addressnum@=\z@
        \loop\ifnum\addressnum@<\addresscount@\advance\addressnum@\@ne
           \footnote{$^{\hbox{\tiny\number\addressnum@}}$%
           \csname @address\number\addressnum@\endcsname
           \csname @curraddr\number\addressnum@\endcsname
           \csname @email\number\addressnum@\endcsname}\repeat
  \ifx\@empty\@date\else \@footnotetext{\@setdate}\fi
  \ifx\@empty\@subjclass\else \@footnotetext{\@setsubjclass}\fi
  \ifx\@empty\@keywords\else \@footnotetext{\@setkeywords}\fi
  \ifx\@empty\thankses\else \@footnotetext{%
    \def\par{\let\par\@par}\@setthanks}\fi
    }%
  \@setcopyright
}

\def\@tmpevenhead{\relax}

\def\cleardoublepage{\clearpage\if@twoside \ifodd\c@page\else
    \let\@tmpevenhead\@evenhead \let\@evenhead\relax\hbox{}\eject 
    \let\@evenhead\@tmpevenhead\if@twocolumn\hbox{}\newpage\fi\fi\fi}

\def\@setcopyright{%
  \let\copyrightyear\currentyear             
  \insert\copyins{\hsize\textwidth
    \parfillskip\z@ \leftskip\z@\@plus.9\textwidth
    \fontsize{6}{7\p@}\normalfont\upshape
    \everypar{}%
    \vskip-\skip\copyins \nointerlineskip
    \noindent\vrule\@width\z@\@height\skip\copyins
    \copyright\copyrightyear\ American Mathematical Society\par
    \kern\z@}%
}

\renewcommand{\@auth}[2][default]{{\raggedleft
        \begingroup
  \fontsize{\@xivpt}{18}\bfseries
  #2\par \endgroup
        \vspace*{2pc}
        }
        \def\@author{#1}%
}

\renewcommand{\@sauth}[1]{{\raggedleft
        \begingroup
  \fontsize{\@xivpt}{18}\bfseries
  #1\par \endgroup
        \vspace*{2pc}
        }
        \def\@author{#1}%
}

\def\@lecture[#1]#2{%
  \addtocounter{lecture}1\relax
  \refstepcounter{chapter}%
\gdef\thelecturename{#1\unskip}\firstlecturefalse
  {\Large\bfseries
   \raggedleft
   \@xp\uppercase\@xp{\thelecturelabel} {\LARGE\thelecturenum}\\
   \vspace*{3pt}%
    #2\unskip
   \endgraf}%
  \let\@secnumber=\thelecturenum
  \@xp\lecturemark\@xp{\thelecturename}%
  \addcontentsline{toc}{chapter}{%
    \thelecturelabel\ \thelecturenum.\ #2}%
  \vspace{10\p@}\noindent}
  
\def\@slecture#1{%
\iffirstlecture
\gdef\thelecturename{#1\unskip}\firstlecturefalse
  {\Large\bfseries
\noindent\thelecturename
   \endgraf}%
  \let\@secnumber=\thelecturenum
  \@xp\slecturemark\@xp{\thelecturename}%
  \addcontentsline{toc}{chapter}{%
    \thelecturename}%
 \vspace{-6\p@}\noindent
\else
\gdef\thelecturename{#1\unskip}\firstlecturefalse
  {\Large\bfseries
   \raggedleft
   \@xp\uppercase\@xp{\thelecturename}
   \endgraf}%
  \let\@secnumber=\thelecturenum
  \@xp\slecturemark\@xp{\thelecturename}%
  \addcontentsline{toc}{chapter}{%
    \thelecturename}%
  \vspace{10\p@}\noindent
\fi}

\makeatother

\usepackage{amssymb}
\usepackage[all]{xy}
\usepackage{enumerate}
\usepackage{ifpdf}

\ifpdf
\usepackage[pdftex]{graphicx}\pdfcompresslevel=9
\usepackage{hyperref}
\else
\usepackage[dvips]{graphicx}
\fi

\numberwithin{section}{chapter}
\numberwithin{equation}{chapter}

\theoremstyle{plain}
\newtheorem{theorem}[equation]{Theorem}
\newtheorem{lemma}[equation]{Lemma}
\newtheorem{proposition}[equation]{Proposition}
\newtheorem{corollary}[equation]{Corollary}

\theoremstyle{definition}
\newtheorem{definition}[equation]{Definition}
\newtheorem{example}[equation]{Example}
\newtheorem{remark}[equation]{Remark}

\newtheorem{convention}[equation]{Convention}
\newtheorem{exercise}{Exercise}

\newcommand{\G}{\mathcal{G}}
\newcommand{\s}{\mathbf{s}}
\renewcommand{\t}{\mathbf{t}}
\renewcommand{\H}{\mathcal{H}}
\newcommand{\Ss}{\mathbb S}
\newcommand{\Tt}{\mathbb T}
\newcommand{\Rr}{\mathbb R}
\newcommand{\Qq}{\mathbb Q}
\newcommand{\Cc}{\mathbb C}
\newcommand{\Zz}{\mathbb Z}

\newcommand{\set}[1]{\left\{#1\right\}}

\newcommand{\eps}{\varepsilon}

\newcommand{\rmap}{\longrightarrow}
\newcommand{\lmap}{\longleftarrow}

\newcommand{\X}{\ensuremath{\mathfrak{X}}}
\newcommand{\F}{\ensuremath{\mathcal{F}}}
\newcommand{\D}{\ensuremath{\mathcal{D}}}
\renewcommand{\P}{\ensuremath{\mathcal{P}}}
\renewcommand{\O}{\ensuremath{\mathcal{O}}}

\newcommand{\su}{\ensuremath{\mathfrak{su}}}

\newcommand{\NN}{\ensuremath{\mathcal{N}}}

\newcommand{\al}{\alpha}
\newcommand{\be}{\beta} 

\newcommand{\Lie}{\mathcal{L}}  
\renewcommand{\gg}{\mathfrak{g}}
\newcommand{\hh}{\mathfrak{h}}
\renewcommand{\L}{\mathbb L}
\newcommand{\tto}{\rightrightarrows}

\DeclareMathOperator{\ad}{ad}
\DeclareMathOperator{\Ad}{Ad}

\renewcommand{\d}{\mathrm{d}}
\DeclareMathOperator{\Der}{Der}

\DeclareMathOperator{\End}{End}
\DeclareMathOperator{\Exp}{Exp}

\DeclareMathOperator{\GL}{GL}

\DeclareMathOperator{\Graph}{Graph}
\DeclareMathOperator{\Ham}{Ham}
\DeclareMathOperator{\Hol}{Hol}
\DeclareMathOperator{\Hom}{Hom}
\DeclareMathOperator{\Ima}{Im}
\renewcommand{\Im}{\Ima}

\DeclareMathOperator{\Ker}{Ker}

\DeclareMathOperator{\Mor}{Mor}

\DeclareMathOperator{\Rep}{Rep}
\DeclareMathOperator{\rank}{rank}

\DeclareMathOperator{\Tr}{Tr}
\DeclareMathOperator{\var}{var}

\DeclareMathOperator{\DR}{\text{dR}}

\newtheorem{clm}{Claim}

\newcommand{\comment}[1]{}  

\begin{document}

\ifpdf
\DeclareGraphicsExtensions{.pdf,.png}
\pdfinfo
{ /Title      (Lectures on Integrability of Lie Brackets)
  /Author     (M. Crainic and R. L. Fernandes)
  /Subject    (Differential Geometry)
  /Keywords   (Lie groupoids, Lie algebroids, integrability)}
\else
\DeclareGraphicsExtensions{.eps,.epsi}
\fi

\mainmatter
\setcounter{page}{1}

\lectureseries[Integrability]{Lectures on Integrability 
  of Lie Brackets}

\auth[M. Crainic and R.L. Fernandes]{Marius Crainic\\ Rui Loja Fernandes}
\address{Department of Mathematics, Utrecht University, 3508 TA Utrecht, 
The Netherlands}
\email{crainic@math.uu.nl}

\address{Departamento de Matem\'{a}tica, 
Instituto Superior T\'{e}cnico, 1049-001 Lisboa, Portugal} 
\email{rfern@math.ist.utl.pt}


\thanks{This research was supported in part by NWO, FCT/POCTI/FEDER, and by 
	grants POCI/MAT/55958/2004 and POCI/MAT/57888/2004.}
	
\setaddress

\tableofcontents

\setcounter{lecture}{0}
\setcounter{chapter}{0}

%

\lecture*{Preface}               
%

The subject of these lecture notes is the problem of integrating
infinitesimal geometric structures to global geometric structures. 
We follow the categorical approach to differential geometry, where 
the infinitesimal geometric structures are called \emph{Lie algebroids} 
and the global geometric structures are called \emph{Lie groupoids}. 
It is also one of our aims to convince you of the advantages of this approach.

You may not be familiar with the language of Lie algebroids or Lie 
groupoids, but you have already seen several instances of the integrability 
problem. For example, the problem of integrating a vector field to a flow, 
the problem of integrating an involutive distribution to a foliation, 
the problem of integrating a Lie algebra to a Lie group, or the problem of 
integrating a Lie algebra action to a Lie group action. In all these special 
cases, the integrability problem always has a solution. However, in general, this
need not be the case and this new aspect of the problem, the
\emph{obstructions to integrability}, is one of the main topics of these
notes. One such example, that you may have seen before, is the problem of
geometric quantization of a presymplectic manifold, where the so-called
prequantization condition appears.

These notes are made up of five lectures. In the first lecture, we
introduce Lie groupoids and their infinitesimal versions. In the
second lecture, we introduce the abstract notion of a Lie algebroid
and discuss how many notions of differential geometry can be described
in this new language. These first two lectures should help you in
becoming familiar with the language of Lie groupoids and Lie
algebroids. The third and fourth lectures are concerned with
the various aspects of the integrability problem. These two lectures
form the core material of this course, and contain a detailed
description of the integrability obstructions. In the last lecture
we consider, as an example, aspects of integrability related
to Poisson geometry. At the end of each lecture, we have include a
few notes with some references to the literature. We warn
you that these notes are not meant to be complete; they are simply a way to
provide you some historical background as well as further material for you to read and
discover. There are also around 100 exercises which are an integral part of the
lectures: you will need to solve the vast majority of the exercises to
get a good feeling of what integrability is all about!

A version of these lecture notes were used for a course we gave at the
Summer School in Poisson Geometry, held at ICTP, Trieste, in the summer of
2005. Due to a lack of time and space, we were not able to include in
these notes all the topics discussed in the course. Topics left out
include aspects of integrability related to cohomology, quantization
and homotopy theory. We have plans to write a book on the
\emph{Geometry of Lie Brackets}, where all those topics that were left
out (and more!) will be included.

\vskip 1 cm

\hskip 5 cm Marius Crainic

\hskip 5 cm Rui Loja Fernandes

\hskip 5 cm Utrecht and Lisbon, November 2006
\newpage

\lecture{Lie Groupoids}            %
\label{groupoids}                  %
%

\section{Why groupoids?}

These lectures are centered around the notion of a
\emph{groupoid}. What are groupoids? What are they good for? These are
basic questions that we will be addressing and we hope that, by the
end of these lectures, we will have convinced you that groupoids are
worth studying.

It maybe a good idea, even before we start with any formal
definitions, to look at an example. You may wish to keep this
kind of example in mind, since it illustrates very nicely many of the basic
abstract concepts we will be introducing in this lecture.
\vskip 15 pt

Let $N$ be a manifold, and suppose that we want to clasify the set of
Riemannian metrics on $N$ (this is obviously too ambitious, but keep
on reading!). This means that our space of objects, which we will
denote by $M$, is the space of metrics on $N$:
\[ M=\set{g:g\text{ is a Riemannian metric on }N}.\]
This space is quite large. In the classification problem it is natural
not to distinguish two metrics that are related by a
diffeomorphism. So let us consider the triples $(g_1,\phi,g_2)$, where
$g_i$ are metrics on $N$ and $\phi$ is a diffeomorphism of $N$ which
relates the two metrics:
\[ \G=\set{((g_2,\phi,g_1):g_2=\phi_*g_1}.\]

Now $\G$ is precisely the prototype of a groupoid.  The word
``groupoid'' is meant to be suggestive of the notion of group. You
will notice that we cannot always compose two elements of $\G$, but
that sometimes we can do it: if $(g_2,\phi,g_1)$ and $(h_2,\psi,h_1)$
are two elements of $\G$, then we can compose them obtaining the new
element $(g_2,\phi\circ\psi,h_1)$, provided $g_1=h_2$. Also, there are
elements which behave like units under this multiplication, namely
$(g,I,g)$, where $I:N\to N$ is the identity map. Finally, each element
$(g_2,\phi,g_1)$ has an inverse, namely $(g_1,\phi^{-1},g_2)$.

As we have mentioned above, in the classification problem, we identify two
metrics $g_1$ and $g_2$ that differ by a diffeomorphism $\phi$ (i.e.,
such that $g_2=\phi_*g_1$). In other words, we are interested in
understanding the quotient space $M/\G$, which we may call the
\emph{moduli space} of Riemman metrics on $M$. On the other hand, if
we want to pay attention to a fixed metric $g$, then we recognize
immediately that the triples $(g,\phi,g)$ (i.e., the set of
diffeomorphisms preserving this metric) 
is just the group of isometries of $(N,g)$. Hence, our groupoid
encodes all the relevant data of our original problem of studying
and classifying metrics on $N$.
\vskip 15 pt

There is nothing special about this example, where we have chosen to
look at metrics on a manifold. In fact, anytime one tries to study and
classify some class of structures on a space there will be a groupoid
around. 

\section{Groupoids}

Here is the shortest definition of a groupoid: 
\begin{definition} 
  \label{groupoid}
  A \textbf{groupoid} $\G$ is a (small) category in which every arrow is
  invertible. 
\end{definition}

The set of morphisms (arrows) of the groupoid will be denoted by the
same letter $\G$, while the set of objects is denoted by $M_{\G}$, or
even by $M$, provided it is clear from the context what the groupoid
is. We call $M$ the base of the groupoid, and we say that $\G$ is
a groupoid over $M$.  

From its very definition, a groupoid $\G$ over $M$ has certain
underlying \textbf{structure maps}: 
\begin{itemize} 
\item the \textbf{source} and the \textbf{target} maps
  \[ \s, \t: \G\rmap M,\]
  associating to each arrow $g$ its source object $\s(g)$ and its target
  object $\t(g)$. Given $g\in \G$, we write $g: x\rmap y$, or
  $x\stackrel{g}{\rmap} y$, or $y\stackrel{g}{\lmap} x$ to indicate
  that $g$ is an arrow from $x$ to $y$.
\item the \textbf{composition map} 
  \[ m: \G_2\rmap \G ,\]
  is defined on the set $\G_2$ of composable arrows:
  \[ \G_2= \set{(g, h)\in \G\times \G: s(g)= t(h)}.\]
  For a pair $(g,h)$ of composable arrows, $m(g,h)$ is the composition 
  $g\circ h$. We also use the notation $m(g,h)=gh$, and sometimes we
  call $gh$ the multiplication of $g$ and $h$.
\item the \textbf{unit map}
  \[ u: M\rmap \G,\]
  which sends $x\in M_{\G}$ to the identity arrow $1_x\in \G$ at $x$.
  We will often identify $1_x$ and $x$. 
\item the \textbf{inverse map}
  \[ i: \G\rmap \G ,\]
  which sends an arrow $g$ to its inverse $g^{-1}$.
\end{itemize}
Of course, these structure maps also satisfy some identities, similar
to the group case, which again are consequences of Definition \ref{groupoid}:
\begin{itemize}
\item \textbf{law of composition}: if
  $x\stackrel{g}{\lmap}y\stackrel{h}{\lmap} z$, then
  $x\stackrel{gh}{\lmap} z$.
\item \textbf{law of associativity}: if $x\stackrel{g}{\lmap}
  y\stackrel{h}{\lmap} z \stackrel{k}{\lmap} u$, then $g(hk)= (gh)k$.
\item \textbf{law of units}: $x\stackrel{1_x}{\lmap} x$ and, for all 
  $x\stackrel{g}{\lmap}y$, $1_xg=g1_y=g$. 
\item \textbf{law of inverses}: if $x\stackrel{g}{\lmap}y$, then 
  $y\stackrel{g^{-1}}{\lmap}x$ and $gg^{-1}= 1_{y}$, $g^{-1}g= 1_x$.
\end{itemize}
Hence, in a more explicit form, here is the long definition of a
groupoid:

\begin{definition} 
  \label{groupoid:alt}
  A \textbf{groupoid} consists of a set $\G$ (of arrows), a set
  $M_{\G}$ (of objects), and maps $\s, \t, u, m$, and $i$ as above,
  satisfying the laws of composition, associativity, units and
  inverses.
\end{definition}

We will be using the following notation for a groupoid $\G$ over
$M$: if $x\in M$, then the sets 
\[ \G(x, -)=\s^{-1}(x),\ \ \G(-, x)= \t^{-1}(x)\]
are called \textbf{the $s$-fiber at $x$}, and \textbf{the $t$-fiber at $x$},
respectively. The inverse map induces a natural bijection between
these two sets:
\[ i: \G(x, -)\rmap \G(-, x) .\]
Given $g: x\rmap y$, the \textbf{right multiplication by $g$} is only
defined on the $s$-fiber at $y$, and induces a bijection 
\[ R_{g}: \G(y, -)\rmap \G(x, -) .\]
Similarly, the \textbf{left multiplication by $g$} induces a map from
the $t$-fiber at $x$ to the $t$-fiber at $y$.

Next, the intersection of the $s$ and $t$-fiber at $x\in M$,
\[ \G_x=\s^{-1}(x)\cap\t^{-1}(x)=\G(x, -)\cap \G(-, x)\]
together with the restriction of the groupoid multiplication, 
is a group called the \textbf{isotropy group at $x$}.

On the other hand, at the level of the base $M$, one has an
equivalence relation $\sim_{\G}$: two objects $x,y\in M$ are said to be
equivalent if there exists an arrow $g\in \G$ whose source is $x$ and whose
target is $y$. The equivalence class of $x\in M$ is called the
\textbf{orbit through $x$}:
\[ \mathcal{O}_x= \{\t(g): g\in\s^{-1}(x) \},\]
and the quotient set
\[ M/\G:= M/\sim_{\G}= \{ \mathcal{O}_x: x\in M \}\]
is called the \textbf{orbit set of $\G$}. 

The following exercise shows that for a groupoid $\G$ there is still
an underlying group around. 

\begin{exercise}
\label{bisections}
  Let $\G$ be a groupoid over $M$. Define a \textbf{bisection} of
  $\G$ to be a map $b:M\to\G$ such that $\s\circ b$ and $\t\circ
  b$ are bijections. Show that any two bisections can be multiplied,
  so that the set of bisections form a group, denoted $\Gamma(\G)$.
\end{exercise}

The groupoids we will be interested in are not just algebraic
objects. Usually we will be interested in comparing two arrows, looking at
neighborhoods of an arrow, etc. 

\begin{definition} 
  A \textbf{topological groupoid} is a groupoid $\G$ whose set of
  arrows and set of objects are both topological spaces, whose
  structure maps $\s, \t, u, m, i$ are all continuous, and such that
  $\s$ and $\t$ are open maps.
\end{definition}

\begin{example}
  Consider the groupoid $\G$ formed by triples $(g_2,\phi,g_1)$ where
  $g_i$ are metrics on $N$ and $\phi$ is a diffeomorphism taking $g_1$
  to $g_2$. The compact-open topology on the space of diffeomorphisms
  of $N$ and on the space of metrics, induces a natural topology on
  $\G$ so that it becomes a topological groupoid.
\end{example}

Note that for a topological groupoid all the $s$ and $t$- fibers are
topological spaces, the isotropy groups are topological groups, and
the orbit set of $\G$ has an induced quotient topology. 

\begin{exercise} 
  For a topological groupoid $\G$ prove that the unit map 
  \[ u:M\to\G \]
  is a topological embedding, i.e., it is a homeomorphism onto its
  image (furnished with the relative topology).
\end{exercise}

Obviously, one can go one step further and set:

\begin{definition}
  A \textbf{Lie groupoid} is a groupoid $\G$ whose set of arrows and
  set of objects are both manifolds, whose structure maps $\s, \t, u,
  m, i$ are all smooth maps and such that $\s$ and $\t$ are
  submersions.
\end{definition}

Note that the condition that $\s$ and $\t$ are submersions ensure that
the $\s$ and $\t$-fibers are manifolds. They also ensure that the
space $\G_2$ of composable arrows is a submanifold of $\G\times
\G$, and the smoothness of the multiplication map $m$ is to be
understood with respect to the induced smooth structure on $\G_2$.

\begin{convention}
  Unless otherwise stated, all our manifolds are second countable and
  Hausdorff. An exception to this convention is the total space of a
  Lie groupoid $\G$ which is allowed to be non-Hausdorff (to understand 
  why, look at Example \ref{exer:non-hausd}). But we will assume that 
  the base manifold $M$ as well as all the $s$-fibers $\G(x, -)$ (and, 
  hence, the $t$-fibers $\G(-,x)$) are Hausdorff.
\end{convention}

\begin{exercise} 
  Given a Lie groupoid $\G$ over $M$ and $x\in M$, prove that:
  \begin{enumerate}[(a)]
  \item the isotropy groups $\G_x$ are Lie groups;
  \item the orbits $\O_x$ are (regular immersed) submanifolds in
  $M$(\footnote{An immersion $i:N\to M$ is called \emph{regular} if for any
  map $f:P\to N$ the composition $i\circ f:P\to M$ is smooth iff
  $f:P\to N$ is smooth.});
  \item the unit map $u:M\to\G$ is an embedding; 
  \item $\t:\G(x, -)\to \O_x$ is a principal $\G_x$-bundle.
  \end{enumerate}
\end{exercise}

\begin{remark}
\label{bisections-remark}
  Let $\G$ be a Lie groupoid over $M$. Define a smooth bisection to be
  a smooth map $b:M\to\G$ such that $\s\circ b$ and $\t\circ b$ are
  diffeomorphisms. If one furnishes the group $\Gamma(\G)$ of smooth
  bisections with the compact-open topology, one can show that
  $\Gamma(\G)$ is a Fr\'echet Lie group. However, very little is known
  about these infinite dimensional Lie groups, so one prefers to study
  the Lie groupoid $\G$, which is a finite dimensional object.
\end{remark}

Since groupoids are categories, a morphism $\G\to\H$ between two
grou\-poids is a functor: to each arrow and each object in $\G$ we
associate an arrow and an object in $\H$, and these two assignments
have to be compatible with the various structure maps. Later we will
see more general notions of morphisms, but for now this suffices:

\begin{definition} 
  Given a groupoid $\G$ over $M$ and a groupoid $\H$ over $N$, a
  \textbf{morphism} from $\G$ to $\H$, consists of a map $\F: \G\rmap \H$
  between the sets of arrows, and a map $f: M\rmap N$ between the sets
  of objects, which are compatible with all the structure maps. A
  morphism between two topological (respectively, Lie) groupoids is a
  morphism whose components are continuous (respectively, smooth).
\end{definition}

The compatibility of $\F$ with the structure maps translate into the
following explicit conditions:
\begin{itemize}
\item if $g: x\rmap y$ is in $\G$, then $\F(g): f(x)\rmap f(y)$ in $\H$.
\item if $g, h\in \G$ are composable, then $\F(gh)=\F(g)\F(h)$. 
\item if $x\in M$, then $\F(1_x)= 1_{f(x)}$. 
\item if $g: x\rmap y$, then $\F(g^{-1})= \F(g)^{-1}$.
\end{itemize}
Notice that the last property actually follows from the first three.

\begin{exercise}
  If $\F:\G\to\H$ is a morphism of Lie groupoids, show that the
  restrictions $\F|_x:\G_x\to\H_{f(x)}$ are morphisms of Lie groups.
\end{exercise}

If $\G$ is a groupoid, a \textbf{subgroupoid} of $\G$ is a pair
$(\H,i)$, where $\H$ is a groupoid and $i:\H\to\G$ is an injective
groupoid homomorphism. In the case of topological (respectively, Lie)
groupoids we require $i$ to be a topological (respectively, smooth)
immersion. A \textbf{wide subgroupoid} is a subgroupoid $\H\subset\G$
which has the same space of units as $\G$.

The following exercise shows that this notion is more subtle for Lie
groupoids than it is for Lie groups.

\begin{exercise}
  Give an example of a Lie groupoid $\G$ and a wide Lie subgroupoid
  $(\H,i)$, such that $i:\H\to\G$ is an embedding, but its image is not
  a closed submanifold of $\G$.\\
  (Hint: Try the action groupoids, introduced later in this lecture.)
\end{exercise}
\vskip 10 pt

Next we observe that groupoids act on fiber spaces: 

\begin{definition}
  Given a groupoid $\G$ over $M$, a \textbf{$\G$-space} $E$ is defined
  by a map $\mu:E\to M$, called the \textbf{moment map}, together with
  a map
  \[ \G\times_M E=\set{(g,e):\s(g)=\mu(e)}\rmap E,\ (g,e)\mapsto ge,\]
  such that the following action identities are satisfied:
  \begin{enumerate}[(i)]
    \item $\mu(ge)= \t(g)$.
    \item $g(he)=(gh)e$, for all $g,h\in\G$ and $e\in E$ for which it
      makes sense;
    \item $1_{\mu(e)}e=e$, for all $e\in E$.
  \end{enumerate}
\end{definition}
\vskip 10 pt

Note that an action of $\G$ on $E$, with moment map $\mu:E\to M$,
really means that to each arrow $g: x\rmap y$ one associates
an isomorphism 
\[ E_x\to E_y,\ e\mapsto ge,\] 
where $E_x=\mu^{-1}(x)$, such that the action identities are satisfied.

Obviously, if $\G$ is a topological (respectively, a Lie groupoid) and
$E$ is a topological space (respectively, a manifold), then we have
the notion of a continuous (respectively, smooth) action.  

What we have defined, is actually the notion of a \emph{left}
$\G$-space. We leave it to reader the task of defining the notion of a
\emph{right} $\G$-space.

\begin{example}
  A groupoid $\G$ acts on itself by left and right multiplication. For
  the left action, we let $\mu=\t:\G\to M$ be the moment map, while
  for the right action, we let $\mu=\s:\G\to M$. These are continuous
  (respectively, smooth) actions if $\G$ is a topological
  (respectively, Lie) groupoid. We will see more examples of actions
  below.
\end{example}

We mention here one more important notion, which is a simple 
extension of the notion of representations of groups, and which
corresponds to linear groupoid actions.

\begin{definition} 
  Given a (topological, Lie) groupoid $\G$ over $M$, a
  \textbf{representation of $\G$} consists of a vector bundle
  $\mu:E\to M$, together with a (continuous, smooth) linear action of
  $\G$ on $E\to M$: for any arrow of $\G$, $g:x\rmap y$, one has an
  induced linear isomorphism $E_x\to E_y$, denoted $v\mapsto gv$, such
  that the action identities are satisfied.
 \end{definition}

Note that, in particular, each fiber $E_x$ is a representation of the 
isotropy Lie group $\G_x$. 

The isomorphism classes of representations of a groupoid $\G$ form a semi-ring
$(\Rep(\G), \oplus, \otimes)$, where, for two representations $E_1$
and $E_2$ of $\G$, the induced actions on $E_1\oplus E_2$ and
$E_1\otimes E_2$ are the diagonal ones. The unit element is given by
the trivial representation, i.e., the trivial line bundle
$\L_{M}$ over $M$ with $g\cdot 1=1$.

\section{First examples of Lie groupoids}
Let us list some examples of groupoids. We start with two extreme
classes of groupoids which, in some sense, are of opposite nature.

\begin{example}[Groups]
  A group is the same thing as a groupoid for which the set of of
  objects contains a single element. So groups are very particular
  instances of groupoids. Obviously, topological (respectively, Lie)
  groups are examples of topological (respectively, Lie) groupoids.
\end{example}

\begin{example}[The pair groupoid]
  At the other extreme, let $M$ be any set. The Cartesian product
  $M\times M$ is a groupoid over $M$ if we think of a pair $(y,x)$ as
  an arrow $x\rmap y$. Composition is defined by:
  \[ (z,y)(y,x)=(z,x).\]
  If $M$ is a topological space (respectively, a manifold), then the
  pair groupoid is a topological (respectively, Lie) groupoid.

  Note that each representation of the pair groupoid is isomorphic
  to a trivial vector bundle with the tautological action. Hence
  \[ \Rep(M\times M)= \mathbb{N},\]
  the semi-ring of non-negative integers.
\end{example}

Our next examples are genuine examples of groupoids, which already
exhibit some distinct features.

\begin{example}[General linear groupoids]
  If $E$ is a vector bundle over $M$, there is an associated general
  linear groupoid, denoted by $GL(E)$, which is similar to the general
  linear group $GL(V)$ associated to a vector space $V$. $GL(E)$ is a
  groupoid over $M$ whose arrows between two points $x$ and $y$
  consist of linear isomorphisms $E_x\rmap E_y$, and the
  multiplication of arrows is given by the composition of maps.
  
  Note that, given a general Lie groupoid $\G$ and a vector bundle $E$
  over $M$, a linear action of $\G$ on $E$ (making $E$ into a
  representation) is the same thing as a homomorphism of Lie groupoids
  $\G\to GL(E)$.
\end{example}

\begin{example}[The action groupoid] 
  Let $G\times M\to M$ be an action of a group on a set $M$. The
  corresponding action groupoid over $M$, denoted $G\ltimes M$, has as
  space of arrows the Cartesian product:
  \[\G=G\times M.\] 
  For an arrow $(g,x)\in\G$ its source and target are given by:
  \[ \s(g,x)=x,\quad \t(g,x)=g\cdot x,\]
  while the composition of two arrows is given by:
  \[ (h,y)(g,x)=(hg,x).\]
  Notice that the orbits and isotropy groups of $\G$ coincide with the
  usual notions of orbits and isotropy groups of the action $G\times
  M\to M$. Also,
  \[ \Rep(G\ltimes M)=\text{Vect}_{G}(M),\]
  the semi-ring of equivariant vector bundles over $M$.
\end{example}

\begin{example}[The gauge groupoid] 
\label{gauge-gpd}
  If $G$ is a Lie group and $P\rmap M$ is a principal $G$-bundle, then
  the quotient of the pair groupoid $P\times P$ by the (diagonal)
  action of $G$ is a groupoid over $P/G= M$, denoted $P\otimes_{G}P$,
  and called the gauge groupoid of $P$. Note that all isotropy groups
  are isomorphic to $G$. Also, fixing a point $x\in M$, a
  representation $E$ of $P\otimes_{G}P$ is uniquely determined by
  $E_x\in \Rep(G)$, and this defines an isomorphism of semi-rings:
  \[ \Rep(P\otimes_{G}P)\sim \Rep(G).\]
  A particular case is when $G= GL_n$ and $P$ is the frame bundle of a
  vector bundle $E$ over $M$ of rank $n$. Then the resulting gauge
  groupoid coincides with $GL(E)$ mentioned above.
\end{example}

\begin{example}[The fundamental groupoid of a manifold]
  Closely related to the pair groupoid of a manifold $M$ is the
  fundamental groupoid of $M$, denoted $\Pi_1(M)$, which consists of
  homotopy classes of paths with fixed end points (we assume that $M$
  is connected). In the light of Example \ref{gauge-gpd}, this is the
  gauge groupoid associated to the universal cover of $M$, viewed as a
  principal bundle over $M$ with structural group the fundamental
  group of $M$- and this point of view also provides us with a smooth 
  structure on $\Pi_1(M)$, making it into a Lie groupoid.
  Note also that, when $M$ is simply-connected, $\Pi_1(M)$
  is isomorphic to the pair groupoid $M\times M$ (a homotopy class of
  a path is determined by its end points). In general, there is an
  obvious homomorphism of Lie groupoids $\Pi_1(M)\rmap M\times M$,
  which is a local diffeomorphism.
\end{example}

\begin{exercise} 
  Show that representations of $\Pi_1(M)$ correspond to vector bundles
  over $M$ endowed with a flat connection.
\end{exercise}

\begin{example}[The fundamental groupoid of a foliation] 
\label{fund-gpd}
  More generally, let $\F$ be a foliation of a space $M$ of class
  $C^r$ ($0\le r\le\infty$). The fundamental groupoid of $\F$, denoted
  by $\Pi_1(\F)$, consists of the leafwise homotopy classes of paths
  (relative to the end points):
  \[ 
  \Pi_1(\F)=\set{[\gamma]:\gamma:[0,1]\to M\text{ a path lying in a leaf}}.
  \]
  For an arrow $[\gamma]\in\Pi_1(\F)$ its source and target are given by:
  \[ \s([\gamma])=\gamma(0),\quad \t([\gamma])=\gamma(1),\]
  while the composition of two arrows is just concatenation of
  paths:
  \[ [\gamma_1][\gamma_2]=[\gamma_1\cdot\gamma_2]. \]
  In this example, the orbit $O_x$ coincides with the leaf $L$ through
  $x$, while the isotropy group $\Pi_1(\F)_x$ coincides with the
  fundamental group $\pi_1(L,x)$.
\end{example}

\begin{exercise}
  Verify that if $\F$ is a foliation of class $C^r$
  ($0\le r\le\infty$) then $\Pi_1(\F)$ is a groupoid of class $C^r$.
\end{exercise}

The next exercise gives a simple example of a non-Hausdorff Lie groupoid.

\begin{exercise}\label{exer:non-hausd}
  Let $\F$ be the smooth foliation of $\Rr^3-\{0\}$ by horizontal
  planes. Check that the Lie groupoid $\Pi_1(\F)$ is not Hausdorff.
\end{exercise}

\section{The Lie algebroid of a Lie groupoid}

Motivated by what we know about Lie groups, it is natural to wonder
what are the infinitesimal objects that correspond to Lie groupoids,
and these ought to be called \emph{Lie algebroids}.

Let us recall how we construct a Lie algebra out of a Lie group: if
$G$ is a Lie group, then its Lie algebra $\text{Lie}(G)$ consists of:
\begin{itemize} 
\item an underlying vector space, which is just the tangent space to
$G$ at the unit element.
\item A Lie bracket on $\text{Lie}(G)$, which comes from the
identification of $\text{Lie}(G)$ with the space $\X_{\text{inv}}(G)$
of right invariant vector fields on $G$, together with the fact that
the space of right invariant vector fields is closed under the usual
Lie bracket of vector fields: 
\[ [\X_{\text{inv}}(G),\X_{\text{inv}}(G)]\subset\X_{\text{inv}}(G).\]
\end{itemize}
Now, back to general Lie groupoids. One remarks two novelties when
comparing with the Lie group case. First of all, there is more than
one unit element. In fact, there is one unit for each point in $M$,
hence we expect a vector bundle over $M$, instead of a vector
space. Secondly, the right multiplication by elements of $\G$ is only
defined on the $s$-fibers. Hence, to talk about right invariant vector
fields on $\G$, we have to restrict attention to those vector fields
which are tangent to the $s$-fibers, i.e., to the sections of the 
sub-bundle $T^{s}\G$ of $T\G$ defined by
\[ T^{s}\G=\Ker(\d\s) \subset T\G .\]
Having in mind the discussion above, we first define $\text{Lie}(\G)$
as a vector bundle over $M$.

\begin{definition} 
  Given a Lie groupoid $\G$ over $M$, we define the vector bundle
  $A=\text{Lie}(\G)$ whose fiber at $x\in M$ coincides with the
  tangent space at the unit $1_x$ of the $s$-fiber at $x$. In other
  words:
  \[ A: = T^{s}\G|_{M}.\]
\end{definition}

Now we describe the Lie bracket of the Lie algebroid $A$. This is, in
fact, a bracket on the space of sections $\Gamma(A)$, and to deduce it
we will identify the space $\Gamma(A)$ with the space of right
invariant vector fields on $\G$.  {To} see this, we observe that the
fiber of $T^{s}(\G)$ at an arrow $h: y\rmap z$ is
\[ T^{s}_{h}\G= T_{h}\G(y, -),\]
so, for any arrow $g: x\rmap y$, the differential of the right
multiplication by $g$ induces a map
\[ R_g: T^{s}_{h}\G\to T^{s}_{hg}\G.\]
Hence, we can describe the space of right invariant vector fields on
$\G$ as:
\[ \X^{s}_{\text{inv}}(\G)=\set{X\in\Gamma(T^{s}\G): X_{hg}=
  R_g(X_h),\ \forall\ (h,g)\in \G_2}.\]
Now, given $\alpha\in \Gamma(A)$, the formula
\[ \tilde{\alpha}_{g}= R_{g}(\alpha_{\t(g)}) \]
clearly defines a right invariant vector field. Conversely, any vector
field $X\in\X^{s}_{\text{inv}}(\G)$ arises in this way: the invariance of $X$
shows that $X$ is determined by its values at the points in $M$:
\[ X_{g}= R_g(X_y),\ \ \ \text{for\ all}\ g: x\rmap y,\]
i.e., $X=\tilde{\alpha}$ where $\alpha:= X|_{M}\in\Gamma(A)$. 
Hence, we have shown that there exists an isomorphism
\begin{equation}
\label{righ-inv} 
\Gamma(A)\stackrel{\sim}{\rmap} \X^{s}_{\text{inv}}(\G),\ \ \ 
\alpha\mapsto \tilde{\alpha}.
\end{equation}
On the other hand, the space $\X^{s}_{\text{inv}}(\G)$ is a Lie
subalgebra of the Lie algebra $\X(\G)$ of vector fields on $\G$ with
respect to the usual Lie bracket of vector fields. This is clear
since the pull-back of vector fields on the $\s$-fibers, along
$R_{g}$, preserves brackets.



\begin{definition} 
  The \textbf{Lie bracket} on $A$ is the Lie bracket on $\Gamma(A)$
  obtained from the Lie bracket on $\X_{\text{inv}}(G)$ under the
  isomorphism (\ref{righ-inv}).
\end{definition}

Hence this new bracket on $\Gamma(A)$, which we denote by $[\cdot,
\cdot]_{A}$ (or simply $[\cdot, \cdot]$, when there is no danger of
confusion) is uniquely determined by the formula:
\begin{equation}  
\label{the-bracket}
\widetilde{[\alpha, \beta]}_{A}= [\tilde{\alpha}, \tilde{\beta}] .
\end{equation}

To describe the entire structure underlying $A$, we need one more 
piece. 

\begin{definition} 
  The \textbf{anchor map} of $A$ is the bundle map
  \[ \rho_{A}: A\to TM \] 
  obtained by restricting $\d\t: T\G\to TM$ to $A\subset T\G$.
\end{definition}

Again, when there is no danger of confusion (e.g., in the discussion
below), we simply write $\rho$ instead of $\rho_{A}$.  The next
proposition shows that the bracket and the anchor are related by a
Leibniz-type identity. We use the notation $\Lie_{X}$ for the Lie
derivative along a vector field.

\begin{proposition} 
  For all $\alpha, \beta\in \Gamma(A)$ and all $f\in C^{\infty}(M)$,
  \[ [\alpha, f\beta]= f[\alpha, \beta]+
  \Lie_{\rho(\alpha)}(f)\beta.\]
\end{proposition}

\begin{proof} 
  We use the characterization of the bracket (\ref{the-bracket}).
  First note that $\widetilde{f\beta}= (f\circ\t)\tilde{\beta}$. From
  the usual Leibniz identity of vector fields (on $\G$),
  \begin{align*}
    \widetilde{[\al,f\beta]}&=[\tilde{\alpha},\widetilde{f\beta}]\\
    &=[\tilde{\alpha}, (f\circ\t)\tilde{\beta}]\\
    &=(f\circ\t)[\tilde{\alpha}, \tilde{\beta}]+ 
    \Lie_{\tilde{\alpha}}(f\circ\t)\tilde{\beta}.
  \end{align*} 
  But, at a point $g: x\rmap y$ in $\G$, we find: 
  \begin{align*}
    \Lie_{\tilde{\alpha}}(f\circ\t)(g)&=\d_{g}(f\circ\t)(\tilde{\alpha}(g))\\
    &=\d_{\t(g)}f\circ \d_g\t(\tilde{\alpha}(g))
    =\Lie_{\rho(\al)}(f)(\t(g)).
  \end{align*}
  Hence, we conclude that:
  \[  
  \widetilde{[\al,f\beta]}=\widetilde{f[\al,\be]}
     +\widetilde{\Lie_{\rho(\al)}(f)\be},
  \]
  and the result follows.
\end{proof}

We summarize this discussion in the following definition:

\begin{definition} 
  The \textbf{Lie algebroid of the Lie groupoid $\G$} is the vector
  bundle $A=\text{Lie}(\G)$, together with the anchor $\rho_{A}:A\to
  TM$ and the Lie bracket $[\cdot, \cdot]_{A}$ on $\Gamma(A)$. 
\end{definition}

\begin{exercise}
  For each example of a Lie groupoid furnished in the paragraph above,
  determine its Lie algebroid. 
\end{exercise}

{To} complete our analogy with the Lie algebra of a Lie group,
and introduce the exponential map of Lie groupoid, 
we set:

\begin{definition}
  \label{flow-al}
  For $x\in M$, we put
  \[ \phi_{\alpha}^{t}(x):= \phi_{\tilde{\alpha}}^{t}(1_x)\in \G\]
  where $\phi_{\tilde{\alpha}}^{t}$ is the flow of the right invariant
  vector field $\tilde{\al}$ induced by $\al$. We call
  $\phi_{\al}^{t}$ \textbf{the flow} of $\al$.
\end{definition}

\begin{remark} 
  \label{rem:exponential}
  The relevance of the flow comes from the fact that it provides the
  bridge between sections of $A$ (infinitesimal data) and elements of
  $\G$ (global data). In fact, building on the analogy with the
  exponential map of a Lie group and Remark \ref{bisections-remark},
  we see that the flow of sections can be interpreted as defining an
  \textbf{exponential map} $\exp:\Gamma(A)\to\Gamma(\G)$ to the group
  of bisections of $\G$:
  \[ \exp(\al)(x)=\phi_{\al}^{1}(x).\]
  This is defined provided $\al$ behaves well enough (e.g., if it has
  compact support). This shows that, in some sense, $\Gamma(A)$ is the
  Lie algebra of $\Gamma(\G)$.
\end{remark}

Using the exponential notation, we can write the flow as
$\phi_{\al}^{t}=\exp(t\al)$.
The next exercise will show that, just as for Lie groups,
$\exp(t\al)$ contains all the information needed to recover the
entire flow of $\tilde{\al}$.

\begin{exercise} 
  For $\alpha\in \Gamma(A)$, show that:
  \begin{enumerate}[(a)]
  \item For all $y\in M$, $\phi_{\alpha}^{t}(y): y\rmap
    \phi_{\rho(\alpha)}^{t}(y)$. 
  \item For all $g: x\rmap y$ in $\G$, $\phi_{\tilde{\alpha}}^{t}(g):
    x\rmap \phi_{\rho(\alpha)}^{t}(y)$, and 
    $\phi_{\tilde{\alpha}}^{t}(g)= \phi_{\alpha}^{t}(y)g$.
  \end{enumerate}
  (here, $\phi_{\rho(\alpha)}^{t}$ is the flow of the vector field
  $\rho(\alpha)$ on $M$). 
\end{exercise}

\section{Particular classes of groupoids}

There are several particular classes of groupoids that deserve special
attention, due to their relevance both in the theory and in various
applications.

First of all, recall that a topological space $X$ is called
$k$-connected (where $k\geq 0$ is an integer) if $\pi_{i}(X)$ is
trivial for all  $0\le i\le k$ and all base points.

\begin{definition}  
  A topological groupoid $\G$ over a space $M$ is called \textbf{source
  $k$-connected} if the $s$-fibers $\s^{-1}(x)$ are $k$-connected for
  every $x\in M$. When $k=0$ we say that $\G$ is a \textbf{$s$-connected
  groupoid}, and when $k=1$ we say that $\G$ is a \textbf{$s$-simply
  connected groupoid}.
\end{definition}

\begin{exercise}
  Show that a source $n$-connected Lie groupoid $\G$ over a
  $n$-connected base $M$, has a space of arrows which is
  $n$-connected.\\ 
  (Hint: Recall that for a Lie groupoid the source map $\s:\G\to M$ is
  a submersion.)
\end{exercise}

Most of the groupoids that we will meet in these lectures are
source-connected. Any Lie groupoid $\G$ has an associated
source-connected groupoid $\G^0$- the $s$-connected component of
the identities. More precisely, $\G^0\subset \G$ consists of those
arrows $g:x\rmap y$ of $\G$ which are in the connected component of
$G(x, -)$ containing $1_x$. 

\begin{proposition}
  \label{prop:0:connected}
  For any Lie groupoid, $\G^0$ is an open subgroupoid of
  $\G$. Hence, $\G^0$ is a $s$-connected Lie groupoid that has the same Lie
  algebroid as $\G$.
\end{proposition}

\begin{proof} 
  Note that right multiplication by an arrow $g:x\rmap y$ is a
  homeomorphism from $\s^{-1}(y)$ to $\s^{-1}(x)$. Therefore, it maps
  connected components to connected components. If $g$ belongs to the
  connected component of $\s^{-1}(x)$ containing $1_x$, then right
  multiplication by $g$ maps $1_y$ to $g$, so it maps the connected
  component of $1_y$ to the connected component of $1_x$. Hence $\G^0$
  is closed under multiplication. Moreover, $g^{-1}$ is mapped to
  $1_x$, so that $g^{-1}$ belongs to the connected component of
  $1_y$, and hence $\G^0$ is closed under inversion. This shows that
  $\G^0$ is a subgroupoid.

  {To} check that $\G^0$ is open, it is enough to observe that there
  exists an open neighborhood $U$ of the identity section $M\subset\G$
  which is contained in $\G^0$. To see that, observe that by the local
  normal form of submersions, each $1_x\in\G$ has an open neighborhood
  $U_x$ which intersects each $s$-fiber in a connected set, i.e.,
  $U_x\subset \G^0$. Then $U=\cup_{x\in M} U_x$ is the desired
  neighborhood.

\end{proof}

\begin{exercise}
  Is this proposition still true for topological groupoids? 
\end{exercise}
\vskip 10 pt

Analogous to simply connected Lie groups and their role in standard Lie theory, groupoids which are $s$-simply connected play a fundamental role in 
the Lie theory of Lie groupoids. 

\begin{theorem}
  \label{thm:1:connected}
  Let $\G$ be a $s$-connected Lie groupoid. There exist a Lie
  groupoid $\tilde{\G}$ and homomorphism $\F:\tilde{\G}\to\G$ such
  that:
  \begin{enumerate}[(i)]
    \item $\tilde{\G}$ is $s$-simply connected.
    \item $\tilde{\G}$ and $\G$ have the same Lie algebroid.
    \item $\F$ is a local diffeomorphism.
  \end{enumerate}
  Moreover, $\tilde{\G}$ is unique up to isomorphism.
\end{theorem}

\begin{proof}
  We define the new groupoid $\tilde{\G}$ by letting $\tilde{\G}(x,
  -)$ be the universal cover of $\G(x, -)$ consisting of homotopy
  classes (with fixed end points) of paths starting at $1_x$. Given
  $g:x\rmap y$, representing a homotopy class $[g]\in\tilde{\G}$, we
  set $\s([g])= x$ and $\t([g])=y$.  Given $[g_1], [g_2]\in
  \tilde{\G}$ composable, we define $[g_1]\cdot [g_2]$ as the homotopy
  class of the concatenation of $g_{2}$ with $R_{g_2(1)}\circ g_1$
  (draw a picture!). 

  It is easy to check that $\tilde{\G}$ is a groupoid. To describe its
  smooth structure, one remarks that $\tilde{\G}=p^{-1}(M)$, where
  $p:\Pi_1(\F(\s))\to \G$ is the source map of the fundamental groupoid
  of the foliation on $\G$ by the fibers of $\s$. Since $p$ is a
  submersion, $\tilde{\G}$ will be smooth and of the same dimension as
  $\G$. It is not difficult to check that the structure maps are also
  smooth. Moreover, the projection $\tilde{\G}\to \G$ which sends $[g]$
  to $g(1)$ is clearly a surjective groupoid morphism ($\G$ is
  $\s$-connected), which is easily seen to be a local
  diffeomorphism. Hence, it induces an isomorphism at the level of
  algebroids. We leave the proof of uniqueness as an exercise.
\end{proof}

At the opposite extreme of $s$-connectedness there are the groupoids that
model the leaf spaces of foliations, known as \'etale groupoids.

\begin{definition} 
  A Lie groupoid $\G$ is called \textbf{\'etale} if its source map
  $\s$ is a local diffeomorphism.
\end{definition}

\begin{example} Here are a few simple examples of \'etale groupoids:
  \begin{enumerate}
  \item The fundamental groupoid $\Pi_1(M)$ of any manifold $M$ is 
  	always \'etale.
  \item The fundamental groupoid of a foliation is not \'etale in
    general. However, let $T$ be a complete transversal, i.e., an 
    immersed submanifold which intersects every leaf, and is transversal 
    at each intersection point (note that $T$ does not have to be connected). 
    Then, each holonomy transformation determines an arrow between points of 
    $T$, and so they form a groupoid $\Hol(T)\tto T$, which is \'etale (\footnote{
    The resulting groupoid is equivalent to the fundamental 
    groupoid of the foliation, in a certain sense that can be made precise.}). 
  \item An action groupoid $\G=G\ltimes M$ is \'etale if and only if $G$ is
    discrete. To see this, just notice that $\s:\G\to M$ is the
    projection $\s(g,m)=m$, so it is a local diffeomorphism iff $G$ is
    discrete. 
\end{enumerate}
\end{example}

\begin{exercise}
  Show that a Lie groupoid is \'etale iff its source fibers are discrete.
\end{exercise}

Analogous to compact Lie groups are proper Lie groupoids, which define
another important class.

\begin{definition} 
  A Lie groupoid $\G$ over $M$ is called \textbf{proper} if the map
  $(\s, \t): \G\to M\times M$ is a proper map.  
\end{definition}

We will see that, for proper groupoids, any representation admits an
invariant metric. The following exercise asks you to prove some other
basic properties of proper Lie groupoids.

\begin{exercise}
  Let $\G$ be a proper Lie groupoid over $M$. Show that:
  \begin{enumerate}[(a)]
    \item All isotropy groups of $\G$ are compact.
    \item All orbits of $\G$ are closed submanifolds.
    \item The orbit space $M/\G$ is Hausdorff.
  \end{enumerate}
\end{exercise}

\begin{example} 
  Let us list some examples of proper groupoids:
  \begin{enumerate}
  \item For a manifold $M$, the pair groupoid $M\times M$ is always
    proper. On the other hand, the fundamental groupoid $\Pi_1(M)$ is
    proper iff the fundamental groups of $M$ are finite.
  \item A Lie group $G$ is proper (as a groupoid) if and only if it is
    compact. 
  \item An action Lie groupoid $G\ltimes M$ is proper if and only if the
    action of $G$ on $M$ is a proper action. Actually, this is just a
    matter of definitions, since a proper action is usually defined as
    one for which the map $G\times M\to M\times M$, $(g,m)\mapsto
    (m,gm)$ is a proper map.
  \end{enumerate}
\end{example}

\begin{exercise}
  Check that groupoids of type $\GL(E)$ are never proper. Choose a
  metric on the vector bundle $E$ and define $O(E)\subset\GL(E)$ to be
  the subgroupoid of isometries of the fibers. Show that $O(E)$ is
  proper.
\end{exercise}

We will say that a Lie groupoid $\G$ is \textbf{source locally
trivial} if the source map $\s:\G\to M$ is a locally trivial
fibration. This implies that the target map is also a locally trivial
fibration. Note that there are examples of proper groupoids for
which the source (or target) map is not a locally trivial fibration. 

\begin{exercise}
  Consider the foliation $\F$ of $\Rr^2-\{0\}$ given by horizontal
  lines. Show that $\Pi_1(\F)$ is neither proper nor source locally trivial.
\end{exercise}

\begin{exercise}
  Consider the foliation $\F$ of $\Rr^3-\{(x,0,0):x\in[0,1]\}$ given by spheres
  around the origin. Show that $\Pi_1(\F)$ is proper but not source locally trivial.
\end{exercise}

\begin{exercise}
  Show that a Lie groupoid which is both proper and \'etale is locally
  trivial. 
\end{exercise}

Groupoids that are both proper and \'etale form a very important class
of groupoids, since they serve as models for orbifolds. 

\section{Notes}

According to Weinstein, the notion of a groupoid was discovered by
Brandt \cite{Bra} in the early twenty century, while studying
quadratic forms over the integers. Groupoids where introduced into
differential geometry by Ehresmann in the 1950's, and he also considered
more general ``structured categories'' (see the comments on his work
in \cite{Ehr}). Already in Ehresmann's work one can find applications
to foliations, fibered spaces, geometry of p.d.e.'s, etc. More or less
at the same time, Grothendieck \cite{Gro} advocated the use of
groupoids in algebraic geometry as the right notion to understand
moduli spaces (in the spirit of the introductory section in this lecture),
and from that the theory of stacks emerged \cite{BFFGK}. 
Important sources of examples, which strongly influenced the theory
of Lie groupoids, comes from Haefliger's approach to transversal geometry of foliations, from Connes' noncommutative geometry
\cite{Connes} and from Poisson geometry (see the announcement \cite{Wein1} and the 
first systematic exposition of symplectic groupoids in \cite{CoDaWe}), with
independent contributions from others (notably Karas\"ev
\cite{Kar} and Zakrzewski \cite{Zak})). 

The book of MacKenzie \cite{Mack1} contains a nice introduction to the 
theory of transitive and locally trivial groupoids, which has now been 
superseeded by his new book \cite{Mack3}, where he also treats non-transitive 
groupoids and double structures. The modern approach to groupoids can
also be found in recent monographs, such as the book by Cannas da Silva and
Weinstein \cite{CaWe} devoted to the geometry of noncommutative
objects, the book by Moerdijk and Mr\v{c}un \cite{MoMr} on Lie
groupoids and foliations, and the book by Dufour and Zung \cite{DuZu}
on Poisson geometry.

Several important aspects of Lie groupoid theory have not made it yet
to expository books. Let us mentioned as examples the theory of proper
Lie groupoids (see, e.g., the papers by Crainic \cite{Cr}, Weinstein
\cite{Wein5} and Zung \cite{Zun}), the theory of differentiable
stacks, gerbes and non-abelian cohomology (see, e.g., the preprints by
Behrend and Xu \cite{BeXu} and Moerdijk \cite{Moer2}).

As time progresses, there are also a few aspects of Lie groupoid
theory that seem to have lost interest or simply disapeared. We feel
that some of these are worth recovering and we mentioned two
examples. A prime example is provided by Haefliger's approach to
integrability and homotopy theory (see the beatifull article
\cite{Haef}). Another remarkable example is the geometric approach to
the theory of p.d.e.'s sketched by Ehresmann's, with roots in
E.~Cartan's pioneer works, and which essentially came to an halt with
the monograph by Kumpera and Spencer \cite{KuSp}.
 

\lecture{Lie Algebroids}           %
\label{algebroids}                 %
%

\section{Why algebroids?}
We saw in Lecture \ref{groupoids} that Lie groupoids are natural
objects to study in differential geometry. Moreover, to every Lie
groupoid $\G$ over $M$ there is associated a certain vector bundle
$A\to M$, which carries additional structure, namely a Lie bracket on
the sections and a bundle map $A\to TM$. If one axiomatizes these
properties, one obtains the abstract notion of a Lie algebroid:

\begin{definition} 
  A \textbf{Lie algebroid} over a manifold $M$ consists of a vector
  bundle $A$ together with a bundle map $\rho_A:A\to TM$ and a Lie
  bracket $[~,~]_A$ on the space of sections $\Gamma(A)$, satisfying
  the Leibniz identity
  \[ [\alpha, f\beta]_A= f[\alpha, \beta]_A+\Lie_{\rho_A(\alpha)}(f)\beta,\]
  for all $\alpha, \beta\in \Gamma(A)$ and all $f\in C^{\infty}(M)$.
\end{definition} 

\begin{exercise} 
  Prove that the induced map $\rho_A:\Gamma(A)\to \X(M)$ is
  a Lie algebra homomorphism.
\end{exercise}

Before we plug into the study of Lie algebroids, we would like to show
that Lie algebroids are interesting by themselves, independently of
Lie grou\-poids. In fact, in geometry one is led naturally to study
Lie algebroids and, in their study, it is useful to integrate them to
Lie groupoids!

In order to illustrate this point, just like we did at the begining
of Lecture \ref{groupoids}, we look again at equivalence problems in
geometry. \'Elie Cartan observed that many equivalence problems in
geometry can be best formulated in terms of coframes. Working out the
coframe formulation, he was able to come up with a method, now called
Cartan's equivalence method, to deal with such problems.

A local version of Cartan's formulation of equivalence problems can be
described as follows: one is given a family of functions
$c^{i}_{j,k}$, $b^{a}_{i}$ defined on some open set $U\subset
\mathbb{R}^n$, where the indices satisfy $1\leq i, j, k\leq r$, $1\leq
a\leq n$ ($n$, $r$ positive integers). Then the problem is:

\textbf{Cartan's problem:} \emph{find a manifold $N$, a coframe $\{\eta^i\}$
on $N$, and a function $h: N\to U$, satisfying the equations:
\begin{align}
  \label{structure:eq:1}
  \d\eta^i&=\sum_{j, k} c^{i}_{j,k}(h) \eta^j\wedge \eta^k, \\
  \label{structure:eq:2}
  \d h^a  &=\sum_{i} b^{a}_{i}(h) \eta^i.
\end{align}
As part of this problem we should be able to answer the following
questions:
\begin{itemize}
\item When does Cartan's problem have a solution?
\item What are the possible solutions to Cartan's problem?
\end{itemize}}

Here is a simple, but interesting, illustrative example:

\begin{example}
  Let us consider the problem of equivalence of metrics in $\Rr^2$
  with constant Gaussian curvature. 
  If $\d s^2$ is a metric in $\Rr^2$ then there exists a diagonalizing
  coframe $\{\eta^1,\eta^2\}$, so that:
  \[ \d s^2=(\eta^1)^2+(\eta^2)^2.\]
  In terms of Cartesian coordinates, this coframe can be written as:
  \[ \eta^1=A\d x+B\d y,\quad \eta^2=C\d x+D\d y,\]
  where $A, B, C, D$ are smooth functions of $(x,y)$ satisfying
  $AD-BC\not=0$. Taking exterior derivatives, we obtain:
  \begin{equation} 
    \label{metric:str:1}
    \d\eta^1=J\eta^1\wedge\eta^2,\quad \d\eta^2=K\eta^1\wedge\eta^2,
  \end{equation}
  where:
  \[ J=\frac{B_x-A_y}{AD-BC},\quad K=\frac{D_x-C_y}{AD-BC}.\]

  If $f$ is any smooth function in $\Rr^2$, the coframe derivatives of
  $f$ are defined to be the coefficients of the differential $\d f$ when
  expressed in terms of the coframe:
  \[ \d f=\frac{\partial f}{\partial \eta^1}\eta^1+
  \frac{\partial f}{\partial \eta^2}\eta^2.\]
  For example, the Gaussian curvature of $\d s^2$ is given in terms of
  the coframe derivatives of the structure functions by:
  \[ \kappa=
  \frac{\partial J}{\partial \eta^2}-\frac{\partial K}{\partial\eta^1}
  -J^2-K^2
  \]
  So we see that for metrics of constant Gaussian curvature, the
  two structure functions $J$ and $K$ are not independent, and we can
  choose one of them as the independent one, say $J$. Then we can
  write:
  \begin{equation} 
    \label{metric:str:2}
    \d J=L\eta^1\wedge\eta^2.
  \end{equation}
  Equations (\ref{metric:str:1}) and (\ref{metric:str:2}) form the
  structure equations of the Cartan problem of classifying metrics of
  constant Gaussian curvature. Using these equations one can find
  normal forms for all metrics of constant Gaussian curvature.

\begin{exercise}
  Show that the metrics of zero Gaussian curvature ($\kappa=0$) can be
  reduced to the following form:
  \[  \d s^2=\d x^2+\d y^2. \]
\end{exercise}

\end{example}
\vskip 20 pt

Obvious \emph{necessary conditions} to solve Cartan's problem can be
obtained as immediate consequences of the fact that $\d^2= 0$ and that 
$\{\eta^i\}$ is a coframe. In fact, a simple computation
gives(\footnote{We will be using the Einstein's summation convention
  without further notice.}):
\begin{equation}
  \label{Cartan:eq:1}
  F^b_i(h){\frac{\partial F^a_j}{\partial h^b}}(h)
  -F^b_j(h){\frac{\partial F^a_i}{\partial h^b}}(h)
  = - c^l_{i,j}(h)\,F^a_l(h),
\end{equation}
\begin{multline}
  \label{Cartan:eq:2}
  F^a_j(h){\frac{\partial c^i_{k,l}}{\partial h^a}}(h)
  +F^a_k(h){\frac{\partial c^i_{l,j}}{\partial h^a}}(h)
  +F^a_l(h){\frac{\partial c^i_{j,k}}{\partial h^a}}(h)=\\
  =-\bigl(c^i_{m,j}(h)c^m_{k,l}(h)
  +c^i_{mk}(h)c^m_{l,j}(h)+c^i_{m,l}(h)c^m_{j,k}(h)\bigr).
\end{multline}

\begin{exercise}
  Take exterior derivatives of the structure equations
  (\ref{structure:eq:1}) and (\ref{structure:eq:2}), and deduce the
  relations (\ref{Cartan:eq:1}) and (\ref{Cartan:eq:2}).
\end{exercise}

Cartan's problem is extremely relevant and suggestive already in the
case where the structure functions $c^i_{j,k}$ are constants (hence
no $h$'s appear in the problem). In this case, condition
(\ref{Cartan:eq:1}) is vacuous, while condition (\ref{Cartan:eq:2}) is
the usual Jacobi identity. Therefore, these conditions precisely mean
that the $c^{i}_{j,k}$'s define an ($r$-dimensional) Lie algebra. Let
us call this Lie algebra $\mathfrak{g}$, so that on a preferred basis
$\{e_i\}$ we have:
\begin{equation}
  \label{bracket} 
  [e_{j}, e_{k}]= c^{i}_{j,k} e_i
\end{equation}
Moreover, a family $\{\eta^i\}$ of one-forms on $N$ can be viewed as a
$\mathfrak{g}$-valued form on $N$, $\eta=\eta^ie_i\in
\Omega^1(N; \mathfrak{g})$, and the equations take the following
global form:
\begin{equation}
  \label{MC}
  \d\eta + \frac{1}{2} [\eta, \eta] = 0 .
\end{equation}
This is just the well known Maurer-Cartan equation! It is the basic
equation satisfied by the Maurer-Cartan form of a Lie group. In other
words, if $G$ is the (unique) simply connected Lie group integrating
$\mathfrak{g}$, then the ``tautological'' one-form
$\eta_{\text{MC}}\in \Omega^1(G; \mathfrak{g})$ is a solution of the
Maurer-Cartan equation. And it is a very special (``universal'') one
since, by a well known result in differential geometry, it classifies
all the solutions to Cartan's problem:

\begin{lemma} 
  Any one-form $\eta\in \Omega^1(N; \mathfrak{g})$ on a simply
  connected manifold $N$, satisfying the Maurer-Cartan equation
  (\ref{MC}) is of type $\eta= f^{*}\eta_{\text{MC}}$ for some smooth
  map $f: N\to G$, which is unique up to conjugation by an element
  in $G$.
\end{lemma}

It is worth keeping in mind what we have actually done to ``solve''
the original equations: the obvious necessary conditions (equation
(\ref{Cartan:eq:2}), in this case) put us into the context of Lie
algebras, then we \emph{integrate} the Lie algebra $\mathfrak{g}$ to
$G$, and finally we pick out the Maurer Cartan form as the solution we
were looking for. Moreover, to prove its universality (see the Lemma),
we had to produce a map $f: N\to G$ out of a form $\eta\in
\Omega^1(N; \mathfrak{g})$.  Viewing $\eta$ as a map $\eta: TN\to
\mathfrak{g}$, we may say that we have ``integrated $\eta$''.
\vskip 10 pt

Now what if we allow the $h$ into the picture, and try to extend this
piece of basic geometry to the general equations? This leads us
immediately to the world of algebroids! Since the structure functions
$c^{i}_{j,k}$ are no longer constant (they depend on $h\in U$), the
brackets they define (\ref{bracket}) make sense provided $\{e_i\}$
depend themselves on $h$. Of course, $\{e_i\}$ can then be viewed as
trivializing sections of an $r$-dimensional vector bundle $A$ over
$U$. On the other hand, the functions $b^{a}_{i}$ are the components
of a map:
\[ 
\Gamma(A)\ni e_{i} \longmapsto 
b^{a}_{i} \frac{\partial}{\partial x^{a}}\in\X(U). 
\] 
The necessary conditions (\ref{Cartan:eq:1}) and (\ref{Cartan:eq:2})
are just the conditions that appear in the definition of a Lie
algebroid. This is the content of the next exercise. 

\begin{exercise} 
  Let $A$ be a Lie algebroid over a manifold $M$. Pick a contractible
  open coordinate neighborhood $U\subset M$, with coordinates $(x^a)$
  and a basis of sections $\{e_i\}$ that trivialize the bundle
  $A|_U$. Also, define structure functions $c^{i}_{j,k}$ and
  $b^{a}_{i}$ by:
  \[ [e_j,e_k]=c^{i}_{j,k}e_i,\quad
  \rho(e_i)=b^{a}_{i}\frac{\partial}{\partial x^a}.\] Check that
  equation (\ref{Cartan:eq:1}) is equivalent to the condition that
  $\rho$ is a Lie algebra homomorphism:
  \[ [\rho(e_j),\rho(e_k)]=\rho([e_j,e_k]),\]
  and that equation (\ref{Cartan:eq:2}) is equivalent to the Jacobi
  identity:
  \[ [e_i,[e_j,e_k]]+[e_j,[e_k,e_i]]+[e_k,[e_i,e_j]]=0.\]
\end{exercise}
\vskip 10 pt

What about the associated global objects and the associated Maurer
Cartan forms? Our objective here is not to give a detailed discussion
of Cartan's problem, but it is not hard to guess that, in the end, you
will find Lie algebroid valued forms and that you will need to
integrate the Lie algebroid, eventually rediscovering the notion of a
\emph{Lie groupoid}.

After this motivating example, it is now time to start our study of
Lie algebroids.

\section{Lie algebroids}

Let $A$ be a Lie algebroid over $M$. We start by looking at
the kernel and at the image of the anchor $\rho:A\to TM$.

If we fix $x\in M$, the kernel $\Ker(\rho_x)$ is naturally a Lie
algebra: if $\alpha,\beta\in\Gamma(A)$ lie in $\Ker(\rho_x)$ when
evaluated at $x$, the Leibniz identity implies that $[\alpha,
f\beta](x)= f(x)[\alpha, \beta](x)$. Hence, there is a well defined
bracket on $\Ker(\rho_x)$ such that 
\[ [\alpha, \beta](x)= [\alpha(x), \beta(x)],\]
for $\alpha$ and $\beta$ as above. 

\begin{definition}
  At any point $x\in M$ the Lie algebra
  \[ \mathfrak{g}_x(A):= \Ker(\rho_x),\]
  is called the \textbf{isotropy Lie algebra} at $x$. 
\end{definition}

\begin{exercise} 
  Let $\G$ be a $s$-connected Lie groupoid over $M$ and let $A$ be its
  Lie algebroid. Show that the isotropy Lie algebra $\mathfrak{g}_x(A)$
  is isomorphic to the Lie algebra of the isotropy Lie group $\G_x$, for
  all $x\in M$.
\end{exercise}

Let us now look at the image of the anchor. This gives a distribution
\[ M\ni x\mapsto \Im(\rho_x)\subset T_xM,\]
of subspaces whose dimension, in general, will vary from point to point.
If the rank of $\rho$ is constant we say that $A$ is a 
\emph{regular Lie algebroid}.

\begin{exercise}
  If $A$ is a regular Lie algebroid show that the resulting
  distribution is integrable, so that $M$ is foliated by immersed
  submanifolds $\mathcal{O}$'s, called \emph{orbits}, satisfying
  $T_x\mathcal{O}=\Im(\rho_x)$, for all $x\in \mathcal{O}$.
\end{exercise}

For a general Lie algebroid $A$ there is still an induced partition of
$M$ by immersed submanifolds, called the \emph{orbits} of $A$. As in the
regular case, one looks for maximal immersed submanifolds
$\mathcal{O}$'s satisfying:
\[ T_x\mathcal{O}=\Im(\rho_x), \]
for all $x\in \mathcal{O}$. Their existence is a bit more delicate in
the non-regular case. One possibility is to use a Frobenius type
theorem for singular foliations. However, there is a simple way of
describing the orbits, at least set theoretically. This uses the
notion of $A$-path, which will play a central role in the
integrability problem to be studied in the next lecture.

\begin{definition} 
  Given a Lie algebroid $A$ over $M$, an \textbf{$A$-path} consists of a
  pair $(a,\gamma)$ where $\gamma:I\to M$ is a path in $M$, $a: I\to
  A$ is a path in $A$(\footnote{Here and below, $I=[0,1]$ will always
  denote the unit interval.}), such that
  \begin{enumerate}[(i)]
  \item $a$ is a path above $\gamma$, i.e., $a(t)\in A_{\gamma(t)}$
  for all $t\in I$.
  \item $\rho(a(t))= \frac{\d\gamma}{\d t}(t)$, for all $t\in I$.
  \end{enumerate}
\end{definition}

Of course, $a$ determines the base path $\gamma$, hence when talking
about an $A$-path we will only refer to $a$. The reason we mention
$\gamma$ in the previous definition is to emphasize the way we think
of $A$-paths: $a$ should be interpreted as an ``$A$-derivative of
$\gamma$'', and the last condition in the definition should be read:
``the usual derivative of $\gamma$ is related to the $A$-derivative by
the anchor map''.

We can now define an equivalence relation on $M$, denoted $\sim_{A}$,
as follows. We say that $x,y\in M$ are equivalent if there exists an
$A$-path $a$, with base path $\gamma$, such that $\gamma(0)= x$ and
$\gamma(1)= y$. An equivalence class of this relation will be
called an \textbf{orbit of $A$}. When $\rho$ is surjective we say
that $A$ is a \textbf{transitive Lie algebroid}. In this case,
each connected component of $M$ is an orbit of $A$. 

\begin{remark} 
  It remains to show that each orbit $\O$ is an immersed submanifold
  of $M$, which integrates $\Im(\rho)$, i.e., $T_x\O=\Im(\rho_x)$ for
  all $x\in \mathcal{O}$. This can be proved using a local normal form
  theorem for Lie algebroids.
\end{remark}

\begin{exercise} 
  Let $\G$ be a $s$-connected Lie groupoid over $M$ and let $A$ be its
  Lie algebroid. Show that the orbits of $\G$ in $M$ coincide with the
  orbits of $A$ in $M$.\\
  (Hint: Check that if $g(t):I\to\G$ is a path that stays in a
  $s$-fiber and starts at $1_x$, then 
  $a(t)=\d_{g(t)}R_{g(t)^{-1}}\cdot\dot{g}(t)$
  is an $A$-path.)
\end{exercise}

As we have seen in the previous lecture, any Lie groupoid has an
associated Lie algebroid. For future reference, we introduce the
following terminology:

\begin{definition} 
  A Lie algebroid $A$ is called \textbf{integrable} if it is
  isomorphic to the Lie algebroid of a Lie groupoid $\G$. For such a
  $\G$, we say that $\G$ integrates $A$.
\end{definition}

Similar to Lie's first theorem for Lie algebras (which asserts that
there is at most one simply-connected Lie group integrating a given
Lie algebra), we have the following theorem:

\begin{theorem}[Lie I] 
  If $A$ is integrable, then there exists an unique (up to isomorphism)
  $s$-simply connected Lie groupoid $\G$ integrating $A$.
\end{theorem}

\begin{proof}
  This follows at once from the previous lecture, where we have shown
  that for any Lie groupoid $\G$ there exists a unique (up to
  isomorphism) Lie groupoid $\widetilde{\G}$ which is $s$-simply
  connected and which has the same Lie algebroid as $\G$ (see Theorem
  \ref{thm:1:connected}).
\end{proof}

Let us turn now to morphisms of Lie algebroids.

\begin{definition}
  Let $A_1\to M_1$ and $A_2\to M_2$ be Lie algebroids. A
  \textbf{morphism of Lie algebroids} is a vector bundle map 
  \[
  \xymatrix{
    A_1\ar[r]^{F}\ar[d]& A_2\ar[d]\\
    M_1\ar[r]_{f}&M_2 
  }\]
  which is compatible with the anchors and the brackets.
\end{definition}

Let us explain what we mean by \emph{compatible}. First of all, we say
that the map $F:A_1\to A_2$ is compatible with the anchors if
\[ \d f(\rho_{A_1}(a))=\rho_{A_2}(F(a)). \]
This can be expressed by the commutativity of the diagram:
\[
\xymatrix{
    A_1\ar[r]^{F}\ar[d]_{\rho_1}& A_2\ar[d]^{\rho_2}\\
    TM_1\ar[r]_{\d f}&TM_2 
}\]
Secondly, we would like to say what we mean by compatibility with
the brackets. The difficulty here is that, in general, sections of
$A_1$ cannot be pushed forward to sections of $A_2$. Instead
we have to work at the level of the pull-back bundle $f^{*}A_2$ (\footnote{an alternative, more intrinsic definition, will be described in Exercise
\ref{morphisms:algebroids}}). First note
that from sections $\al$ of $A_1$ or $\al'$ of $A_2$, we can produce 
new sections $F(\al)$ and $f^{*}(\al')$ of $f^{*}A_2$ by:
\[ F(\al)= F\circ \al,\ f^{*}(\al')= \al'\circ f.\]
Now, given any section $\al\in\Gamma(A_1)$, we can express its image
under $F$ as a (non-unique) finite combination 
\[ F(\al)=\sum_{i} c_i f^{*}(\al_i) ,\]
where $c_i\in C^\infty(M_1)$ and $\al_i\in\Gamma(A_2)$.
By compatibility with the brackets we mean that, if
$\al,\be\in\Gamma(A_1)$ are sections such that their images are
expressed as finite combinations as above, then their bracket is a
section whose image can be expressed as:
\begin{multline}
\label{eq:compatibility}
  F([\al,\be]_{A_1})=
  \sum_{i,j} c_i c'_{j}f^{*}[\al_i,\be_{j}]_{A_2}+\\
  +\sum_j \Lie_{\rho(\al)}(c'_{j})f^{*}(\be_{j})-
  \sum_i \Lie_{\rho(\be)}(c_{i})f^{*}(\al_{i}).           
\end{multline}
Notice that, in the case where the sections $\al,\be\in\Gamma(A_1)$
can be pushed forward to sections $\al',\be'\in\Gamma(A_2)$, so that 
$F(\al)=\al'\circ f$ and $F(\be)=\be'\circ f$, this
just means that:
\[ F([\al,\be]_{A_1})=[\al',\be']_{A_2}\circ f. \]

The following exercises should help make you familiar with the notion
of a Lie algebroid morphism.

\begin{exercise}
  Check that condition (\ref{eq:compatibility}) is independent of the
  way one expresses the image of the sections under $F$ as finite
  combinations.
\end{exercise}

\begin{exercise}
  Let $A$ be a Lie algebroid. Show that a path $a:[0,1]\to A$ is an
  $A$-path iff the map $a\d t:TI\to A$ is a morphism of Lie algebroids.
\end{exercise}

\begin{exercise}
  \label{ex:morphisms}
  Show that if $\F:\G\to\H$ is a homomorphism of Lie groupoids, then it
  induces a Lie algebroid homomorphism $F:A\to B$ of their Lie algebroids.
\end{exercise}

Again, just like in the case of Lie algebras, under a suitable
assumption, we can integrate morphisms of Lie algebroids to morphisms
of Lie groupoids:

\begin{theorem}[Lie II]
  \label{thm:Lie:II}
  Let $F: A\to B$ be a morphism of integrable Lie algebroids, and
  let $\G$ and $\H$ be integrations of $A$ and $B$. If $\G$ is
  $\s$-simply connected, then there exists a (unique) morphism of Lie
  groupoids $\F:\G\to \H$ integrating $F$.
\end{theorem}

You may try to reproduce the proof that you know for the Lie algebra
case. You will also find a proof in the next lecture (see Exercise
\ref{exer:Lie:II}). The next exercise shows how to define the
exponential map for Lie algebroids/groupoids.

\begin{exercise}
  Let $A$ be the Lie algebroid of a Lie groupoid $\G$. Use Lie II to
  define the exponential map $\exp:\Gamma_c(A)\to\Gamma(\G)$, taking
  sections of compact support to bisections of $\G$. How does this
  relate to the construction hinted at in Remark \ref{rem:exponential}?
\end{exercise}

\section{First examples of Lie algebroids}

Let us present now a few basic examples of Lie algebroids.

\begin{example}[tangent bundles] 
  One of the extreme examples of a Lie algebroid over $M$ is the
  tangent bundle $A=TM$, with the identity map as anchor, and the
  usual Lie bracket of vector fields. Here the isotropy Lie algebras
  are trivial and the Lie algebroid is transitive.

\begin{exercise} 
  Prove that the Lie algebroids of the pair groupoid $M\times M$ and
  of the fundamental groupoid $\Pi_1(M)$ are both isomorphic to $TM$.
\end{exercise}
\end{example}

\begin{example}[Lie algebras] 
  At the other extreme, any Lie algebra $\gg$ is a Lie algebroid over a
  singleton. Here there is only one isotropy Lie algebra which
  coincides with $\gg$. Obviously, any Lie group with Lie algebra $\gg$
  gives a Lie groupoid integrating $A$.
\end{example}

Both these examples can be slightly generalized. For example, the
tangent bundle can be generalized as follows:

\begin{example}[foliations]
  Let $A\subset TM$ be an involutive subbundle, i.e., constant rank
  smooth distribution which is closed for the usual Lie bracket. This
  gives a Lie algebroid (in fact a Lie subalgebroid of $TM$) over $M$,
  with anchor map the inclusion, and the Lie bracket the restriction
  of the usual Lie bracket of vector fields. 

  Recall that, by the Frobenius Integrability Theorem, $A$
  determines a foliation $\F$ of a manifold $M$ (and conversely, every
  foliation determines a Lie algebroid $T\F\subset TM$). In fact, $\F$
  is just the orbit foliation of $A$ (see the discussion above). On
  the other hand, since the anchor is injective, the isotropy Lie
  algebras $\gg_x$ are all trivial.
\end{example}

\begin{exercise}
  Prove that the Lie algebroid of the fundamental groupoid $\Pi_1(\F)$
  is isomorphic to $T\F$. Can you give another example of a groupoid
  integrating $T\F$?
\end{exercise}

On the other hand, Lie algebras can be generalized as follows:

\begin{example}[bundles of Lie algebras] 
  A bundle of Lie algebras over $M$ is a vector bundle $A$ over $M$
  together with a Lie algebra bracket $[\cdot, \cdot]_x$ on each fiber
  $A_x$, which varies smoothly with respect to $x$ in the sense that
  if $\alpha, \beta\in \Gamma(A)$, then $[\alpha, \beta]$ defined by
  \[ [\alpha, \beta](x)= [\alpha(x), \beta(x)]_{x}\]
  is a smooth section of $A$. Note that this notion is weaker then
  that of Lie algebra bundle, when one requires that $A$ is locally
  trivial as a bundle of Lie algebras (in particular, all the Lie
  algebras $A_x$ should be isomorphic).
  
  Since the anchor is identically zero, the orbits of $A$ are the
  points of $M$, while the isotropy Lie algebras are the fibers $A_x$.
  It is easy to see that a bundle of Lie algebras over $M$ is
  precisely the same thing as a Lie algebroid over $M$ with zero
  anchor map.
\end{example}

A general Lie algebroid can be seen as combining aspects from both the
previous two examples. In fact, take a Lie algebroid $A$ over $M$ and
fix an orbit $i:\mathcal{O}\hookrightarrow M$. In the following exercise
we ask you yo check that the bracket restricts to a bracket on 
$A_{\mathcal{O}}= A|_{\mathcal{O}}:=i^*A$.

\begin{exercise}
Let $\al,\be\in\Gamma(A_\mathcal{O})$ be local sections, and
$\widetilde{\al},\widetilde{\be}\in\Gamma(A)$ be any choice of
local extensions. Show that:
\[ [\al,\be]:={[\widetilde{\al},\widetilde{\be}]}|_{\mathcal{O}},\]
is well-defined, i.e., it is independent of the choice of extensions.\\
(Hint: Since $\mathcal{O}$ is a leaf, the vector fields $\rho(\al)$
and $\rho(\be)$ are tangent to $\mathcal{O}$.)
\end{exercise}

Since the restriction of the anchor gives a map $\rho_{\mathcal{O}}:A_\mathcal{O}\to T\mathcal{O}$,
we obtain a Lie algebroid structure on $A_{\mathcal{O}}$ over $\mathcal{O}$. 
On the other hand, one has an induced bundle of Lie algebras:
\[ \mathfrak{g}_{\mathcal{O}}(A)=\Ker(\rho_{\mathcal{O}}),\]
whose fiber at $x$ is the isotropy Lie algebra $\mathfrak{g}_x(A)$.
All these fit into a short exact sequence of algebroids over $\mathcal{O}$:
\[ 0\rmap \mathfrak{g}_{\mathcal{O}}(A)\rmap A_{\mathcal{O}}\rmap
T\mathcal{O} \rmap 0 .\]
The isotropy bundle $\mathfrak{g}_{\mathcal{O}}(A)$ is a Lie algebra
bundle, as indicated in the following exercise.

\begin{exercise}
  \begin{enumerate}[(a)]
  \item Show that a bundle of Lie algebras $\gg_M\to M$ is a Lie
    algebra bundle iff there exists a 
    connection $\nabla: \X(M)\otimes \Gamma(\gg_M)\to
    \Gamma(\gg_M)$ satisfying
    \[ \nabla_{X}([\alpha, \beta])= [\nabla_{X}(\alpha),\beta]+ 
        [\alpha, \nabla_{X}(\beta)]\]
    for all $\alpha, \beta\in \Gamma(\gg_M)$.
  \item Deduce that each isotropy bundle $\mathfrak{g}_{\mathcal{O}}(A)$
    is a Lie algebra bundle. In particular, all the isotropy Lie algebras
    $\mathfrak{g}_{x}(A)$, with $x\in \mathcal{O}$, are isomorphic.
  \item[ ] (Hint: for the first part, use parallel transport. For the second
    part, use a splitting of $\rho_{\mathcal{O}}: A|_{\mathcal{O}}\to
    T\mathcal{O}$ and the bracket of $A_{\mathcal{O}}$ to produce a connection).
  \end{enumerate}
\end{exercise}

\begin{example}[vector fields] 
  It is not difficult to see that Lie algebroid structures on the
  trivial line bundle over $M$ are in 1-1 correspondence with vector
  fields on $M$. Given a vector field $X$, we denote by $A_X$ the
  induced Lie algebroid. Explicitly, as a vector bundle, $A_X=\L=
  M\times\mathbb{R}$ is the trivial line bundle, while the anchor is
  given by multiplication by $X$, and the Lie bracket of two section
  $f,g\in \Gamma(A_X)=C^{\infty}(M)$ is defined by:
  \[ [f, g]= f \Lie_{X}(g)- \Lie_{X}(f)g.\]

\begin{exercise} 
  Find out how the flow of $X$ defines an integration of $A_X$.
\end{exercise}
\end{example}

\begin{example}[action Lie algebroid] 
  Generalizing the Lie algebroid of a vector field, consider an
  infinitesimal action of a Lie algebra $\mathfrak{g}$ on a manifold
  $M$, i.e., a Lie algebra homomorphism $\rho: \mathfrak{g}\to
  \X(M)$. The standard situation is when $\mathfrak{g}$ is the Lie
  algebra of a Lie group $G$ that acts on $M$. Then
  \[ \rho(v)(x):= \frac{\d}{\d t}\exp(tv) x,\quad (v\in \mathfrak{g}, x\in M)\]
  defines an infinitesimal action of $\mathfrak{g}$ on $M$. We say
  that an infinitesimal action of $\mathfrak{g}$ on $M$ is integrable
  if it comes from a Lie group action.

  Given an infinitesimal action of $\mathfrak{g}$ on $M$, we form a Lie
  algebroid $\mathfrak{g}\ltimes M$, called the action Lie algebroid,
  as follows. As a vector bundle, it is the trivial vector bundle
  $M\times \mathfrak{g}$ over $M$ with fiber $\mathfrak{g}$, the
  anchor is given by the infinitesimal action, while the Lie bracket is
  uniquely determined by the Leibniz identity and the condition that
  \[ [c_v, c_w]=c_{[v, w]},\]
  for all $v, w\in \mathfrak{g}$, where $c_v$ denotes the constant
  section of $\mathfrak{g}$.
  
\begin{exercise} This exercise discusses the relationship between the
  integrability of Lie algebra actions and the corresponding action
  Lie algebroid: 
  \begin{enumerate}[(a)]
  \item Given an action of a Lie group $G$ on $M$, show that the
    corresponding action Lie groupoid $G\ltimes M$ over $M$ has Lie
    algebroid the action Lie algebroid $\mathfrak{g}\ltimes M$.
    Hence, if an infinitesimal action of $\mathfrak{g}$ on $M$ is
    integrable, then the Lie algebroid $\mathfrak{g}\ltimes M$ is
    integrable.
  \item Find an infinitesimal action of a Lie algebra $\mathfrak{g}$
    which is not integrable but which has the property that the Lie algebroid 
    $\mathfrak{g}\ltimes M$ is integrable. 
  \item[] (Hint: think of vector fields!)
\end{enumerate}
\end{exercise}

\begin{exercise}
  Show that \emph{all} action Lie algebroids $\mathfrak{g}\ltimes M$, arising
  from infinitesimal Lie algebra actions are integrable.\\
  (Hint: think again of what happens for vector fields!).
\end{exercise}
\end{example}

\begin{example}[Two forms]
  \label{ex:2-forms}
  Any closed 2-form $\omega$ on a manifold $M$ has an associated Lie
  algebroid, denoted $A_{\omega}$, and defined as follows. As a vector
  bundle,
  \[ A= TM\oplus\L,\]
  the anchor is the projection on the first component, while the
  bracket on sections $\Gamma(A_\omega)\simeq\X(M)\times C^\infty(M)$
  is defined by:
  \[ [(X,f),(Y,g)]= ([X,Y],\Lie_{X}(g)-\Lie_{Y}(f)+\omega(X,Y)) .\]

  \begin{exercise} 
    Given a 2-form on $M$, check that the previous formulas make
    $A_{\omega}$ into a Lie algebroid if and only if $\omega$ is closed.
  \end{exercise}

\end{example}

\begin{example}[Atiyah sequences] 
  Let $G$ be a Lie group. To any principal $G$-bundle $P$ over $M$
  there is an associated Lie algebroid over $M$, denoted $A(P)$, and
  defined as follows. As a vector bundle, $A(P):= TP/G$ (over $P/G=
  M$). The anchor is induced by the differential of the projection
  from $P$ to $M$. Also, since the sections of $A(P)$ correspond to
  $G$-invariant vector fields on $P$, we see that there is a canonical
  Lie bracket on $\Gamma(A(P))$.  With these, $A(P)$ becomes a Lie
  algebroid. 

  \begin{exercise}
    Show that $A(P)$ is just the Lie algebroid of the gauge groupoid
    $P\otimes_G P$.
  \end{exercise}

  Note that $A(P)$ is transitive, i.e., the anchor map is
  surjective. We denote by $P[\mathfrak{g}]$ the kernel of the anchor
  map.  Hence we have a short exact sequence
  \[ 0\rmap P[\mathfrak{g}]\rmap A(P)\rmap TM \rmap 0,\] 
  known as the Atiyah sequence associated to $P$. Of course, this is
  obtained by dividing out the action of $G$ on the exact sequence of
  vector bundles over $P$: $P\times \mathfrak{g}\rmap TP\rmap TM$.  In
  particular, $P[\mathfrak{g}]$ is the Lie algebra bundle obtained by
  attaching to $P$ the adjoint representation of $G$.

  One of the interesting things about the Atiyah sequence is its
  relation with connections and their curvatures.  First of all,
  connections on the principal bundle $P$ are the same thing as
  splittings of the Atiyah sequence.  Indeed, a left splitting of the
  sequence is the same thing as a bundle map $TP\to P\times
  \mathfrak{g}$ which is $G$-invariant, and which is left inverse
  to the infinitesimal action $\mathfrak{g}\to TP$, in other words,
  a connection 1-form. By standard linear algebra, left splittings of
  a short exact sequence are in 1-1 correspondence with right
  splittings. Henceforth, we will identify connections with left/right
  splittings.

  Next, given a connection whose associated right splitting is denoted
  by $\sigma: TM\to A(P)$, the curvature of the connection is a 2-form
  on $M$ with values in $P[\mathfrak{g}]$, and can be described using
  $\sigma$ and the Lie bracket on $A(P)$:
  \[ \Omega_{\sigma}(X,Y)= [\sigma(X),\sigma(Y)]-\sigma([X,Y]).\]

  Using these as motivations, transitive (i.e., with surjective anchor)
  Lie algebroids $A$ are also called abstract Atiyah sequence, while
  splittings of their anchor map are called connections. The sequence
  associated to $A$ is, of course:
  \[ 0\rmap \Ker(\rho)\rmap A\rmap TM\rmap 0.\]
  Note that not all abstract Atiyah sequences come from principal
  bundles.

\begin{exercise}
  Give an example of a transitive Lie algebroid which is not
  the Atiyah sequence of some principal bundle.
\end{exercise}
\end{example}

\begin{example}[Poisson manifolds]
  \label{ex:Poisson}
  Any Poisson structure on a manifold $M$ induces a Lie algebroid
  structure on $T^*M$ as follows.  Let $\pi$ be the Poisson bivector
  on $M$, which is related to the Poisson bracket by $\{f,g\}=\pi(\d
  f,\d g)$. Also, the hamiltonian vector field $X_f$ associated to a
  smooth function $f$ on $M$ is given by $X_f(g)= \{f, g\}$. We use
  the notation
  \[ \pi^{\sharp}: T^*M\to TM \]
  for the map defined by $\beta(\pi^{\sharp}(\alpha))= \pi(\alpha,
  \beta)$.  The Lie algebroid structure on $T^*M$ has $\pi^{\sharp}$
  as anchor map, and the Lie bracket is defined by
  \[ [\alpha, \beta]= \Lie_{\pi^{\sharp}(\alpha)}(\beta)- 
  \Lie_{\pi^{\sharp}(\beta)}(\alpha)- \d(\pi(\alpha, \beta)).\] 
  We call this the \textbf{cotangent Lie algebroid} of the Poisson
  manifold $(M,\pi)$.

  \begin{exercise} 
    \label{ex:Koszul}
    Show that this Lie algebroid structure on $T^*M$ is the unique
    one with the property that the anchor maps $\d f$ to $X_{f}$ and
    $[\d f,\d g]=\d\{f,g\}$, for all $f, g\in C^{\infty}(M)$.
  \end{exercise}
\end{example}

\begin{example}[Nijenhuis tensors] 
  Given a bundle map $\mathcal{N}: TM\to TM$, recall that its Nijenhuis
  torsion, denoted $T_\mathcal{N}\in \Gamma(\wedge^2T^*M\otimes TM)$, 
  is defined by
  \[ T_\mathcal{N}(X,Y)=[\mathcal{N}X,\mathcal{N}Y]-
  \mathcal{N}[\mathcal{N}X,Y]-\mathcal{N}[X,\mathcal{N}Y]+\mathcal{N}^2[X,Y],\]
  for $X,Y\in\X(M)$. When $T_\mathcal{N}=0$ we call $\mathcal{N}$ a
  \textbf{Nijenhuis tensor}. 
  {To} any Nijenhuis tensor $\mathcal{N}$, there is associated a new
  Lie algebroid structure on $TM$: the anchor is given by 
  $\rho(X)=\mathcal{N}(X)$, while the Lie bracket is defined by
  \[ 
  [X,Y]_\mathcal{N}:=[\mathcal{N}X,Y]+[X,\mathcal{N}Y]-
  \mathcal{N}([X,Y]). 
  \]
  This kind of Lie algebroid plays an important role in the theory of
  integrable systems, in the study of (generalized) complex  structures,
  etc.
  \begin{exercise}
    Thinking of a Lie algebroid as a generalized tangent bundle,
    extend this construction to any Lie algebroid $A\to M$: given a bundle
    map $\mathcal{N}:A\to A$ over the identity, such that its Nijenhuis torsion
    vanishes, i.e.,
    \[
    [\mathcal{N}\al,\mathcal{N}\be]_A-\mathcal{N}[\mathcal{N}\al,\be]_A-
    \mathcal{N}[\al,\mathcal{N}\be]_A+\mathcal{N}^2[\al,\be]_A=0, 
    \]
    for all $\al,\be\in\Gamma(A)$, show that there exists a new Lie
    algebroid $A_\mathcal{N}$ associated with $\mathcal{N}$.
  \end{exercise}
\end{example}

\section{Connections}

Concepts in Lie algebroid theory arise often as generalizations of
standard notions both of Lie theory and/or of differential geometry. This
is related with the two extreme examples of Lie algebroids: Lie
algebras and tangent bundles. This dichotomy affects also notation and
terminology as we will see now.

Having the example of the tangent bundle in mind, let us interpret a
general Lie algebroid as describing a space of ``generalized vector
fields on $M$'' ($A$-fields), so that the anchor map relates these
generalized vector fields to the usual vector fields. This leads
immediately to the following notion of connection:

\begin{definition} 
  Given a Lie algebroid $A$ over $M$ and a vector bundle $E$ over $M$,
  an \textbf{$A$-connection} on $E$ is a bilinear map
  \begin{align*}
    &\Gamma(A)\times \Gamma(E)\to \Gamma(E), \\
    &(\alpha, s)\mapsto\nabla_{\alpha}(s),
  \end{align*}
  which is $C^{\infty}(M)$ linear on $\alpha$, and which satisfies the
  following Leibniz rule with respect to $s$:
  \[ \nabla_{\alpha}(fs)= f\nabla_{\alpha}(s)+
  \Lie_{\rho(\alpha)}(f) s.\] 
  The \textbf{curvature} of the connection is the
 map:
  \begin{align*}
    &R_\nabla:\Gamma(A)\times \Gamma(A)\to \Hom(\Gamma(E), \Gamma(E)),\\
    &R_\nabla(\al,\be)(X)=\nabla_{\alpha}\nabla_{\beta}X-
    \nabla_{\beta}\nabla_{\alpha}X-\nabla_{[\alpha, \beta]}X.
  \end{align*}
  The connection is called \textbf{flat} if $R_\nabla=0$. 
\end{definition}

\begin{exercise}
Check that $R_{\nabla}(\al, \be)(X)$ is $C^{\infty}(M)$-linear in
$\al$, $\be$ and $X$. In other words,
\[ R_{\nabla}\in \Gamma(\wedge^2A^*\otimes \End(E)).\]
\end{exercise}

Note that a $TM$-connection is just an ordinary connection, and all
notions we have introduced (and that we will introduce!) reduce to
well-known notions of ordinary connection theory.

Given an $A$-path $a$ with base path $\gamma: I\to M$,
and $u:I\to E$ a path in $E$ above $\gamma$, then
the derivative of $u$ along $a$, denoted $\nabla_{a}u$,
is defined as usual: choose a time dependent section $\xi$ of $E$
such that $\xi(t,\gamma(t))= u(t)$, then
\[ 
\nabla_{a}u(t)= \nabla_{a}\xi^{t}(x)+\frac{\d\xi^{t}}{\d t}(x),
\text{ at } x=\gamma(t).
\]

\begin{exercise}
  \label{parallel-transport} 
  For an $A$-connection $\nabla$ on a vector bundle $E$ define
  parallel transport along $A$-paths. When is a connection complete
  (i.e., when is parallel transport defined for every $A$-path)?
\end{exercise}

The following exercise gives an alternative aproach to connections 
in terms of the principal frame bundle. It also suggests how to define
Ehresmann $A$-connections on any principal $G$-bundle.

\begin{exercise}
Let $\nabla$ be an $A$-connection on a vector bundle $E$ of rank $r$. 
If $P:=\GL(E)\to M$ is the bundle of linear frames on $E$, 
show that $\nabla$ induces a smooth bundle map $h:p^*A\to TP$, such that:
\begin{enumerate}[(i)]
\item $h$ is horizontal, i.~e., the following diagram commutes:
\[
\xymatrix{
p^*A\ar[r]^{h}\ar[d]_{\widehat{p}}& TP \ar[d]^{p_*} \\ A\ar[r]_{\#}
&TM }
\]
\item $h$ is $\GL(r)$-invariant, i.~e., we have
\[ h(u g,a)=(R_g)_*h(u,a), \qquad \text{ for all }g\in \GL(r);\]
\end{enumerate}
Conversely, show that every smooth bundle map $h:p^*A\to TP$ satisfying
(i) and (ii) induces a connection $\nabla$ on $E$.
\end{exercise}

Of special importance are the $A$-connections on $A$. If $\nabla$ is
an $A$-connection on $A$, we can define its \textbf{torsion} to be the
$C^\infty(M)$-bilinear map:
  \begin{align*}
    &T_\nabla:\Gamma(A)\times \Gamma(A)\to \Gamma(A),\\
    &T_\nabla(\al,\be)=\nabla_{\alpha}\beta-\nabla_{\beta}\alpha
    -[\alpha,\beta].
  \end{align*}
We say that $\nabla$ is torsion free if $T_{\nabla}= 0$.

\begin{exercise}
Check that $T_{\nabla}(\al, \be)$ is $C^{\infty}(M)$-linear in
$\al$ and $\be$. In other words,
\[ T_{\nabla}\in \Gamma(\wedge^2A^*\otimes A).\]
\end{exercise}

We will use later some special connections $A$-connections which arise
once an ordinary ($TM$-) connection on $A$ is fixed. These constructions
are explained in the following two exercises:

\begin{exercise}
  \label{exer:connections}
  Let $A$ be a Lie algebroid and let $\nabla$ be a $TM$-connection on
  the vector bundle $A$. Show that:
  \begin{enumerate}[(a)] 
  \item The following formula defines an $A$-connection on the vector
    bundle $A$: 
    \[ \nabla_{\alpha}\beta\equiv \nabla_{\rho(\alpha)}\beta.\]
  \item The following formula defines an $A$-connection on the vector
    bundle $A$: 
    \[ \overline{\nabla}_{\alpha}\beta\equiv
    \nabla_{\rho(\beta)}\alpha+ [\alpha, \beta].\]
  \item The following formula defines an $A$-connection on the vector
    bundle $TM$:
    \[ \overline{\nabla}_{\alpha}X\equiv \rho(\nabla_{X}\alpha)+
    [\rho(\alpha), X] .\]
  \end{enumerate}
  Note that the
  last two connections are compatible with the anchor:
  \[\overline{\nabla}_{\alpha}\rho(\beta)=
  \rho(\overline{\nabla}_{\alpha}\beta).\]
\end{exercise}

\begin{exercise}[Levi-Civita connections]
Let $\langle\cdot, \cdot\rangle$ be a metric on $A$. We say that an
$A$-connection $\nabla$ on $A$ is compatible with the metric if
\[ 
\Lie_{\rho(\al)}(\langle\beta,\gamma\rangle)= 
\langle\nabla_{\al}\be,\gamma\rangle+ 
\langle\be,\nabla_{\al}\gamma\rangle,
\]
for all sections $\al, \be, \gamma$ of $A$. Show that $A$ admits a
unique torsion free $A$-connection compatible with the metric.
\end{exercise}

Finally we remark that one can give a more intrinsic description
of algebroid homomorphisms by using $A$-connections.

\begin{exercise}
\label{morphisms:algebroids}
Let $A_{1}$ and $A_2$ be Lie algebroids over $M_1$, and $M_2$, respectively,
and let $F:A_1\to A_2$ be a bundle map covering $f: M_1\to M_2$
which is compatible with the anchor.      
\begin{enumerate}[(a)]
\item For a connection $\nabla$ on $A_2$, we denote by the same letter the
pull-back of $\nabla$ via $F$ (a connection in $f^{*}A_2$). 
Show that the expression
\[ R_F(\al, \be):= \nabla_{\al}(F(\be))- 
                              \nabla_{\be}(F(\al))- 
                              F([\al, \be])- 
                              T_{\nabla}(F(\al), F(\be))\]
defines an element
\[ R_F\in \Gamma(\wedge^2A_{1}^{*}\otimes f^{*}A_2)\]
which is independent of the choice of $\nabla$.
\item Show that $F$ is a Lie algebroid homomorphism if and only of $R_F= 0$.
\item Deduce that, given a manifold $M$ and a Lie algebra $\mathfrak{g}$, an
algebroid morphism $TM\to \mathfrak{g}$ is the same thing as a $1$-form
$\omega\in \Omega^1(M; \mathfrak{g})$ satisfying the Maurer-Cartan equation 
$\d\omega+ \frac{1}{2}[\omega, \omega]= 0$. What does $\omega$ integrate to?
\item Back to the general context, make sense of a similar Maurer-Cartan
formula
\[ R_F=\d_{\nabla}(F)+ \frac{1}{2}[F, F]_{\nabla}.\]
\end{enumerate}
\end{exercise}

\section{Representations}

While the notion of connection is motivated by the tangent bundle, the
following notion is motivated by Lie theory:

\begin{definition} 
  A \textbf{representation} of a Lie algebroid $A$ over $M$ consists of a
  vector bundle $E$ over $M$ together with a flat connection
  $\nabla$, i.e., a connection such that
  \[ \nabla_{\alpha}\nabla_{\beta}- \nabla_{\beta}\nabla_{\alpha}=
  \nabla_{[\alpha, \beta]} ,\]
  for all $\alpha, \beta\in \Gamma(A)$.
\end{definition}

The isomorphism classes of representations of a Lie algebroid $A$ form
a semiring $(\Rep(A),\oplus,\otimes)$, where we use the direct sum and
tensor product of connections to define the addition and the
multiplication: given connections $\nabla^{i}$ on $E^i$, $i=1, 2$,
$E^1\oplus E^2$ and $E^1\otimes E^2$ have induced connections
$\nabla^{\oplus}$, $\nabla^{\otimes}$, given by
\begin{align*}
\nabla^{\oplus}_{\alpha}(s_1, s_2)&= 
(\nabla^{1}_{\alpha}(s_1),\nabla^{2}_{\alpha}(s_2)),\\
\nabla^{\otimes}_{\alpha}(s_1\otimes s_2)&= 
\nabla^{1}_{\alpha}(s_1)\otimes s_2+s_1\otimes \nabla^{2}_{\alpha}(s_2).
\end{align*}
The identity element is the trivial line bundle $\mathbb{L}$ over $M$
with the flat connection $\nabla_{\alpha}:= \Lie_{\alpha}$.

The notion of a Lie algebroid representation is the infinitesimal
counterpart of a representation of a Lie groupoid.  To make this more
precise, we denote by $\mathfrak{gl}(E)$ the Lie algebroid of $\GL(E)$. 
This Lie algebroid is related to the Lie algebra 
$\Der(E)$ of derivations on $E$.

\begin{definition} 
Given a vector bundle $E$ over $M$, a derivation of $E$ is a pair $(D,
X)$ consisting of a vector field $X$ on $M$ and a linear map
\[ D: \Gamma(E)\to \Gamma(E),\]
satisfying the Leibniz identity
\[ D(fs)= fDs+ \Lie_{X}(f)s,\]
for all $f\in C^{\infty}(M)$, $s\in \Gamma(E)$. We denote by
$\Der(E)$ the Lie algebra of derivations of $E$, where the Lie
bracket is given by
\[ [(D,X),(D',X')]=(DD'-D'D,[X,X']).\]
\end{definition}

\begin{lemma} 
  Given a vector bundle $E$ over $M$, the Lie algebra of sections of the
  algebroid $\mathfrak{gl}(E)$ is isomorphic to the Lie algebra
  $\Der(E)$ of derivations on $E$, and the anchor of
  $\mathfrak{gl}(E)$ is identified with the projection $(D, X)\mapsto X$.
\end{lemma}

\begin{proof}
  Given a section $\alpha$ of $\mathfrak{gl}(E)$,
  its flow (see Definition \ref{flow-al}) gives linear maps
  \[ \phi_{\al}^{t}(x): E_x\to E_{\phi_{\rho\al}^{t}(x)}.\]
  In particular, we obtain a map at the level of sections,
  \begin{align*}
    (\phi_{\al}^{t})^{*} &: \Gamma(E)\to \Gamma(E), \\
    (\phi_{\al}^{t})^{*}(s)&= 
    \phi_{\alpha}^{t}(x)^{-1}(s(\phi_{\rho\al}^{t}(x)).
  \end{align*}
  The derivative of this map induces a derivation $D_{\alpha}$ on
  $\Gamma(E)$:
  \[ D_{\al}s=\left.\frac{\d}{\d t}\right|_{t=0} (\phi_{\al}^{t})^{*}(s).\]
  Conversely, given a derivation $D$, we find a $1$-parameter group
  $\phi_{D}^{t}$ of automorphisms of $E$, sitting over the flow
  $\phi_{X}^{t}$ of the vector field $X$ associated to $D$, as the
  soluion of the equation
  \[ D s=\left.\frac{\d}{\d t}\right|_{t=0} (\phi_{D}^{t})^{*}(s).\]
  Viewing $\phi_{D}^{t}(x)$ as an element of $\G$, and differentiating
  with respect to $t$ at $t= 0$, we obtain a section $\al_{D}$ of
  $\mathfrak{gl}(E)$. Clearly, the correspondences $\alpha\mapsto
  D_{\alpha}$ and $D\mapsto \alpha_{D}$ are inverse to each other,
  hence we have an isomorphism between $\Gamma(\mathfrak{gl}(E))$ and
  $\Der(E)$. One way to see that this preserves the bracket is by
  using the description of Lie brackets in terms of
  flows. Alternatively, one remarks that the correspondences we have
  defined are local. Hence we may assume that $E$ is trivial as a
  vector bundle, in which case the computation is simple and can be
  left to the reader.
\end{proof}

\begin{exercise}
  Show that a representation $(E,\nabla)$ of a Lie algebroid $A$ is
  the same thing as a Lie algebroid homomorphism $\nabla:A\to \mathfrak{gl}(E)$.
\end{exercise}

The infinitesimal version of a Lie groupoid homomorphism $\G\to\GL(E)$
is a Lie algebroid homomorphism $A\to \mathfrak{gl}(E)$ (Exercise
\ref{ex:morphisms}). By the previous exercise, such a homomorphism is
the same thing as a flat $A$-connection on $E$, so we deduce the following.

\begin{corollary}
  Let $\G$ be a Lie groupoid over $M$ and let $A$ be its Lie algebroid.
  Any representation $E$ of $\G$ can be made into a representation
  of $A$ with corresponding $A$-connection defined by
  \[ \nabla_{\alpha}s(x)= \left.\frac{\d}{\d t} \phi_{\alpha}^{t}(x)^{-1}
    s(\phi_{\rho(\alpha)}^{t}(x))\right|_{t=0}.\]
  Moreover, this construction defines a homomorphism of semi-rings
  \[ \Phi: \Rep(\G)\to \Rep(A),\]
  which is an isomorphism if $\G$ is $\s$-simply connected.
\end{corollary}

\begin{proof}  
  We first interpret the groupoid action of $\G$ on $E$ as a groupoid morphism
  $\F:\G\to\GL(E)$. Passing to algebroids, we obtain an algebroid
  homomorphism $F:A\to\mathfrak{gl}(E)$. By the exercise
  above, this is the same thing as a flat $A$-connection $\nabla$ on
  $E$: $\nabla_{\al}= D_{F(\al)}$.  The actual formula for $\nabla$
  follows from the fact that the flow of $F(\al)$ is obtained by applying
  $F$ to the flow of $\al$. Using the definition of the derivation
  induced by a section of $\mathfrak{gl}(E)$ (see the previous proof),
  we conclude that
  \[ 
  \nabla_{\al}s= \left.\frac{\d}{\d t}\right|_{t=0}(F\circ
  \phi_{\al}^{t})^{*}(s),
  \]
  which is a compact version of the formula in the statement.
\end{proof}

A notion which actually generalizes both notions from geometry and Lie
theory, is the notion of \textbf{Lie algebroid cohomology with
coefficients} in a representation $E$. Given a representation $E$,
defined by a flat $A$-connection $\nabla$, we introduce the de Rham
complex of $A$ with coefficients in $E$ as follows. A $p$-differential
$A$-form with values in $E$ is an alternating $C^{\infty}(M)$-multilinear map 
\[ 
\omega:
\underbrace{\Gamma(A)\times\cdots\times\Gamma(A)}_{p\ \text{times}}
\to \Gamma(E).
\]
We denote by 
\[ \Omega^{p}(A;E)= \Gamma(\wedge^p A^*\otimes E)\]
the set of all $p$-differential $A$-forms with values in $E$. The
differential 
\[ \d: \Omega^{p}(A;E)\to \Omega^{p+1}(A;E)\]
is defined by the usual Koszul-type formula
\begin{multline}
\label{eq:differential}
\d \omega(\al_0,\dots,\al_p)=\sum_{k=0}^{p}
(-1)^k\nabla_{\al_k}(\omega(\al_0,\dots,\widehat{\al}_k,\dots,\al_p))\\
+\sum_{k<l}(-1)^{k+l+1}\omega([\al_k,\al_l],\al_0,\dots,\widehat{\al}_k,\dots,\widehat{\al}_l,\dots,\al_p),
\end{multline}
where $\al_0,\dots,\al_p\in\Gamma(A)$. 
You should check that $\d^2=0$. The resulting cohomology is denoted
$H^\bullet(A;E)$ and is called the \textbf{Lie algebroid cohomology with
coefficients} in the representation $E$. In the particular 
case of trivial coefficients, i.e., where $E$ is the trivial line
bundle over $M$ with $\nabla_{\alpha}= \Lie_{\rho(\alpha)}$, we talk
about the de Rham complex $\Omega^\bullet(A)$ of $A$ and the Lie
algebroid cohomology $H^\bullet(A)$ of $A$.

\begin{exercise} 
\label{ex:cohomology}
The notion of Lie algebroid cohomology generalizes
many well-known cohomology theories. In particular:
\begin{enumerate}[(a)]
\item When $A=TM$, check that $H(TM;E)=H_{\DR}(M;E)$ is the usual de
 Rham cohomology;
\item When $A=\gg$ is a Lie algebra, check that $H(A;E)=H(\gg;E)$
  is the usual Chevalley-Eilenberg Lie algebra cohomology;
\item When $A=T\F$ is the Lie algebroid of a foliation, check that 
  $H(A;E)=H(\F;E)$ is the usual foliated cohomology.
\end{enumerate}
\end{exercise}

The last example in the previous exercise shows that Lie algebroid
cohomology is, in general, infinite dimensional, and quite hard to
compute. With the exception of transitive Lie algebroids, there is no
known effective method to compute Lie algebroid cohomology.

Just like in the case of Lie groups, we can relate the Lie algebroid
cohomology of a Lie algebroid $A$, associated with a Lie groupoid
$\G$, with the de Rham cohomology of invariant forms on the groupoid.

\begin{exercise}
Let $\G$ be an $s$-connected Lie groupoid with Lie algebroid
$A$. Define a right-invariant differential form on $\G$ to be a $s$-foliated
differential form $\omega$ which satisfies:
\[ R_g^*\omega=\omega, \quad \forall g\in\G.\]
We denote by $\Omega^\bullet_{\text{inv}}(\G)$ the space of right
invariant forms on $\G$.
\begin{enumerate}[(a)]
\item Show that there is a natural isomorphism between
  $\Omega^\bullet_{\text{inv}}(\G)$ and the de Rham complex
  $\Omega^\bullet(A)$ of $A$. Conclude that the Lie algebroid
  cohomology of $A$ is isomorphic to the invariant de Rham cohomology
  of $\G$;
\item Generalize this isomorphism to any coefficients.
\end{enumerate}
\end{exercise}

\section{Notes}
The concept of Lie algebroid was introduced by Jean Pradines in 1966--68 
who, in a serious of notes \cite{Pra1,Pra2,Pra3,Pra4}, introduced a Lie theory
for Lie groupoids. Related purely algebraic notions, in particular the
notion of a Lie pseudo-algebra, were introduced much before by many authors,
and under different names (see the historical notes in \cite{Mack2}
and \cite[Chapter 3.8]{Mack3}).

The theory of Lie algebroids only took off in the late 1980's with the 
works of Almeida and Molino \cite{Mol} on developable foliations and the 
works of Mackenzie on connection theory (see the account in his first book \cite{Mack1}). 
These works were devoted almost exclusively to transitive Lie algebroids, and it 
was Weinstein \cite{Wein1} and Karas\"ev \cite{Kar}, who understood first the 
need to study non-transitive algebroids, namely from the appearence of the 
cotagent Lie algebroid of a Poisson manifold (\cite{CoDaWe}). 

The notion of morphism between two Lie algebroids is first appropriately understood in
the work of Higgins and Mackenzie \cite{HiMa}, where one can also find the description
in terms of connections. An alternative treatment using
the supermanifold formalism is due to Va\u{\i}ntrob \cite{Vain}. The theory of connections
played a large motivation in Mackenzie's approach to Lie groupoid and algebroid theory. 
A geometric approach to the theory of connections on Lie algebroids was given by 
Fernandes in \cite{Fer1,Fer2}. 

Representations of Lie algebroids where introduced first for transitive Lie algebroids by 
Mackenzie \cite{Mack1}, and they appear in various different contexts such us in the 
study of cohomological invariants attached to Lie algebroids (see, e.g., 
\cite{ELW,Cr,Fer2}).
The cohomology theory of Lie algebroids was started by Mackenzie \cite{Mack1} (transitive
case) inspired by the cohomology theory of Lie pseudoalgebras. Since then, many authors 
have developed many aspects of the theory (see, e.g., \cite{Ginz,Hueb1,WeXu,Xu2}) though 
computations are in general quite difficult.

The connection between Lie algebroids and Cartan's equivalence method is not well-known, 
and deserves to be explored. We have learned it from Bryant \cite[Appendix A]{Bry}.

\lecture{Integrability: Topological Theory}%
\label{integrability:I}            %
%
\section{What is integrability?}
In the previous lecture we saw many examples of Lie algebroids. Some
of these were \emph{integrable} Lie algebroids, i.e., they were isomorphic
to the Lie algebroid of some Lie groupoid. We are naturally
led to ask:
\begin{itemize}
\item Is every Lie algebroid integrable? 
\end{itemize}
It may come as a surprise to you, based on your experience with Lie
algebras or with the Frobenius integrability theorem, that the answer
is \emph{no}. This is a subtle (but important) phenomenon, which was
overlooked for sometime (see the notes at the end of this lecture). 

Let us start by giving a few examples of non-integrable Lie
algebroids, though at this point we cannot fully justify why they are
not integrable (this will be clear later on).

Our first source of examples comes from the Lie algebroids associated with
closed 2-forms.

\begin{example}
  \label{ex:2-forms:cont}
  Recall from the previous lecture (see Example \ref{ex:2-forms}) that
  every closed 2-form 
  $\omega\in\Omega^2(M)$ determines a Lie algebroid structure on 
  $A_\omega=TM\oplus\L$. If $M$ is simply-connected and the cohomology
  class of $\omega$ is integral, we have a prequantization principal
  $\Ss^1$-bundle which gives rise to a Lie groupoid $\G_\omega$
  integrating $A_\omega$. The groupoid $\G_\omega$ being
  transitive, all its isotropy Lie groups are isomorphic, and one canl
  see that they are in fact isomorphic to $\Rr/\Gamma_\omega$, where
  $\Gamma_\omega$ is the group of spherical periods of $\omega$:
    \[
    \Gamma_\omega=\set{\int_\gamma\omega:\gamma\in\pi_2(M)}\subset\Rr.
    \]

  When $\omega$ is not integral, we don't have a prequantization
  bundle, and the construction of $\G_\omega$ fails. However, we will
  see later in this lecture that, if $A_\omega$ integrates to a source
  simply-connected Lie groupoid $\G_\omega$, then its isotropy groups
  must still be isomorphic to $\Rr/\Gamma_\omega$. One can then show
  that $A_\omega$ is integrable iff the group of spherical periods
  $\Gamma_\omega\subset\Rr$ is a discrete subgroup.
  
  Let us take, for example, $M=\Ss^2\times\Ss^2$ with
  $\omega=dS\oplus\lambda\, dS$, where $dS$ is the standard area form on
  $\Ss^2$. If we choose $\lambda\in\Rr-\Qq$, then the group of
  spherical periods is $\Gamma_\omega=\Zz\oplus\lambda\Zz$, so that
  $\Rr/\Gamma_\omega$ is non-discrete. Therefore, the corresponding 
  Lie algebroid is non-integrable.
\end{example}

Another source of examples of non-integrable Lie algebroids is Poisson
geometry, since there are Poisson manifolds whose cotangent
Lie algebroids (see Example \ref{ex:Poisson}) are not integrable.

Let us recall that the dual $\gg^*$ of a finite dimensional Lie algebra has a
natural linear Poisson structure, namely the Kostant-Kirillov-Souriau
Poisson structure, which is defined by:
\[ 
  \{f_1,f_2\}(\xi)=\langle \xi,[\d_{\xi}f_1,\d_{\xi}f_2]\rangle,
  \quad f_1,f_2\in C^\infty(\gg^*),\ \xi\in\gg^*.
\]
The corresponding cotangent Lie algebroid integrates to the Lie
groupoid $\G=T^*G$, where $G$ is a Lie group with Lie algebra
$\gg$. The source and target maps $\s,\t:T^*G\to\gg^*$ are the
right and left translation maps trivializing the cotangent bundle:
\[ 
\s(\omega_g)=(\d_e R_g)^*\omega_g,\quad \t(\omega_g)=(\d_e L_g)^*\omega_g,
\]
while composition is given by
\[
\omega_g\cdot\eta_h=(\d_{hg}R_{g^-1})^*\eta_h.
\]

\begin{exercise}
  Consider the trivialization $T^*G=G\times \gg^*$ obtained by right
  translations. Determine the new expressions of $\s$, $\t$ and the
  product under this isomorphism, and check that the orbits of $\G$
  are precisely the coadjoint orbits. Note that these coincide with
  the symplectic leaves of $\gg^*$.
\end{exercise}

Though linear Poisson structures are always integrable, one can
easily produce non-integrable Poisson structures by slightly modifying
them. Here are two examples.

\begin{example}
  \label{ex:su(n)}
  Let us take $\gg=\su(n)$, with $n>2$. We use the $\Ad$-invariant inner
  product on $\gg$ defined by:
  \[ \langle X,Y\rangle=\Tr(XY^*),\]
  to identify $\gg^*\simeq \gg$. Observe that the spheres
  $\Ss_r=\set{X\in\gg:||X||=r}$ are collections of (co)adjoint orbits,
  and hence are Poisson submanifolds of $\gg^*$, with the KKS Poisson
  structure. As we will see later, it turns out that (the cotangent Lie
  algebroids of) the Poisson manifolds $\Ss_r$ are not integrable.
\end{example}

In the case $n=2$, the spheres are symplectic submanifolds and so they
turn out to be integrable. However, we can still change slightly the
Poisson bracket and obtain another kind of non-integrable Poisson structure.

\begin{example}
  \label{ex:su(2)}
  Let us take then $\gg=\su(2)$. If we identify $\su(2)^*\simeq\Rr^3$,
  the Poisson bracket in euclidean coordinates $(x,y,z)$ is determined
  by:
  \[ \{x,y\}=z,\quad \{y,z\}=x,\quad\{z,x\}=y.\]
  The symplectic leaves of this Poisson structure (the coadjoint
  orbits) are the spheres $r^2=x^2+y^2+z^2=C$ and the origin. In fact,
  any smooth function $f=f(r)$ is a Casimir (Poisson commutes with any
  other function).
  
  Let us choose a smooth function of the radius $a(r)$, such that
  $a(r)>0$ for $r>0$. Then we can define a new rescaled Poisson
  bracket:
  \[ \{f,g\}_a:=a\{f,g\}.\]
  This bracket clearly has the same symplectic foliation, while the
  symplectic area of each leaf $x^2+y^2+z^2=r^2$ is rescaled by a
  factor of $1/a(r)$. We will see later that this new Poisson
  structure is integrable iff the symplectic area function $A(r)= 4\pi
  r/a(r)$ has no critical points. For example, if $a(r)=e^{r^2/2}$
  then $A(r)$ has a critical point at $r=1$, and the corresponding
  Poisson manifold is not integrable.
\end{example}

\vskip 20 pt

As the previous examples illustrate, it is not at all obvious when
a Lie algebroid is integrable, and even ``reduction'' may turn an integrable
Lie algebroid to a non-integrable one. Hence, it is quite important to
find the answer to the following question:
\begin{itemize}
\item What are the precise obstructions to integrate a Lie algebroid?
\end{itemize}
This is the problem we shall address in the remainder of this lecture
and in the next lecture.

\section{Integrating Lie algebras}

Before we look at the general integrability problem for Lie
algebroids, it is worth to consider the special case of Lie
algebras:

\begin{theorem}[Lie III]
Every finite dimensional Lie algebra is isomorphic to the Lie algebra
of some Lie group.
\end{theorem}

In most texts in Lie theory, the proof of Lie's third theorem uses the
structure theory of Lie algebras. The usual strategy is to prove first Ado's
Theorem, stating that every Lie algebra has a faithful finite
dimensional representation, from which it follows that there exists a
matrix Lie group (not necessarily simply connected) with the given Lie
algebra. Ado's Theorem, in turn, can be easily proved for semi-simple
Lie algebras, follows from some simple structure theory for solvable
Lie algebras, and then extends to any Lie algebra using the Levi
decomposition. There is, however, a much more direct geometric
approach to Lie's third theorem, which we will sketch now, and which
will be extended to Lie algebroids in later sections.

The main idea is as follows: Suppose $\gg$ is a finite Lie algebra
which integrates to a connected Lie group $G$. Denote by $P(G)$ the space of
paths in $G$ starting at the identity $e\in G$, with the
$C^2$-topology:
\[ P(G)=\set{g: [0,1]\to G|~ g\in C^2,\ g(0)=e}.\]
Also, denote by $\sim$ the equivalence relation defined by
$C^1$-homotopies in $P(G)$ with fixed end-points. Then we have a
standard description of the simply-connected Lie group integrating
$\gg$ as
\[ \widetilde{G}= P(G)/\sim.\]
Let us be more specific about the product in $\widetilde{G}$: given
two paths $g_1, g_2\in P(G)$ we define
\[
g_1\cdot g_2(t)\equiv \left\{
\begin{array}{ll}
g_2(2t),\qquad& 0\le t\le \frac{1}{2},\\ \\
g_1(2t-1)g_2(1),\qquad & \frac{1}{2}< t\le 1.
\end{array}
\right.
\]
Note, however, that this multiplication can take one out of $P(G)$,
since the composition will be a path which is only \emph{piecewise}
$C^1$, so we need:

\begin{exercise}
  Show that any element in $P(G)$ is equivalent to some $g(t)$ with
  derivatives vanishing at the end-points, and if $g_1$ and $g_2$ have this
  property, then $g_1\cdot g_2\in P(G)$.
\end{exercise}

This will give us a multiplication in $P(G)$, which is associative up
to homotopy, so we get the desired multiplication on the quotient
space which makes $\widetilde{G}$ into a (topological) group. Since
$\widetilde{G}$ is the universal covering space of the manifold $G$,
there is also a smooth structure in $\widetilde{G}$ which makes it
into a Lie group. In this way, we have recovered the simply-connected
Lie group integrating $\gg$ \emph{assuming} that $\gg$ is integrable.

Now, any $G$-path $g$ defines a path $a:I\to\gg$ by differentiation
and right translations:
\begin{equation}
  \label{eq:G-path:gg-path}
  a(t)=\left.\frac{\d}{\d s}g(s)g(t)^{-1}\right|_{s=t}.
\end{equation}
Let us denoted by $P(\gg)$ the space of paths $a:[0,1]\to \gg$ with
the $C^1$-topology.

\begin{exercise}
  Show that the map $P(G)\to P(\gg)$ just defined is a homeomorphism.
\end{exercise}

Using this bijection, we can transport our equivalence relation  and
our product in $P(G)$ to an equivalence relation and a product in
$P(\gg)$. The next lemma gives an explicit expression for the
equivalence relation in $P(\gg)$:

\begin{lemma}
  Two paths $a_0,a_1\in P(\gg)$ are equivalent iff there exists a
  homotopy $a_\epsilon\in P(\gg)$, $\epsilon\in[0,1]$, joining $a_0$ to
  $a_1$, such that 
  \begin{equation}
    \label{variation:pre}
    \int_{0}^{1}B_{\epsilon}(s)\cdot\frac{\d a_\epsilon}{\d\epsilon}(s)\, ds=0,
    \quad \forall \epsilon\in[0,1],
  \end{equation}
  where $B_\epsilon(t)\in GL(\gg)$, for each $\epsilon$, is the
  solution of the initial value problem:
  \[ 
  \left\{
    \begin{array}{l}
      \left.\frac{\d}{\d s} B_\epsilon(s)B_\epsilon(t)^{-1}\right|_{s=t}
      =\ad(a_\epsilon(t)),\\
      \\
      B_\epsilon(0)=I.\\
    \end{array}
  \right.
  \]
\end{lemma}

\begin{exercise}
  Prove this lemma.\\
  (Hint: see the proof of Proposition \ref{equivalence} below, which
  generalizes this to any Lie algebroid.)
\end{exercise}

\begin{exercise}
Let $g_1(t)$ and $g_2(t)$ be $G$-paths and denote by $a_1(t)$ and
$a_2(t)$ the corresponding $\gg$-paths defined by
(\ref{eq:G-path:gg-path}). Show that $\gg$-path associated with the
concatenation $g_1\cdot g_2$ is given by:
\begin{equation}
  \label{eq:concatenation:gg-paths}
  a_1\cdot a_2(t)\equiv \left\{
    \begin{array}{ll}
      2a_2(2t),\qquad& 0\le t\le \frac{1}{2},\\ \\
      2a_1(2t-1),\qquad & \frac{1}{2}< t\le 1.
    \end{array}
  \right.
\end{equation}
\end{exercise}

The Lemma and the previous exercise show that the equivalence relation
and the product in $P(\gg)$ can be defined exclusively in terms of
data in $\gg$, and involve no reference at all to the Lie group
$G$. Hence, we can make a fresh start, \emph{without} assuming that
$\gg$ is the Lie algebra of some $G$:

\begin{definition}
  Let $\gg$ be a finite dimensional Lie algebra. Then we define the
  topological group
  \[ \G(\gg)=P(\gg)/\sim,\]
  where:
  \begin{enumerate}[(i)]
  \item Two paths $a_0,a_1\in P(\gg)$ are equivalent iff there exists a
    homotopy $a_\epsilon\in P(\gg)$, $\epsilon\in[0,1]$, joining $a_0$ to
    $a_1$, such that relation (\ref{variation:pre}) is satisfied.
  \item The product is $[a_1][a_2]=[a_1\cdot a_2]$, where the dot
  denotes composition of paths, which is defined by
  (\ref{eq:concatenation:gg-paths}).
  \item The topology on $\G(\gg)$ is the quotient topology.  
  \end{enumerate}
\end{definition}

From what we saw above, if $\gg$ is integrable, then $\G(\gg)$ is the
unique simply-connected Lie group with Lie algebra $\gg$. Of course,
we know that all finite dimensional Lie algebras are integrable. But
if we want to prove this, what is left is to prove is to show that the
smooth Banach structure on $P(\gg)$ descends to a smooth structure on
the quotient:

\begin{theorem}
  Let $\gg$ be a finite dimensional Lie algebra. Then $\G(\gg)$ is
  a simply-connected Lie group integrating $\gg$.
\end{theorem}

We will not give a proof here, since we will eventually prove a much
more general version of this result, valid for Lie algebroids.

\begin{exercise}
  Use the description of $\G(\gg)$ to prove Lie's 2nd theorem for Lie
  algebras: if $\phi:\gg\to\hh$ is a Lie algebra homomorphism, then
  there exists a unique Lie group homomorphism $\Phi:\G(\gg)\to H$,
  with $\d\Phi=\phi$, for any Lie group $H$ integrating $\hh$.
\end{exercise}

\section{Integrating Lie algebroids}

Motivated by what we just did for Lie algebras/Lie groups we set:

\begin{definition}
  Let $\G$ be a Lie groupoid. A \textbf{$\G$-path} is a path $g:
  [0,1]\to \G$ such that $\s(g(t))=x$, for all $t$, and $g(0)= 1_{x}$
  (i.e., a path lying in a $\s$-fiber of $\G$ and starting at the unit
  of the fiber). We denote by $P(\G)$ the space of $\G$-paths,
  furnished with the $C^2$-topology.
\end{definition}

Let $A$ be a Lie algebroid which is integrable to a Lie groupoid $\G$.
Let us denote by $\sim$ the equivalence relation defined by
$C^1$-homotopies in $P(\G)$ with fixed end-points. Then we can
describe the $s$-simply connected Lie groupoid integrating $A$:
\[ \widetilde{\G}=P(\G)/\sim. \]
The source and target maps are the obvious ones(\footnote{Note that
  the symbols $\s$ and $\t$ have different meanings on different sides
  of these relations!}):
\[ \s([g])=\s(g(0)),\quad \t([g])=\t(g(1)),\]
and for two paths $g_1,g_2\in P(\G)$ which are composable,
i.e., such that $\s(g_1(1))=\t(g_2(0))$, we define 
\[
g_1\cdot g_2(t)\equiv \left\{
\begin{array}{ll}
g_2(2t),\qquad& 0\le t\le \frac{1}{2},\\ \\
g_1(2t-1)g_2(1),\qquad & \frac{1}{2}< t\le 1.
\end{array}
\right.
\]
Note that we are just following the same strategy as we did for the
integrability of Lie algebras. Of course, the same comments about
concatenation taking us out of $P(\G)$ apply, but:
\begin{enumerate}[(a)]
\item any element in $P(\G)$ is equivalent to some $g(t)$ with
  derivatives vanishing at the end-points, and
\item if $g_1$ and $g_2$ have this property, then $g_1\cdot g_2\in
  P(\G)$. 
\end{enumerate}
Therefore, the multiplication will be well defined and associative up
to homotopy, so we get the desired multiplication on the quotient
space which makes $\widetilde{\G}$ into a (topological) groupoid. The
construction of the smooth structure on $\widetilde{\G}$ is similar to
the construction of the smooth structure on the universal cover of a
manifold (see, also, the proof of Theorem \ref{thm:1:connected}).

Now, just like in the case of Lie algebras, we have:

\begin{proposition}
  \label{A-G-paths}
  If $\G$ integrates the Lie algebroid $A$, then there is a
  homeomorphism $D^{R}: P(\G)\to P(A)$ between the space of $\G$-paths,
  and the space of $A$-paths. 
\end{proposition}

The map $D^{R}$ will be called the \textbf{differentiation of $\G$-paths}, and
its inverse will be called the \textbf{integration of $A$-paths}.

\begin{proof}
  Any $\G$-path $g:I\to \G$ defines an $A$-path $D^{R}(g):I\to A$,
  where
  \begin{equation} 
  \label{D-to-R}
  (D^{R}g)(t)= (\d R_{g(t)^{-1}})_{g(t)} \dot{g}(t) \ ,
  \end{equation}
  (here, for $h: x\to y$ an arrow in $\G$, $R_{h}: \s^{-1}(y)\to s^{-1}(x)$
  is the right multiplication by $h$).  Conversely, any $A$-path $a$
  arises in this way: we first integrate (using \textsc{Lie II}) the Lie
  algebroid morphism $TI\to A$ defined by $a$, and then we notice that
  any Lie groupoid homomorphism $\phi:I\times I\to \G$, from the pair
  groupoid into $\G$, is of the form $\phi(s,t)=g(s)g^{-1}(t)$, for some
  $\G$-path $g$.
\end{proof}

Our next task is to transport the equivalence relation from
$\G$-paths to $A$-paths. For that we will need the concept of \emph{flow
of a section} of a Lie algebroid: if $\al\in\Gamma(A)$ then its flow
is the unique 1-parameter group of Lie algebroid automorphisms
$\phi_\al^t:A\to A$ such that:

\begin{equation}
  \label{Lie-flows}
  \left.\frac{\d}{\d t}\right|_{t=0}\phi_\al^t(\be)=[\al,\be].
\end{equation}

\begin{exercise}
  Give an explicit construction of the flow $\phi_\al^t:A\to A$ of a
  section $\al\in\Gamma(A)$, and show that they are Lie algebroid
  automorphisms that cover the flow $\phi_{\rho(\al)}^t:M\to M$ of the
  vector field $\rho(\al)$.\\
  (Hint: Consider the derivation $\ad_\al\in\Der(A)$ defined by
  $\ad_\al(\be)=[\al,\be]$.
  Identifying $\ad_\al$ with a section of $\mathfrak{gl}(A)$, take
  $\phi_\al^t$ to be the flow of $\ad_\al$ in the sense of the last
  lecture.)
\end{exercise}

\begin{exercise}
  Define the flow $\phi_{\al_s}^{t,s}$ of a \emph{time-dependent}
  section $\al_s\in\Gamma(A)$.
\end{exercise}

Henceforth, by a \emph{variation of $A$-paths} we mean a map
\[a_{\epsilon}(t)= a(\epsilon, t): I\times I\to A\]
such that $a_{\epsilon}$ is a family of $A$-paths of class $C^2$ on
$\epsilon$, with the property that the base paths
$\gamma_{\epsilon}(t)= \gamma (\epsilon, t): I\times I\to M$ have
fixed end points.  When $a_{\epsilon}=D^R g_{\epsilon}$, the family
$g_{\epsilon}$ does not necessarily give a homotopy between $g_{0}$
and $g_{1}$, because the end points $g_{\epsilon}(1)$ may depend on
$\epsilon$. The following propositions describes two distinct ways of
controlling the variation $\frac{\d}{\d\epsilon} g_{\epsilon}(1)$, both
depending only on infinitesimal data (i.e., Lie algebroid data).

\begin{proposition}
  \label{equivalence}
  Let $A$ be any Lie algebroid and $a=a_{\epsilon}\in P(A)$ a variation of
  $A$-paths.
  \begin{enumerate}[(i)]
  \item If $\nabla$ is an $TM$-connection on $A$, the solution
    $b=b(\epsilon, t)$ of the differential equation
    \begin{equation}\label{diffeq}
      \partial_{t}b-\partial_{\epsilon} a= T_{\nabla}(a, b),\quad
      b(\epsilon,0)=0,
    \end{equation}
    does not depend on $\nabla$. Moreover, $\rho(b)=
    \frac{\d}{\d\epsilon}\gamma$.
  \item If $\xi_{\epsilon}$ are time depending sections of $A$ such
    that $\xi_{\epsilon}(t, \gamma_{\epsilon}(t))= a_{\epsilon}(t)$,
    then $b(\eps,t)$ is given by
    \begin{equation}
      \label{variation}
      b(\epsilon, t) = \int_{0}^{t}\phi_{\xi_{\epsilon}}^{t, s}
      \frac{\d \xi_{\epsilon}}{\d\epsilon} (s, \gamma_{\epsilon}(s))\d s,
    \end{equation}
    where $\phi_{\xi_{\epsilon}}^{t, s}$ denotes the flow of the
    time-dependent section $\xi_{\epsilon}$.
  \item If $\G$ integrates $A$ and $g_{\epsilon}$ are the $\G$-paths
    satisfying $D^{R}(g_{\epsilon})= a_{\epsilon}$, then $b=
    D^{R}(g^{t})$, where $g^t$ are the paths in $\G$: $\eps\to
    g^{t}(\epsilon)=g(\epsilon, t)$.
  \end{enumerate}
\end{proposition}

This motivates the following definition:

\begin{definition}
  We say that two $A$-paths $a_0$ and $a_1$ are \textbf{$A$-homotopic}
  and we write $a_0\sim a_1$, if there exists a variation
  $a_{\epsilon}$ with the property that $b$ insured by Proposition
  \ref{equivalence} satisfies $b(\epsilon, 1)= 0$ for all $\epsilon
  \in I$.
\end{definition}

If $A$ admits an integration $\G$, then the isomorphism $D^{R}:
P(\G)\to P(A)$ transforms the usual homotopy into $A$-homotopy. Also,
since $A$-paths should be viewed as algebroid morphisms, the
pair $(a,b)$ defining the $A$-homotopy should be viewed as a true
homotopy
\[ a \d t + b \d\epsilon  : TI\times TI \to A.\]
in the world of Lie algebroids.

\begin{exercise}
  Show that equation (\ref{diffeq}) is equivalent to the condition
  that $a\d t+b\d\epsilon:TI\times TI\to A$ is a morphism of Lie
  algebroids.\\
  (Hint: See Exercise \ref{morphisms:algebroids}.)
\end{exercise}

\begin{proof}[Proof of Proposition \ref{equivalence}] 
  Obviously, (i) follows from (ii). To prove (ii), let
  $\xi_{\epsilon}$ be as in the statement, and let $\eta$ be given by
  \[ \eta (\epsilon, t, x)= \int_{0}^{t}
  \phi_{\xi_{\epsilon}}^{t, s} \frac{\d\xi_{\epsilon}}{\d\epsilon}(s,
  \phi_{\rho(\xi_{\epsilon})}^{s, t}(x))\d s \in A_x. \] 
  We may assume that $\xi_{\epsilon}$ has compact support. We note
  that $\eta$ coincides with the solution of the equation
  \begin{equation}
    \label{eq.equiv}
    \frac{\d\eta}{\d t}- \frac{\d\xi}{\d\epsilon}= [\eta, \xi]\ .
  \end{equation}
  with $\eta(\epsilon, 0)= 0$. Indeed, since
  \[ \eta (\epsilon, t, -)= \int_{0}^{t}
  (\phi_{\xi_{\epsilon}}^{s,
  t})^{*}(\frac{\d\xi_{\epsilon}^{s}}{\d\epsilon})\d s \in \Gamma(A),\]
  equation (\ref{eq.equiv}) immediately follows from the basic formula
  (\ref{Lie-flows}) for flows. Also, $X=\rho(\xi)$ and $Y=\rho(\eta)$ satisfy
  a similar equation on $M$, and since we have
  $X(\epsilon,t,\gamma_{\epsilon}(t))= \frac{\d\gamma}{\d t}$, it follows
  that $Y(\epsilon, t, \gamma_{\epsilon}(t))=
  \frac{\d\gamma}{\d\epsilon}$. In other words, $b(\epsilon, t)=
  \eta(\epsilon, t, \gamma(\epsilon, t))$ satisfies $\rho(b)=
  \frac{\d\gamma}{\d\epsilon}$. We now have
  \begin{equation*}
    \partial_{t}b= \nabla_{\frac{\d\gamma}{\d t}}\eta+
    \frac{\d\eta}{\d t}=
    \nabla_{\rho(\xi)}\eta+ \frac{\d\eta}{\d t}.
  \end{equation*}
  Subtracting from this the similar formula for $\partial_{\epsilon}a$
  and using (\ref{eq.equiv}) we arrive at
  \begin{equation*} 
    \partial_{t}b-\partial_{\epsilon}a=
    \nabla_{\rho(\xi)}\eta-\nabla_{\rho(\eta)}\xi+[\eta,\xi]= 
    T_{\nabla}(\xi,\eta).
  \end{equation*}

  We are now left with proving (iii). Assume that $\G$ integrates $A$ and
  $g_{\epsilon}$ are the $\G$-paths satisfying $D^{R}(g_{\epsilon})=
  a_{\epsilon}$. The formula of variation of parameters applied to the
  right-invariant vector field $\xi_{\epsilon}$ shows that
  \begin{align*}
    \frac{\partial g(\eps,t)}{\partial \eps}&= \int_0^t
    (\d\varphi^{t,s}_{\xi_{\epsilon}})_{g(\eps,s)}
    \frac{\d\xi_{\epsilon}^s}{\d\eps}(g(\eps,s))\d s\\ &=
    (\d R_{g(\eps,t)})_{\gamma_\eps(t)}
    \int_0^t\phi^{t,s}_{\xi_{\epsilon}}
    \frac{\d\xi_{\epsilon}^s}{\d\eps}(\gamma_\eps(s))\d s.
  \end{align*}
  But then:
  \[
  D^R(g^t)=\int_0^t\phi^{t,s}_{\xi_{\epsilon}}
  \frac{\d\xi_{\epsilon}^s}{\d\eps}(\gamma_\eps(s))\d s=b(\eps,t).
  \]
\end{proof}

Let us now turn to the task of transporting the composition from
$\G$-paths to $A$-paths. This is an easy task if we use the properties
of $A$-homotopies given in the following exercise:

\begin{exercise}
  \label{homotopy properties}
  Let $A$ be a Lie algebroid. Show that:
  \begin{enumerate}[(a)]
  \item If $\tau: I\to I$ is a smooth change of parameter, then any
    $A$-path $a$ is $A$-homotopic to its re-parameterization
    $a^{\tau}(t)\equiv\tau'(t)a(\tau(t))$.
  \item Any $A$-path $a_{0}$ is $A$-homotopic to a smooth (i.e., of
    class $C^{\infty}$) $A$-path.
  \item If two smooth $A$-paths $a_{0}$ and $a_{1}$ are $A$-homotopic,
    then there exists a smooth $A$-homotopy between them.
  \end{enumerate}
\end{exercise}

For a Lie algebroid $\pi:A\to M$ we say that two $A$-paths $a_0$ and
$a_1$ are composable if they have the same end-points, i.e.,
$\pi(a_{0}(1))= \pi(a_{1}(0))$. In this case, we define their
concatenation by
\[
a_{1}\odot a_{0}(t)\equiv \left\{
\begin{array}{ll}
2a_{0}(2t),\qquad& 0\le t\le \frac{1}{2},\\ \\
2a_{1}(2t-1),\qquad & \frac{1}{2}< t\le 1.
\end{array}
\right.
\]
This is essentially the multiplication that we need. However,
$a_{1}\odot a_{0}$ is only \emph{piecewise} smooth. One way around
this difficulty is allowing for $A$-paths which are \emph{piecewise}
smooth. Instead, we choose a cutoff function $\tau\in C^\infty(\Rr)$
with the following properties:
\begin{enumerate}[(a)]
\item $\tau(t)=1$ for $t\ge 1$ and $\tau(t)=0$ for $t\le 0$;
\item $\tau'(t)>0$ for $t\in]0,1[$.
\end{enumerate}
For an $A$-path $a$ we denote by $a^{\tau}$ its re-parameterization
$a^{\tau}(t):=\tau'(t)a(\tau(t))$. We now define the multiplication of
composable $A$-paths by
\[  a_{1}\cdot a_{0}\equiv a_{1}^{\tau}\odot a_{0}^{\tau}\in P(A).\]
According to Exercise \ref{homotopy properties} (a), the product
$a_{0}\cdot a_{1}$ is equivalent to $a_{0}\odot a_{1}$ whenever $a_{0}(1)=
a_{1}(0)$.  It follows that the quotient $\G(A)$ is a groupoid with
this product together with the natural structure maps: 
\begin{itemize}
\item the source and target $\s,\t:\G(A)\to M$ map a class $[a]$ to
  its end-points $\pi(a(0))$ and $\pi(a(1))$, respectively; 
\item the unit section $u:M\to P(A)$ maps $x$ to the class $[0_x]$ of
  the constant trivial path above $x$; 
\item the inverse $i:P(A)\to P(A)$ maps a class $[a]$ to the class
  $[\overline{a}]$ of its opposite path, which is defined by
  $\overline{a}(t)=-a(1-t)$.
\end{itemize}
The groupoid $\G(A)$ will be called the \textbf{Weinstein groupoid} of
the Lie algebroid $A$.

\begin{theorem}
  \label{Weinstein}
  Let $A$ be a Lie algebroid over $M$. Then the quotient
  \[ \G(A)\equiv P(A)/\sim \]
  is an $\s$-simply connected topological groupoid independent of the
  choice of cutoff function. Moreover, whenever $A$ is integrable,
  $\G(A)$ admits a smooth structure which makes it into the unique
  $\s$-simply connected Lie groupoid integrating $A$.
\end{theorem}

\begin{proof}
  If we take the maps on the quotient induced from the structure maps
  defined above, then $\G(A)$ is clearly a groupoid. Note that the
  multiplication on $P(A)$ was defined so that, whenever $\G$
  integrates $A$, the map $D^{R}$ of Proposition \ref{A-G-paths}
  preserves multiplications. Hence the only thing we still have to
  prove is that $\s,\t: \G(A)\to M$ are open maps.
  
  {To} prove this we show that for any two $A$-homotopic $A$-paths
  $a_{0}$ and $a_{1}$, there exists a homeomorphism $T: P(A)\to P(A)$
  such that $T(a)\sim a$ for all $a$'s, and $T(a_{0})=a_{1}$. We can
  construct such a $T$ as follows: we let $\eta=\eta(\epsilon,t)$ be a
  family of time dependent sections of $A$ which determines the
  equivalence $a_{0}\sim a_{1}$ (see Proposition \ref{equivalence}),
  so that $\eta(\epsilon, 0)=\eta(\epsilon,1)= 0$ (we may assume
  $\eta$ has compact support, so that all the flows involved are
  everywhere defined). Given an $A$-path $b_0$, we consider a time
  dependent section $\xi_0$ so that $\xi_0(t,\gamma_0(t))= b(t)$ and
  denote by $\xi$ the solution of equation (\ref{eq.equiv}) with
  initial condition $\xi_{0}$.  If we set
  $\gamma_{\epsilon}(t)=\Phi_{\rho(\eta_{t})}^{\epsilon, 0}(\gamma_0(t))$
  and $b_{\epsilon}(t)= \xi_{\epsilon}(t, \gamma_{\epsilon}(t))$, then
  $T_{\eta}(b_0)\equiv b_1$ is homotopic to $b_0$ via $b_{\epsilon}$,
  and maps $a_{0}$ into $a_{1}$.
\end{proof}

The following exercises provides further evidence of the importance of
the notion of $A$-homotopy and the naturality of the construction of
the groupoid $\G(A)$:

\begin{exercise}
  \label{exer:hol:hom}
  Let $E$ be a representation of the Lie algebroid $A$. Show that if
  $a_0$ and $a_1$ are $A$-homotopic paths from $x$ to $y$, then the
  parallel transports $\tau_{a_0},\tau_{a_1}:E_x\to E_y$ along $a_0$
  and $a_1$ coincide.  Conclude that every representation
  $E\in\Rep(A)$ determines a representation of $\G(A)$, which
  in the integrable case is the induced smooth representation.
\end{exercise}

\begin{exercise}
  \label{exer:Lie:II}
  Show that every algebroid homomorphism $\phi:A\to B$ determines a
  continuous groupoid homomorphism $\Phi:\G(A)\to \G(B)$. If $A$ and
  $B$ are integrable, show that $\Phi$ is smooth and $\Phi_*=\phi$.
  Finally, use this to prove Lie II (Theorem \ref{thm:Lie:II}). 
\end{exercise}

\begin{exercise}
  \label{exer:exponential}
  Show that, for any Lie algebroid $A$, there exists an exponential
  map $\exp:\Gamma_c(A)\to\Gamma(\G(A))$, which generalizes the
  exponential map in the integrable case, and for any $\al,\be\in
  \Gamma_c(A)$ satisfies:
  \[ \exp(t\al)\exp(\be)\exp(-t\al)= \exp(\phi^t_\al \be) ,\]
  where $\phi^t_\al$ denotes the flow of $\al$.
\end{exercise}

\section{Monodromy}

In this section we will introduce the \emph{monodromy groups} that
control the integrability of a Lie algebroid $A$.

Let us assume first that $A$ is an integrable Lie algebroid and that
$\G$ is a source 1-connected Lie groupoid integrating $A$. If $x\in
M$, there are two Lie groups that integrate the isotropy Lie algebra
$\gg_x$: 
\begin{itemize}
\item The isotropy Lie groups $\G_x:=\s^{-1}(x)\cap\t^{-1}(x)$.
\item The 1-connected Lie group $G_x:=\G(\gg_x)$ with Lie algebra
  $\gg_x$.
\end{itemize}
Let us denote by $\G_x^0$ the connected component of $\G_x$ containing
the identity element. By simple ordinary Lie theory theory arguments,
there exists a subgroup $\tilde{\NN}_x\subset Z(G_x)$ of the center
of $G_x$ such that:
\[ \G_x^0\simeq G_x/\tilde{\NN}_x.\]
Note that $\tilde{\NN}_x$ can be identified with $\pi_1(\G_x^0)$ and
that it is a discrete subgroup of $G_x$. The group $\tilde{\NN}_x$
will be called the \emph{monodromy group} of $A$ at $x$. We will show
now that one can define these monodromy groups even when $A$ is
non-integrable, but then they may fail to be discrete. This lack of
discreteness is the clue to understand the non-integrability of $A$.

Let $\pi:A\to M$ be \emph{any} Lie algebroid. Notice that the isotropy
group of the Weinstein groupoid $\G(A)$ at $x$ is formed by the equivalence
classes of \textbf{$A$-loops} based at $x$:
\[ \G_x(A)=\set{[a]\in\G:\pi(a(0))=\pi(a(1))=x}.\]
We emphasize that \emph{only} the base path of $a$ is a loop and, in
general, we will have $a(0)\not=a(1)$. Also, the base loop must lie inside the
orbit $\O_x$ through the base point. 

It will be also convenient to consider the \textbf{restricted isotropy
group} formed by those $A$-loops whose base loop is contractible in the
orbit $\O_x$:
\[ \G_x(A)^0=\set{[a]\in\G_x(A):\gamma\sim\ast\text{ in }O_x}.\]
It is clear that $\G_x(A)^0$ is the connected component of the
identity of the isotropy group $\G_x(A)$. Moreover, the map
$[a]\mapsto[\gamma]$, associating to an $A$-homotopy class of
$A$-paths the homotopy class of its base path gives a short exact
sequence:
\[ 
\xymatrix{1\ar[r]& \G_x(A)^0\ar[r]& \G_x(A)\ar[r]& \pi_1(\O_x)\ar[r]& 1.}
\]

We know that when $A$ is integrable, $\G_x(A)^0$ is a connected Lie
group integrating the isotropy Lie algebra $\gg_x=\gg_x(A)$.  We would
like to understand this restricted isotropy Lie group for a general,
possibly non-integrable, Lie algebroid. Recall from Lecture
\ref{algebroids}, that over an orbit $\O_x$ of any Lie algebroid, the
anchor gives a short exact sequence of Lie algebroids:
\[
\xymatrix{0\ar[r]& \gg_{\O_x}\ar[r]& A|_{\O_x}\ar[r]^{\rho}& T\O_x\ar[r]& 0.}
\]
If we think of this sequence as a fibration, the next proposition
shows that there exists (the first terms of) an associated homotopy
long exact sequence:

\begin{proposition}
  \label{monodromy-seq}
  There exists a homomorphism $\partial: \pi_{2}(\O_x)\to \G(\gg_x)$
  which makes the following sequence exact:
  \[
  \xymatrix{\cdots \ar[r]& \pi_{2}(\O_x)\ar[r]^{\partial}& \G(\gg_x)
  \ar[r]& \G_x(A)\ar[r]& \pi_1(\O_x).}
  \] 
\end{proposition}

\begin{proof}
  {To} define $\partial$ let $[\gamma]\in \pi_{2}(\O_x)$ be
  represented by some smooth path $\gamma: I\times I\to \O_x$ which
  maps the boundary into $x$. We choose a morphism of algebroids
  \[ a\d t+ b\d\epsilon : TI\times TI\to A_{\O_x} \]
  (i.e., $(a,b)$ satisfies equation (\ref{diffeq})) which lifts
  $\d\gamma: TI\times TI\to T\O_x$ via the anchor, and such that
  $a(0,t)$, $b(\epsilon,0)$, and $b(\epsilon,1)$ vanish. This is
  always possible: we can take
  $b(\epsilon,t)=\sigma(\frac{\d}{\d\epsilon}\gamma(\epsilon,t))$ where
  $\sigma: T\O_x\to A_{\O_x}$ is any splitting of the anchor map, and
  take $a$ to be the unique solution of the differential equation
  (\ref{diffeq}) with the initial conditions $a(0,t)= 0$. Since
  $\gamma$ is constant on the boundary, $a_{1}= a(1,-)$ stays inside
  the Lie algebra $\gg_{x}$, i.e., defines a $\gg_{x}$-path
  $a_{1}:I\to\gg_{x}$.  Its integration (see Proposition
  \ref{A-G-paths} applied to the Lie algebra $\gg_{x}$) defines a path
  in $\G(\gg_x)$, and its end point is denoted by $\partial(\gamma)$.

  We need to check that $\partial$ is well defined. For that we assume
  that
  \[ 
  \gamma^{i}=\gamma^{i}(\epsilon,t):I\times I\to\O_x,\quad i\in \{0,1\}
  \]
  are two homotopic paths relative to the boundary, and that
  \[ 
  a^{i}\d t+ b^i\d\epsilon: TI\times TI\to A_{\O_x},\quad i\in \{0,1\}, 
  \]
  are lifts of $\d\gamma^{i}$ as above. We prove that the paths
  $a^{i}(1,t)$ ($i\in \{0, 1\}$) are homotopic as $\gg_x$-paths.

  By hypothesis, there is a homotopy
  $\gamma^{u}=\gamma^{u}(\epsilon,t)$ ($u\in I$) between $\gamma^{0}$
  and $\gamma^{1}$. We choose a family $b^{u}(\epsilon, t)$ joining
  $b^{0}$ and $b^{1}$, such that $\rho(b^{u}(\epsilon, t))=
  \frac{\d\gamma^u}{\d\epsilon}$ and $b^{u}(\epsilon, 0)=
  b^{u}(\epsilon, 1)= 0$. We also choose a family of sections $\eta$
  depending on $u$, $\epsilon$ and $t$, such that
  \[
  \eta^{u}(\epsilon,t,\gamma^{u}(\epsilon, t))= b^{u}(\epsilon,t),
  \text{ with }\eta=0 \text{ when } t= 0,1.
  \]
  As in the proof of Proposition \ref{equivalence}, let $\xi$ and
  $\theta$ be the solutions of
  \[
  \left\{
   \begin{array}{l}
   \frac{\d\xi}{\d\epsilon}-\frac{\d\eta}{\d t}= [\xi, \eta],\text{ with }
   \xi=0 \text{ when } \epsilon = 0, 1,\\  \\
   \frac{\d\theta}{\d\epsilon}-\frac{\d\eta}{\d u}=[\theta, \eta],\text{ with }
   \theta=0 \text{ when } \epsilon = 0, 1.
   \end{array}
  \right.\]
  Setting $u=0,1$ we get
  \[ 
  a^{i}(\epsilon,t)=\xi^{i}(\epsilon,t,\gamma^{i}(\epsilon, t)), \quad
  i=0,1.
  \]
  On the other hand, setting $t=0,1$ we get $\theta=0$ when $t=0,1$.
  A brief computation shows that
  $\phi\equiv\frac{\d\xi}{\d u}-\frac{\d\theta}{\d t}-[\xi,\theta]$
  satisfies
  \[
  \frac{\d\phi}{\d\epsilon}=[\phi, \eta],
  \]
  and since $\phi=0$ when $\epsilon=0$, it follows that
  \[ \frac{\d\xi}{\d u}- \frac{\d\theta}{\d t}= [\xi,\theta] .\]
  If in this relation we choose $\epsilon=1$, and use
  $\theta^u(1,t)=0$ when $t=0,1$, we conclude that $a^{i}(1,t)=
  \xi^i(1,t,\gamma^{i}(1, t))$, $i=0,1$, are equivalent as
  $\gg_x$-paths.

  Finally, to check that the sequence is exact we only need to check
  exactness at the level $\G(\gg_x)$. However, it is clear from the
  definition of $\partial$ that its image is exactly the subgroup of
  $\G(\gg_x)$ which consists of the equivalence classes
  $[a]\in\G(\gg_x)$ of $\gg_x$-paths with the property that, as an
  $A$-path, $a$ is equivalent to the trivial $A$-path.
\end{proof}

The previous proposition motivates our next definition:

\begin{definition}
  \label{tilde-N-groups} 
  The homomorphism $\partial: \pi_{2}(\O_x)\to \G(\gg_x)$ of
  Proposition \ref{monodromy-seq} is called the \textbf{monodromy
  homomorphism} of $A$ at $x$. Its image:
  \[ 
  \tilde{\NN}_x(A)=\set{[a]\in\G(\gg_x): a\sim 0_x\text{ as an
  $A$-path}},
  \]
  is called the \textbf{monodromy group} of $A$ at $x$.
\end{definition}

The reason for using the tilde in the notation for the monodromy group
is explained in the next exercise.

\begin{exercise}
  \label{N-versus-tilde-N}
  Show that $\tilde{\NN}_{x}(A)$ is a subgroup of $\G(\gg_x)$ contained
  in the center $Z(\G(\gg_x))$, and its intersection with the
  connected component $Z(\G(\gg_x))^{0}$ of the center is isomorphic
  to 
  \begin{equation}
  \label{N-groups}
  \NN_{x}(A)=\set{v\in Z(\gg_x): v\sim 0_x\text{ as $A$-paths}}\subset\gg_x(A).
  \end{equation}
  (Hint: Use the fact that for any $g\in \tilde{\NN}_{x}(A)\subset
  \G(\gg_x)$ which can be represented by a $\gg_{x}$-path $a$,
  parallel transport $T_{a}=Ad_g:\gg_x\to \gg_x$ along $a$ is the
  identity. Then apply the exponential map $\exp: Z(\gg_{x})\to
  Z(\G(\gg_x))^{0}$).
\end{exercise}

\begin{exercise}
  Check that if $x$ and $y$ lie in the same orbit of $A$ then there
  exists a canonical isomorphism $\tilde{\NN}_{x}(A)\simeq\tilde{\NN}_{y}(A)$.
\end{exercise}

It follows from Proposition \ref{monodromy-seq} that the restricted
isotropy group is given by:
\begin{equation}
\label{id-comp}  
\G_x(A)^0=\G(\gg_x)/\tilde{\NN}_x(A).
\end{equation}
This leads immediately to the following result:

\begin{proposition}
  For any Lie algebroid $A$, and any $x\in M$, the following are
  equivalent:
  \begin{enumerate}[(i)]
  \item $\G_x(A)^0$ is a Lie group with Lie algebra $\gg_x$.
  \item $\tilde{\NN}_{x}(A)$ is closed;
  \item $\tilde{\NN}_{x}(A)$ is discrete;
  \item $\NN_{x}(A)$ is closed;
  \item $\NN_{x}(A)$ is discrete.
  \end{enumerate}
\end{proposition}

\begin{proof}
  We just need to observe that the group $\tilde{\NN}_{x}(A)$ (and hence
  $\NN_{x}(A)$) is countable since it is the image under $\partial$ of
  $\pi_1(\O)$, which is always a countable group.
\end{proof}

At this point, we notice that we have an obvious \emph{necessary} condition
for a Lie algebroid to be integrable: if $A$ is integrable, then each
$\G_x(A)^0$ is a Lie group, and hence the monodromy groups must be
discrete. This is enough to explain the non-integrability in the
examples at the beginning of this lecture. But before we can do that,
we need to discuss briefly how the monodromy groups can be explicitly
computed in many examples.

Again, let us consider the short exact sequence of an orbit
\[
\xymatrix{0\ar[r]& \gg_{\O}\ar[r]& A|_{\O}\ar[r]^{\rho}& T\O\ar[r]& 0.}
\]
and any linear splitting $\sigma: T\O\to A_{\O}$ of $\rho$. The
\emph{curvature} of $\sigma$ is the element $\Omega_{\sigma}\in
\Omega^{2}(\O;\gg_{\O})$ defined by:
\[ \Omega_{\sigma}(X, Y):= \sigma ([X, Y])- [\sigma(X), \sigma(Y)]\ .\]
In favorable cases, the computation of the monodromy can be reduced to the
following 

\begin{lemma}
  \label{compute}
  If there is a splitting $\sigma$ with the property that its
  curvature $\Omega_{\sigma}$ is $Z(\gg_\O)$-valued, then
  \[ \tilde{\NN}_{x}(A)\simeq \NN_{x}(A)= \{ \int_{\gamma} \Omega_{\sigma}:
  [\gamma]\in \pi_{2}(\O, x)\}\subset Z(\gg_x)\]
  for all $x\in \O$.
\end{lemma}

\begin{remark}
  Note that $Z(\gg_{\O})$ is canonically a flat vector bundle over
  $\O$. The corresponding flat connection can be expressed with the
  help of the splitting $\sigma$ as
  \[ \nabla^\sigma_{X}\alpha= [\sigma(X), \alpha ],\]
  and it is easy to see that the definition does not depend on
  $\sigma$. In this way $\Omega_{\sigma}$ appears as a $2$-cohomology
  class with coefficients in the local system defined by $Z(\gg_{\O})$
  over $\O$, and then the integration is just the usual pairing
  between cohomology and homotopy. In practice one can always avoid
  working with local coefficients: if $Z(\gg_{\O})$ is not already
  trivial as a vector bundle, one can achieve this by pulling back to
  the universal cover of $\O$ (where parallel transport with respect
  to the flat connection gives the desired trivialization).

  We should specify what we mean by integrating forms with
  coefficients in a local system. Assume $\omega\in \Omega^2(M; E)$ is
  a 2-form with coefficients in some flat vector bundle
  $E$. Integrating $\omega$ over a 2-cycle $\gamma: \Ss^2\to M$ means
  (i) taking the pull-back $\gamma^*\omega\in
  \Omega^2(\Ss^2;\gamma^*E)$, and (ii) integrate $\gamma^*\omega$ over
  $\Ss^2$. Here $\gamma^*E$ should be viewed as a flat vector bundle
  over $\Ss^2$ for the pull-back connection. Notice that the
  connection enters the integration part, and this matters for the
  integration to be invariant under homotopy.
\end{remark}

\begin{proof}[Proof of Lemma \ref{compute}]
  We may assume that $\O=M$, i.e., $A$ is transitive. In agreement
  with the Remark above, we also assume for simplicity that
  $Z(\mathfrak{g})$ is trivial as a vector bundle (here and below
  $\gg:=\gg_\O$).
  
  The formula above defines a connection $\nabla^{\sigma}$ on the
  entire $\mathfrak{g}$. We use $\sigma$ to identify $A$ with
  $TM\oplus \mathfrak{g}$ so the bracket becomes
  \[ [(X, v), (Y, w)]= ([X, Y], [v, w]+ \nabla_{X}^{\sigma}(w)- 
  \nabla_{Y}^{\sigma}(v)- \Omega_{\sigma}(X, Y)) .\]
  Now choose some connection $\nabla^{M}$ on $M$ and consider the
  connection $\nabla=(\nabla^{M},\nabla^{\sigma})$ on $A=TM\oplus
  \mathfrak{g}$. Note that
  \[ T_{\nabla}((X, v), (Y, w))= (T_{\nabla^{M}}(X, Y),
  \Omega_{\sigma}(X, Y)- [v, w]) \]
  for all $X, Y\in TM$, $v, w\in \mathfrak{g}$. This shows that the two
  $A$-paths $a$ and $b$ in Proposition \ref{equivalence}, will take
  the form $a=(\frac{\d\gamma}{\d t},\phi)$,
  $b=(\frac{\d\gamma}{\d\epsilon},\psi)$, where $\phi$ and $\psi$ are paths in
  $\mathfrak{g}$ satisfying
  \[ 
  \partial_{t}\psi-\partial_{\epsilon}\phi= 
  \Omega_{\sigma}(\frac{\d\gamma}{\d t},\frac{\d\gamma}{\d\epsilon})-
  [\phi, \psi].
  \] 
  Now we only have to apply the definition of $\partial$: Given
  $[\gamma]\in \pi_{2}(M,x)$, we choose the lift $a\d t+ b\d\epsilon$ of
  $\d\gamma$ with $\psi=0$ and  
  \[ 
  \phi=-\int_{0}^{\eps}
  \Omega_{\sigma}(\frac{\d\gamma}{\d t},\frac{\d\gamma}{\d\epsilon}).
  \]
  Then $\phi$ takes values in $Z(\mathfrak{g}_{x})$ and we obtain
  $\partial [\gamma]=[\int_{\gamma}\Omega_{\sigma}]$.

\end{proof}

We can now justify why the examples of Lie algebroids given at the
beginning of this lecture are non-integrable.

\begin{example}
  \label{ex:2-forms:cont:2}
  Let $\omega\in\Omega^2(M)$ be a closed two-form, and consider the Lie 
  algebroid $A_{\omega}$ (see Examples \ref{ex:2-forms} and
  \ref{ex:2-forms:cont}). This is a transitive Lie algebroid and its
  anchor fits into the short exact sequence
  \[
  \xymatrix{0\ar[r]& \L\ar[r]& A_\omega\ar[r]^{\rho}& TM\ar[r]& 0.} 
  \]
  Using the obvious splitting of this sequence, Lemma \ref{compute}
  tells us that the monodromy group
  \[ \NN_{x}(A_{\omega})= \set{ \int_{\gamma} \omega:
    [\gamma]\in \pi_{2}(M,x)}\subset \mathbb{R} \]
  is just the group $\Gamma_\omega$ of spherical periods of $\omega$. 
  
  Assume that $A_\omega$ is integrable. Then the restricted isotropy
  group $\G_x(A_\omega)^0$ is a Lie group, and so $\Gamma_\omega$ must be
  discrete. Hence, $A_\omega$ is non-integrable if the group of spherical
  periods $\Gamma_\omega$ is non-discrete.
\end{example}

\begin{example}
  \label{ex:su(2):cont} 
  Let us came back to the example of $\su(2)$ with the modified Poisson
  structure $\{~,~\}_a$. (Example \ref{ex:su(2)}). 
  Under the identification of $\su(2)^*\simeq \Rr^3$, with coordinates
  $(x,y,z)$, the Lie bracket of the cotangent Lie algebroid of $\su(2)^*$ is:
  \[ [dx,dy]=dz,\quad [dy,dz]=dx,\quad [dz,dx]=dy,\]
  while the anchor $\rho: T^*\Rr^3\to T\Rr^3$ is given by
  \[ \rho(dx)=X,\quad \rho(dy)=Y,\quad \rho(dz)=Z,\]
  where $X$, $Y$, and $Z$ are the infinitesimal generators of
  rotations around the coordinate axis:
  \[
  X=z\frac{\partial}{\partial y}-y\frac{\partial}{\partial z}, \quad 
  Y=x\frac{\partial}{\partial z}-z\frac{\partial}{\partial x}, \quad
  Z=y\frac{\partial}{\partial x}-x\frac{\partial}{\partial y}.
  \]
  
  After we rescale the Poisson bracket by the factor $a(r)$,
  the new anchor $\rho_a:T^*\Rr^3\to T\Rr^3$ is related to the old
  one by:
  \[ \rho_a=a\rho.\]

  \begin{exercise}
    Determine the Lie bracket on 1-forms for the rescaled Poisson
    structure.
  \end{exercise}

  We now compute the monodromy of this new Lie algebroid. We restrict
  to a leaf $\Ss_{r}^{2}$, with $r>0$, and we pick the splitting of
  $\rho_a$ defined by 
  \begin{align*} 
    \sigma(X)=\frac{1}{a}(dx-\frac{x}{r^2}\bar{n}),&\qquad
    \sigma(Y)=\frac{1}{a}(dy-\frac{y}{r^2}\bar{n}),\\
    \sigma(Z)&=\frac{1}{a}(dz-\frac{z}{r^2}\bar{n}),
  \end{align*}
  where $\bar{n}= \frac{1}{r}(xdx+ydy+zdz)$. Then we obtain the
  center-valued 2-form:
  \[ \Omega_{\sigma}= \frac{a'(r)r- a(r)}{a(r)^2r^4}\omega\, \bar{n}\]
  where $\omega=x\,dy\wedge dz+y\,dz\wedge dx+z\,dx\wedge dy$.
  Since $\int_{\Ss^2_r}\omega=4\pi r^3$ it follows that the monodromy
  groups are given by:
  \[  
  \NN_{(x,y,z)}\simeq 4\pi \frac{A'(r)}{r}\Zz\bar{n}\subset\Rr\bar{n},
  \]
  where $A'(r)= 4\pi \frac{a'(r)-a(r)}{a(r)^2}$ is the derivative of
  the symplectic areas. So the monodromy might vary in a non-trivial
  fashion, even nearby regular leaves.

  \begin{exercise}
  Verify that
  \[
  \G_{(x,y,z)}= \left\{
    \begin{array}{ll}
      SU(2) \qquad&\  \text{if}\ r=0 \text{ and }a(r)\not= 0,\\
      \Rr^3 \qquad&\  \text{if}\ r=0 \text{ and }a(r)=0,\\
      \Ss^1 \qquad&\  \text{if}\ r\not=0\ \text{ and }  A'(r)\not=0,\\
      \Rr \qquad&\  \text{if}\ r\not=0 \text{ and }A'(r)=0.\\
    \end{array}
  \right.
  \]
  Also, explain why we obtain a non-integrable Lie algebroid whenever the
  symplectic area undergoes a critical point (i.e, $A'(r)=0$).
  \end{exercise}
\end{example}

\begin{exercise}
  Determine the isotropy groups $\G_x(A)$ for the Lie algebroid of
  Example \ref{ex:su(n)} and explain why $A$ is not integrable.
\end{exercise}

\begin{exercise}
  Give an example where a splitting as in Lemma \ref{compute} does not
  exist.\\
  (Hint: Verify that the groups $\NN_{x}$ and $\tilde{\NN}_{x}$ are
  distinct for the cotangent Lie algebroid of the Poisson manifold
  $M=\mathfrak{su}^*(3)$, at points lying in a orbit of dimension 4.)
\end{exercise}

\section{Notes}
Special instances of the integrability problem, for special classes of
Lie algebroids, are well-known and have a positive solution. For example, 
we have:
\begin{enumerate}[(i)]
\item For algebroids over a point (i.e., Lie algebras), the
integrability problem is solved by Lie's third theorem on the
integrability of (finite dimensional) Lie algebras by Lie groups;
\item For algebroids with zero anchor map (i.e., bundles of Lie
algebras), it is Douady-Lazard \cite{DoLa} extension of Lie's third
theorem which ensures that the Lie groups integrating each Lie algebra
fiber fit into a smooth bundle of Lie groups;
\item For algebroids with injective anchor map (i.e., involutive
distributions $\F\subset TM$), the integrability problem is solved
by Frobenius' integrability theorem.
\end{enumerate}
The integrability problem for general Lie algebroids goes back to
Pradines original works \cite{Pra1,Pra2,Pra3,Pra4}. These notes contain 
little proofs, and he made the erroneous statement that all Lie algebroids
are integrable. Almeida and Molino in \cite{Mol} gave the first example of
a non-integrable Lie algebroid in connection with developability of foliations.
Other fundamental examples were discovered over many years by different authors, 
such as \'Elie Cartan's infinite continuous groups (Singer and Sternberg, \cite{Sing}), 
the integrability of infinitesimal actions on manifolds (Palais,
\cite{Palais}, Moerdijk and Mr\v{c}un \cite{MoMrc}), abstract Atiyah sequences 
and transitive Lie algebroids (Almeida and Molino \cite{Mol}; Mackenzie \cite{Mack1}), 
of Poisson manifolds (Weinstein, \cite{Wein1}) and of algebras of vector fields (Nistor, \cite{Ni}).

The strategy to integrate Lie algebras to Lie groups that we have explained above
appears in the book of Duistermaat and Kolk \cite{DuKo} (with minor changes). Alan Weinstein 
told one of us (RLF) about the possibility of applying this strategy to integrate Lie algebroids 
during a short visit to Berkeley, in late 2000. In \cite{CaFe}, Cattaneo and Felder describe 
an approach to symplectic groupoids via an infinite dimensional symplectic reduction (the Poisson 
Sigma model), and which can be seen as a special instance of this strategy (however they fail to
recognize the monodromy groups). Finally, Severa describes in the preprint \cite{Sev} a similar 
strategy and makes some comments that are reminiscent of monodromy, without giving details.

The monodromy groups were introduced by us in \cite{CrFe1}. Of course these
groups have appeared before (in desguised form) in connection with some special 
classes of Lie algebroids. One example is the class of transitive groupoids, 
where the monodromy is equivalent to Mackenzie's cohomological class \cite{Mack1} 
that represents the obstruction to integrability. Another example, is provided by
regular Poisson manifolds, where the monodromy groups appeared in the work of 
Alcade Cuesta and Hector \cite{AlHe}.

\lecture{Integrability: Smooth Theory}%
\label{integrability:II}           %
%

In the previous lecture, for \emph{any} Lie algebroid $A$, we have
constructed a topological groupoid $\G(A)$, by taking the equivalence
classes of $A$-paths modulo $A$-homotopy. Moreover, we noticed that
an obvious \emph{necessary} condition for a Lie algebroid to be
integrable is that the monodromy groups be discrete. In this lecture,
we use the monodromy groups to give a complete answer to the problem
we raised at the beginning of that lecture:
\begin{itemize}
\item What are the precise obstructions to integrate a Lie algebroid?
\end{itemize}
In the next lecture, we will give applications of the
integrability criteria.

\section{The Main Theorem}
Here is one of the possible statements of our main integrability
criteria:

\begin{theorem} 
  For a Lie algebroid $A$, the following statements are equivalent:
  \begin{enumerate}[(a)]
  \item $A$ is integrable.
  \item $\G(A)$ is smooth.
  \item the monodromy groups $\NN_{x}(A)$ are locally uniformly discrete.
  \end{enumerate}
  Moreover, in this case, $\G(A)$ is the unique $s$-simply connected Lie
  groupoid integrating $A$. 
\end{theorem}

The precise meaning of (c) is the following: for any $x\in M$, there
is an open set $U\subset A$ containing $0_x$, such that
$\NN_{y}(A)\cap U=\{0\}$, for all $y$'s close enough to $x$.  
One can have a better understanding of the nature of this condition 
if one splits it into two conditions: a longitudinal one and a
transversal one, relative to the orbits of $A$. To do this, we fix a
norm on $A$, i.e., a norm on each fiber $A_x$ which varies 
continuously with respect to $x$(\footnote{For example, one can use a
simple partition of unit argument to construct even a smooth field of
norms.}), and we use it to measure the discretness of the monodromy
group $\NN_x(A)$: 
\[ r_{\NN}(x)= d(0, \NN_{x}(A)-\{0\}).\]
Here we adopt the convention $d(0,\emptyset)=+\infty$. Note that the
group $\NN_{x}$ is discrete iff $r_{\NN}(x)>0$. Now we
can restate our theorem as follows: 

\begin{theorem} 
  \label{thm:integrability:main}
  A Lie algebroid $A$ over $M$ is integrable if and only if,
  for all $x\in M$, the following conditions hold:
  \begin{enumerate}[(i)]
  \item $r_{\NN}(x)>0$;
  \item $\liminf_{y}r_{\NN}(y)>0$, where the limit is over $y\to
  x$ with $y$ outside the orbit through $x$.
\end{enumerate}
\end{theorem}

Notice that, since $\NN_{x}\simeq\NN_{y}$ whenever $x$ and $y$ are in
the same orbit, condition (i) can be seen as a longitudinal
obstruction, while condition (ii) can be seen as a transversal
obstruction.  We should emphasize, also, that the transversal
obstruction is \emph{not} a stable condition: it can hold in a deleted
neighborhood of a point without beeing true at the point. This comes
from the fact that, in general, the map $r_{\NN}$ is neither upper nor
lower semicontinuous. Also, the two integrability obstructions are
independent. 

Here are a few examples to ilustrate how wild $r_{\NN}$ can behave.

\begin{example}
  \label{ex:first:obstruction}
  Consider the Lie algebroid $A_{\omega}$ associated to a $2$-form
  $\omega$ on $M$ (cf.~Examples \ref{ex:2-forms}, \ref{ex:2-forms:cont} 
  and  \ref{ex:2-forms:cont:2}). We deduce that, as promised, $A_{\omega}$ is
  integrable if and only if $\Gamma_{\omega}$ is discrete. Note that,
  in this case, the second obstruction is void, since our Lie
  algebroid is transitive.
\end{example}

\begin{example}
  Let us look one more time at the example of $\su^*(2)$ with the
  modified Poisson structure $\{~,~\}_a$ (cf.~Examples
  \ref{ex:su(2)} and \ref{ex:su(2):cont}). If we choose the
  standard euclidian norm on the fibers $A=T^*(\su^*(2))$, the discussion
  in Example \ref{ex:su(2):cont} shows that:
  \[
  r_{\NN}(x,y,z)= \left\{
    \begin{array}{ll}
      4\pi \frac{A'(r)}{r}\Zz\bar{n}, \quad&\  \text{if}\ A'(r) \not=
      0,\\
      \\
      \infty, \quad&\  \text{otherwise}.\\
    \end{array}
  \right.
  \]
  In this case, $r_{\NN}$ is not upper semicontinuous at the critical
  points of the symplectic area $A(r)$. Also, the theorem above implies that
  the integrability of the underlying cotangent Lie algebroid is
  equivalent to the fact that the symplectic area has no critical
  points. In this example, only the second obstruction is violated.
\end{example}

\begin{example}
  Consider the manifold $M=\Ss^2\times \mathbb{H}$, where $\mathbb{H}$
  denotes the quaternions.  We take the Lie algebroid $A$ over $M$
  which, as a vector bundle, is trivial of rank $3$, with a fixed basis
  $\{e_1, e_2, e_3\}$. The Lie bracket is defined by:
  \[ [e_1, e_2]=e_3,\quad [e_2, e_3]=e_1,\quad [e_3, e_1]=e_2.\]
  {To} define the anchor, we let $v_1$, $v_2$ and $v_3$ be the vector
  fields on $\Ss^2$ obtained by restriction of the infinitesimal
  generators of rotations around the coordinate axis (the $X$, $Y$ and
  $Z$ in Example \ref{ex:su(2)}). We also consider the vector fields
  on $\mathbb{H}$ corresponding to multiplication by $i$, $j$ and $k$,
  denoted $w_1$, $w_2$ and $w_3$. The anchor of $A$ is then defined by
  \[ \rho(e_i)= (v_i, \frac{1}{2}w_i).\]
  A simple computation shows that this Lie algebroid has only two
  leaves. More precisely,
  \begin{enumerate}[(a)]
  \item the anchor is injective on $\Ss^2\times (\mathbb{H}-\{0\})$.
  \item $\Ss^{2}\times \{0\}$ is an orbit.
  \end{enumerate}
  In particular, the monodromy groups at all points outside
  $\Ss^{2}\times \{0\}$ are zero. On the other hand, when restricted to
  $\Ss^{2}\times \{0\}$, $A$ becomes the algebroid on $\Ss^2$ associated
  to the area form. Hence we obtain
  \[
  r_{\NN}(x)= \left\{
    \begin{array}{ll}
      r_0 \quad&\ \text{if}\ x\in \Ss^2\times\{0\},\\ 
      \\
      \infty \quad&\ \text{otherwise},\\
    \end{array}
  \right.
  \]
  where $r_{0}>0$. In particular, $r_\NN$ is not lower semicontinuous.
  However, in this case, there are no obstructions and $A$ is integrable.
\end{example}

\section{Smooth structure}

In order to describe the smooth structure on $\G(A)$, we first
describe the smooth structure on $P(A)$. We consider the larger space
$\tilde{P}(A)$ of all $C^1$-curves $a: I\to A$ with base path
$\gamma=\pi\circ a$ of class $C^2$. It has an obvious structure of
Banach manifold. 

\begin{exercise}
Show that the tangent space $T_{a}(\tilde{P}(A))$ consists of all
$C^0$-curves $U:I\to TA$ such that $U(t)\in T_{a(t)}A$. Using a $TM$-connection
$\nabla$ on $A$, check that such curves can be viewed as pairs
$(u,\phi)$ formed by a curve $u:I\to A$ over $\gamma$ and a curve
$\phi:I\to TM$ over $\gamma$ (the \emph{vertical} and \emph{horizontal}
component of $U$).
\end{exercise}

From now on, we fix a $TM$-connection $\nabla$ on $A$ and we use the
description of $T_{a}(\tilde{P}(A))$ given in this exercise. Also, we
will use the associated $A$-connections $\overline{\nabla}$ that were
discussed in Exercise \ref{exer:connections}.

\begin{lemma}
  \label{P=Banach} 
  $P(A)$ is a (Banach) submanifold of $\tilde{P}(A)$. Moreover, its
  tangent space $T_{a}P(A)$ consists of those paths $U= (u, \phi)$
  with the property that
  \[
  \rho(u)= \overline{\nabla}_{a}\phi.
  \]
\end{lemma}

\begin{proof}
  We consider the smooth map $F:\tilde{P}(A)\to \tilde{P}(TM)$ given
  by
  \[ \ F(a)=\rho(a)-\frac{d}{dt}\pi\circ a. \]
  Clearly $P(A)= F^{-1}(Q)$, where $Q$ is the submanifold of
  $\tilde{P}(TM)$ consisting of zero paths. Fix $a\in P(A)$, with base
  path $\gamma=\pi\circ a$, and let us compute the image of $U=
  (u,\phi)\in T_{a}\tilde{P}(A)$ by the differential
  \[ (dF)_{a}: T_{a}\tilde{P}(A)\to  T_{0_{\gamma}}\tilde{P}(TM) \ .\]
  The result will be some path $t\mapsto (dF)_{a}\cdot U(t)\in
  T_{0_{\gamma(t)}}TM$, hence, using the canonical splitting
  $T_{0_x}TM\cong T_{x}M\oplus T_{x}M$, it will have a horizontal and
  vertical component. We claim that for any connection $\nabla$, if
  $(u, \phi)$ are the components of $U$, then
  \[
  ((dF)_{a}\cdot U)^{\text{hor}}= \phi,\quad ((dF)_{a}\cdot
  U)^{\text{ver}}= \rho(u) -\overline{\nabla}_{a}\phi.
  \]
  Note that this immediately implies that $F$ is transverse to $Q$, so
  the assertion of the proposition follows. Since this decomposition
  is independent of the connection $\nabla$ and it is local (we can
  look at restrictions of $a$ to smaller intervals), we may assume
  that we are in local coordinates, and that $\nabla$ is the standard
  flat connection. We let $x=(x^{1},\dots,x^{n})$ denote local
  coordinates on $M$, and we denote by $\frac{\partial}{\partial
  x^{i}}$ the horizontal basis of $T_{0_x}TM$, and by
  $\frac{\delta}{\delta x^{i}}$ the vertical basis.  Also, we denote
  by $\set{e_{1},\dots,e_ {k}}$ a (local) basis of $A$ over this
  chart. The anchor and the bracket of $A$ decompose
  as(\footnote{Recall that we use the Einstein sum convention.})
  \[ 
  \rho(e_{p})= b^{i}_{p}\frac{\partial}{\partial x^{i}},\qquad
  [e_p, e_q]= c_{pq}^{r} e_{r},
  \] 
  and an $A$-path $a$ can be written $a(t)= a^{p}(t)e_p$.  A simple
  computation shows that the horizontal component of
  $(dF)_{a}(u,\phi)$ is $\phi^{i}\frac{\partial}{\partial x^{i}}$,
  while its vertical component is
  \[
  \left(
    -\dot{\phi}^{j}(t)+ u^{p}(t) b_{p}^{j}(\gamma(t))
    +a^{p}(t)\phi^{i}(t)\frac{\partial b_{p}^{j}}{\partial
      x^i}(\gamma(t)) \right) \frac{\delta}{\delta x^{j}}.
  \] 
  That this is precisely $\rho(u) - \overline{\nabla}_{a}\phi$ immediately
  follows by computing
  \begin{align*}
    \overline{\nabla}_{e_{p}}\frac{\partial}{\partial x^{i}} & = 
    \rho(\nabla_{\frac{\partial}{\partial x^{i}}}e_{p})-
    \left[\frac{\partial}{\partial x^{i}}, \rho(e_p)\right]\\ & =
    -\frac{\partial b_{p}^{j}}{\partial x^i}
    \frac{\partial}{\partial x^{j}}
  \end{align*}
\end{proof}

Finally, we can be more precise about the smoothness of $\G(A)$.

\begin{definition} 
  We say that $\G(A)$ is smooth if it admits a smooth structure with the
  property that the projection $P(A)\to \G(A)$ is a submersion.
\end{definition}

Note that, by the local form of submersions, it follows that the
smooth structure on $\G(A)$ will be unique if it exists. Invoking the
local form of submersions requires some care, since we are in the
infinite dimensional setting. However, there are no problems in our
case since, as we shall see in the next section, the fibers of the
projection map are submanifolds of $P(A)$ of finite codimension.

\begin{exercise} 
  Check that if $\G(A)$ is smooth, then $\G(A)$ is a Lie groupoid.\\
  (Hint: At the level of $P(A)$, the structure maps (i.e., concatenation,
  inverses, etc.) are smooth).
\end{exercise}

\section{A-homotopy revisited}

Before we can proceed with the proof of our main result, we need a
better control on the $A$-homotopy equivalence relation defining the
groupoid $\G(A)$. In this section we will show that we can describe
the equivalence classes of the $A$-homotopy relation $\sim$
as the leaves of a foliation $\F(A)$ on $P(A)$, of finite
codimension. Hence, $\G(A)$ is the leaf space of a foliation $\F(A)$
on $P(A)$ of finite codimension.

Using that the basic connection $\overline{\nabla}$ associated to
$\nabla$ commutes with the anchor (see Exercise
\ref{exer:connections}), there is a simple way of constructing vector
fields on $P(A)$. To make this more explicit, for a path $\gamma$ in $M$,
we consider the space
\[ 
\tilde{P}_{\gamma}(A)= 
\{ b\in \tilde{P}(A): b(0)= 0, b(t)\in A_{\gamma(t)}\}.
\]

\begin{definition} 
  Let $\nabla$ be a connection on $A$. Given an $A$-path $a_0$ with base
  path $\gamma_0$, and $b_0\in \tilde{P}_{\gamma_0}(A)$, we define
  \[ X_{b_0, a_0}\in T_{a_0}\tilde{P}(A) \]
  as the tangent vector with components $(u,\phi)$ relative to the
  connection $\nabla$ given by
  \[
  u= \overline{\nabla}_{a_0}b_0,\quad \phi=\rho(b_0).
  \]
\end{definition}

The following lemma gives a geometric interpretation of these vector
fields, showing that they are ``tangent'' to $A$-homotopies:

\begin{lemma}
  \label{X-b-a:indep} 
  The vectors $X_{b_0, a_0}$ are tangent to $P(A)$ and they do not
  depend on the choice of the connection. More precisely, choosing $a$
  and $b$ satisfying the homotopy condition (\ref{diffeq}):
  \begin{equation}
    \label{diffeq:2}
    \partial_{t}b-\partial_{\epsilon} a= T_{\nabla}(a, b),\quad
    b(\epsilon,0)=0,
  \end{equation}
  and such that $a(0,t)=a_0(t)$, $b(0,t)= b_0(t)$, then
  \begin{equation}
    \label{X-b-a}
    X_{b_0,a_0}={\left.\frac{\d}{\d\epsilon}\right|}_{\epsilon= 0}a_{\epsilon}
    \in T_{a_{0}}P(A).
  \end{equation}
\end{lemma}

\begin{proof} 
The components of $X_{b_0,a_0}$ clearly satisfy the condition from
Lemma \ref{P=Banach}, hence we obtain a vector tangent to $P(A)$.
Choosing $a$ and $b$ as in the statement, we compute the vertical and
horizontal components $(u,\phi)$ with respect to $\nabla$ of the
tangent vector
\[ 
 {\left. \frac{\d}{\d\epsilon}\right|}_{\epsilon= 0} a_{\epsilon}(t)
 \in T_{a_{0}(t)}A .
\] 
We find:
\begin{align*}
u&= {\left.\partial_{\epsilon}a\right|}_{\epsilon= 0}
    =\partial_t b_0-T_{\nabla}(a_0, b_0)=  \overline{\nabla}_{a_0}b_0,\\
\phi&={\left. \frac{\d}{\d\epsilon}\right|}_{\epsilon= 0}\gamma_{\epsilon}(t)
    =\rho(b_{0}(t)),
\end{align*}
which are precisely the components of $X_{b_0,a_0}$.
\end{proof}

\begin{remark}
\label{X-a-b:flow} 
One should point out that, in the previous lemma, given $a_0$
and $b_0$, one can allways extend them to families $a_{\epsilon}$ and
$b_{\epsilon}$ satisfying the required conditions. For instance, write
$b_0(t)=\eta(t,\gamma_{0}(t))$ for some time dependent section
$\eta=\eta(t,x)$ with $\eta_0= 0$ (one may also assume that $\eta$ has
compact support) and put 
\[ 
\gamma(\epsilon, t)= 
\phi_{\rho(\eta_{t})}^{\epsilon,0}(\gamma_0(t)),\quad 
b(\epsilon,t)=  \eta(t,\gamma(\epsilon,t)).\]
Next, one solves the equation in (\ref{diffeq:2}), with the unknown
$a= a(\epsilon, t)$, satisfying the initial condition $a(0,t)= a_{0}(t)$. 

Of course, the solution $a$ can be described without any reference to
the connection, and this is what we have seen in Proposition
\ref{equivalence}: write $a_{0}(t)=\xi_{0}(t, \gamma_0(t))$ for some
time dependent $\xi_0$, then consider the $(\epsilon, t)$-dependent
solution $\xi$ of
\begin{equation}
\label{X-a-b-eq}
\frac{\d\xi}{\d\epsilon}- \frac{\d\eta}{\d t}= [\xi,\eta], \quad
\xi(0,t,x)=\xi_0(t,x),
\end{equation}
and, finally, 
\[ a(\epsilon,t)=\xi(\epsilon,t,\gamma(\epsilon,t)). \]
\end{remark}
\vskip 15 pt

Next, we move towards our aim of realizing homotopy classes of
$A$-paths as leaves of a foliation on $P(A)$. Due to the definition of
the homotopy, we have to restrict to paths $b_0$ vanishing at the end
points. So, for a path $\gamma$ in $M$, consider
\[ \tilde{P}_{0,\gamma}(A)= \{b\in \tilde{P}_{\gamma}(A): b(1)= 0 \}.\]

\begin{definition} 
  We denote by $\D$ the distribution on $P(A)$ whose
  fiber at $a\in P(A)$ is given by
  \[ \D_a\equiv \set{ X_{b, a}\in T_{a}P(A): b\in
    \tilde{P}_{0,\gamma}(A)}.
  \]
\end{definition}

Now we have:

\begin{lemma} 
A family $\{a_{\epsilon}\}$ of $A$-paths is a homotopy if and only if the
path $I\ni \epsilon \mapsto a_{\epsilon}\in P(A)$ is tangent to $\D$. 
In particular, two $A$-paths are homotopic if and only if they are 
connected by a path tangent to $\D$.
\end{lemma}

\begin{proof}
We will continue to use the notations from Proposition \ref{equivalence}
(in particular, we assume that a connection $\nabla$ on $A$ has been
fixed).
Exactly as in the proof of Lemma \ref{X-b-a:indep}, we compute the
components of the derivatives (this time at arbitrary $\epsilon$'s)
\[ 
  \frac{\d}{\d\epsilon} a_{\epsilon}\in T_{a_{\epsilon}}P(A)
\] 
and we find that this derivative coincides with 
$X_{b_{\epsilon},a_{\epsilon}}$.  Hence, if $\{a_{\epsilon}\}$ is a
homotopy, then $b_{\epsilon}(1)=0$ and we deduce that $\epsilon\mapsto
a_{\epsilon}$ is tangent to $\D$. 

For the converse, we fix $\epsilon$ and we want to prove that
$b_{\epsilon}(1)=0$. But the assumption implies that there exists
$c\in \tilde{P}_{0,\gamma_{\epsilon}}(A)$ (remember that we have
fixed $\epsilon$) such that $X_{c,a_{\epsilon}}= X_{b_{\epsilon},
a_{\epsilon}}$. Writing the components with respect to a connection,
we find that
\[ \overline{\nabla}_{a}(b_{\epsilon}-c)= 0.\]
Since $b_{\epsilon}(0)-c(0)=0$, we see that $b_{\epsilon}= c$. In
particular, $b_{\epsilon}(1)=c(1)=0$.
\end{proof}

Next we show that the distribution $\D$ is integrable. Set:
\[ 
  P_{0}\Gamma(A):= \set{ I\ni t\mapsto \eta_{t}\in \Gamma(A): 
  \eta_{0}= \eta_{1}=0, \eta\text{ is of class } C^2\text{ in } t}.
\]
Then $P_{0}\Gamma(A)$ is a Lie algebra with Lie bracket the pointwise
bracket induced from $\Gamma(A)$. For each $\eta\in P_{0}\Gamma(A)$
define a vector field $X_{\eta}$ on $P(A)$ as follows. If $a\in P(A)$,
let $b=\eta(t,\gamma(t))$, where $\gamma$ is the base path of $a$, and
set:
\[ {\left. X_\eta\right|}_ {a}:= X_{b,a}\in T_{a}P(A).\]

\begin{exercise} 
Show that, $X_\eta$ can be written in terms of flows as:
\[
{\left. X_\eta\right|}_ {a}
={\left. \frac{\d}{\d\epsilon}\right|}_{\epsilon= 0}
\phi_{\eta_{t}}^{\epsilon, 0}(a(t))+\frac{\d\eta_{t}}{\d t}(\gamma(t)).
\]
(Hint: Show that once $\xi_{0}$ has been chosen, then $\xi$ is given
by:
\[ 
\xi(\epsilon,t)=
\int_{0}^{\epsilon} (\phi_{\eta_{t}}^{\epsilon', \epsilon})^{*}
                     \frac{\d\eta}{\d t}(t)\d\epsilon'
                   + (\phi_{\eta_{t}}^{0, \epsilon})^{*}\xi_{0}\]
(cf.~Remark \ref{X-a-b:flow}). Deduce the explicit integral formula:
\[
a_{\epsilon}(t)= \int_{0}^{\epsilon}
                    \phi_{\eta_{t}}^{\epsilon, \epsilon'}\left(
                    \frac{\d\eta_{t}}{\d t}
                    (\phi_{\rho(\eta_{t})}^{\epsilon',  0}(\gamma_{0}(t))
                    \right)\d\epsilon'
                  + \phi_{\eta_{t}}^{\epsilon, 0}(a_{0}(t))
\]
(cf.~Proposition \ref{equivalence}). Now differentiate at $\epsilon=0$
to find ${X_\eta}|_a$.)
\end{exercise}
 
Now we can show that $\D$ is an involutive distribution so that there
exists a foliation $\F(A)$ integrating it: 

\begin{lemma} 
  The vector fields $X_{\eta}$ span the distribution $\D$. Moreover,
  the map
  \[  P_{0}\Gamma(A) \to \X(P(A)),\quad \eta\mapsto X_{\eta} \]
  defines an action of the Lie algebra $P_{0}\Gamma(A)$ on $P(A)$.
  In particular, $\D$ is involutive.
\end{lemma}

\begin{proof} 
  The first assertion is clear since for any $b\in\tilde{P}_{0,\gamma}(A)$
  we can find $\eta\in P_{0}\Gamma(A)$ such that $b(t)= \eta(t,\gamma(t))$.
  The proof of the second assertion follows by a computating in local
  coordinates. 
\end{proof}

Finally, we can put together the main properties of the foliation
$\F(A)$ integrating $\D$:

\begin{theorem}
  \label{G-as-fol}
  For a Lie algebroid $A$, there exists a foliation $\F(A)$ on $P(A)$
  such that:
  \begin{enumerate}
  \item[(i)] $\F(A)$ is a foliation of finite codimension equal to
    $n+k$, where $n=\dim M$ and $k=\rank A$.
  \item[(ii)] Two $A$-paths are equivalent (homotopic) if and only if
    they are in the same leaf of $\F(A)$.
\end{enumerate}
\end{theorem}

\begin{proof} 
  We already know that $\F(A)$ is a foliation whose leaves are
  precisely the homotopy classes of $A$-paths.

  To determine the codimension of $\D_a=T_a\F(A)$, we use a connection
  $\nabla$ to construct a surjective map
  \[ \nu_{a}: T_{a}P(A)\rmap A_{\gamma(1)}\oplus T_{\gamma(0)}M\]
  whose kernel is precisely $\D_a$.  Explicitely, given a vector
  tangent to $P(A)$ with components $(u, \phi)$, we put
  \[ \nu_{a}(u, \phi)= (b(1), \phi(0)),\]
  where $b$ is the solution of the equation $\overline{\nabla}_{a}(b)=
  u$ with initial condition $b(0)= 0$ (which can be expressed in terms
  of the parallel transport along $a$ with respect to
  $\overline{\nabla}$).  Assume now that $\nu_{a}(u, \phi)= 0$,
  i.e., $b(1)= 0$ and $\phi(0)= 0$.  Since $(u, \phi)$ is tangent to
  $P(A)$, we must have $\rho(u)=
  \overline{\nabla}_{a}(\phi)$. Replacing $u$ by
  $\overline{\nabla}_{a}(b)$ we deduce that
  \[ \overline{\nabla}_{a}(\rho(b)- \phi)= 0.\]
  Since $\rho(b)- \phi$ vanishes at the initial point, we deduce that
  $\phi= \rho(b)$. Hence $(u,\phi)$ gives precisely the vector $X_{b,
  a}$ tangent to $\F(A)$. We deduce that the kernel of $\nu_{a}$ is
  $\D_a$. The surjectivity of $\nu_{a}$ is easily checked: given
  $b_1\in A_{\gamma(1)}$ and $\phi_0\in T_{\gamma(0)}M$, we denote by
  $\phi$ the solution of the equation $\overline{\nabla}_{a}(\phi)= 0$
  with initial condition $\phi_0$ and we choose any path $b: I\rmap A$
  above $\gamma$ joining $0$ and $b_1$. Then
  $(\overline{\nabla}_{a}(b), \phi+ \rho(b))$ is tangent to $P(A)$ and
  is mapped by $\nu_{a}$ into $(b_1, \phi_0)$.  This proves (i).

\end{proof}

\section{The exponential map}

Recall that for a Lie group $G$ with Lie algebra $\mathfrak{g}$, the
exponential map $\exp:\mathfrak{g}\to G$ is a local diffeomorphism around
the origin. Hence the smooth structure on $G$ around the identity
element can be obtained from $\mathfrak{g}$ and the exponential map;
using then right translations, the same is true around each point in
$G$.  With this in mind, our plan is to construct a similar
exponential map for algebroids, and then use it to obtain a
smooth structure on $\G(A)$.

Let $\G$ be a Lie groupoid with Lie algebroid $A$. The exponential
map that we will use will require the choice of a connection $\nabla$
on the vector bundle $A$. Let $\nabla$ be such a connection, and take
the pull-back of $\nabla$ along the target map
\[ \t: \G(x,-)\to M.\]
This will define a connection on the vector bundle $\t^{*}A$ restricted to
$\G(x,-)$, which is isomorphic to $T\G(x,-)$. Hence we obtain a
connection, denoted by $\nabla^{x}$, on the $s$-fiber $\G(x,-)$. We
consider the associated exponential, defined on a
neighborhood of the origin in $T_{x}\G(x,-)=A_x$, which by abuse of
notation, we write it as if it was defined on the entire $A_x$:
\[ \Exp_{\nabla^x}: A_x\rmap \G(x, -).\]
Putting together these exponential maps we obtain a global exponential
map 
\[ \overline{\Exp}_{\nabla}: A\to \G\]
Again, in spite of the notation, this is only defined on an open
neighborhood of the zero section. This exponential map has the desired
property:

\begin{lemma} 
  $\overline{\Exp}_{\nabla}$ is a diffeomorphism from an open
  neighborhood of the zero section in $A$ to an open neighborhood of $M$
  in $\G$.
\end{lemma}

\begin{proof} 
  Abreviate $E=\overline{\Exp}_{\nabla}$ and $E^x=
  \overline{\Exp}_{\nabla^x}$.  It is enough to prove that, for an
  arbitrary $x\in M$, the differential of $E$ at $0_x\in A$ is a linear
  isomorphism from $T_{0_x}A$ to $T_{1_x}\G$.  To see this, we first
  remark that $(dE_x)$ is a linear isomorphism from $A_x$ into
  $T_{1_x}\G(x,-)$. Secondly, we remark that we have a commutative
  diagram whose horizontal lines are short exact sequences:
  \[
  \xymatrix{ 
    A_{x}\ar[r]\ar[d]_{(dE^x)} & T_{0_x}A\ar[r]^-{(d\pi)}\ar[d]_-{(dE)} 
    & T_xM\ar[d]^-{\text{Id}}\\
    T_x\G(x,-)\ar[r] & T_x\G\ar[r]^-{(d\s)} & T_xM 
  }
  \]
  The desired conclusion now follows from a simple diagram chasing.
\end{proof}

Next, we define a version of $\overline{\Exp}_{\nabla}$ which makes
sense on the (possibly non-smooth) $\G(A)$, and which coincides with
the construction above when $\G=\G(A)$ is smooth.

Let $A$ be any Lie algebroid and let $\nabla$ be a connection on
$A$. We denote by the same letter the induced $A$-connection on $A$:
$\nabla_{\alpha}\beta=\nabla_{\rho(\alpha)}\beta$.  Similar to the
classical case, we consider $A$-geodesics with respect to $\nabla$,
i.e., $A$-paths $a$ with the property that $\nabla_{a}a=0$. Still as
in the classical case (and exactly by the same arguments), for each
$v\in A_x$ there is a unique geodesic $a_v$ with $a_v(0)=v$. Also, for
$v$ close enough to zero, $a_v$ is defined on the entire interval $I$.
Hence we obtain a map $v\mapsto a_v$, denoted
\[ \Exp_{\nabla}: A\to P(A)\]
and defined on an open neighborhood of $M$ in $A$. Taking homotopy
classes of $A$-paths, we obtain a map $v\mapsto[a_v]$ denoted
\[ \overline{\Exp}_{\nabla}: A\rmap \G(A).\]

Next, let us check that the two versions of $\overline{\Exp}_{\nabla}$
are compatible.

\begin{lemma} 
  If $\G(A)$ is smooth, then the two constructions above for
  $\overline{\Exp}_{\nabla}: A\to \G(A)$ coincide.
\end{lemma}

\begin{proof} 
  To emphasize that the first construction of
  $\overline{\Exp}_{\nabla}$ uses the smooth structure of $\G(A)$, let
  us assume that $\G$ is an arbitrary $s$-connected Lie groupoid
  integrating $A$; at the end we will set $\G=\G(A)$. Again, we
  abreviate the notations for the exponential maps associated to
  $\nabla$: $E: A\to \G(A)$ and $E': A\to \G$. The key remark now is:
  if $g: I\to\G(x, -)$ is a geodesic with respect to $\nabla^x$, then
  the $A$-path $a=D^R(g)$ is an $A$-geodesic with respect to
  $\nabla$. This shows that the two exponential maps coincide modulo
  the isomorphism $D^R$ between $\tilde{\G}= P(\G)/\sim$ and $\G(A)$.
\end{proof}

As a consequence we obtain

\begin{corollary} 
If a Lie algebroid $A$ is integrable, then
\[ \overline{\Exp}_{\nabla}: A\to \G(A)\]
is injective around the zero section.
\end{corollary}

We will see later in the lecture, in Theorem \ref{criteria-exp}, that
the converse is also true. At this point, we can detect the monodromy
groups $\NN_x(A)$ as the obvious obstructions to the integrability of
$A$. Indeed, since $\overline{\Exp}_{\nabla}$ restricted to $\gg_x$ is
the composition of the exponential map of $\gg_x$ with the obvious map
$i: \G(\gg_x)\to \G(A)_x$, we will have 
$\overline{\Exp}_{\nabla}(v_x)= 1_x$ for all $v_x\in \NN_{x}(A)$ in
the domain of the exponential map.  We conclude:

\begin{corollary} 
  If $\G(A)$ is smooth (i.e., if $A$ is integrable), then the mo\-no\-dromy
  groups $\NN_{x}(A)$ are locally uniformly discrete.
\end{corollary}

Next, to understand how the exponential map can be used to obtain a
smooth structure on $\G(A)$ (when it is injective), we have to
understand its behaviour at the level of tha paths space.

\begin{proposition} 
  For any (local) connection $\nabla$ on $A$, the exponential map
  $\Exp_{\nabla}: A\to P(A)$ is transverse to $\F(A)$.
\end{proposition}

\begin{proof}
  We assume, for simplicity, that we are in local coordinates and that
  $\nabla$ is the trivial flat connection (this is actually all we
  will use for the proof of the main theorem, and this in turn will
  imply the full statement of (iii)).  Also, we only need to show
  that $\Exp_{\nabla}(A)$ is transverse to $\F(A)$ at any trivial
  $A$-path $a=O_{x}$ over $x\in M$. Now, the equations for the
  geodesics show that if $(u,\phi)$ is a tangent vector to
  $\Exp_{\nabla}(A)$ at $a$ then we must have:
  \[ \dot{\phi}^i=b^i_p(x)u^p,\quad \dot{u}^p=0.\]
  Therefore, we see that:
  \[ T_{a}\Exp_{\nabla}(A)=\set{(u,\phi)\in T_{a}P(A): u(t)=u_0, 
    \phi(t)=\phi_0+t\rho(u_0)}.\] 
  Suppose now that a tangent vector $(u,\phi)$ belongs to this $n+k$
  dimensional space and is also tangent to $\F(A)$. Since $(u,\phi)$
  is tangent to $\F(A)$ we must have $\phi_0=0$ and $b(1)=0$ (see the
  proof of Theorem \ref{G-as-fol}). Therefore, we must have $\phi_0=0$
  and $u_0=0$, so $(u,\phi)$ is the null tangent vector. This shows
  that $\Exp_{\nabla}(A)$ is transverse to $\F(A)$ at $0_x$, for any
  $x$.
\end{proof}

Finally, we have:

\begin{theorem}
  \label{criteria-exp} 
  $\G(A)$ is smooth if and only if the exponential map is injective
  around the zero section. Moreover, in this case $\G(A)$ is a Lie
  groupoid integrating $A$.
\end{theorem}

\begin{proof} 
  Assuming injectivity, we will prove that $\G(A)$ is smooth. We have to
  prove that the foliation $\F(A)$ on $P(A)$ is simple.  For this it
  suffices to show that, for each $a\in P(A)$, there exists
  $S_a\subset P(A)$ which is transverse to $\F(A)$, and which
  intersects each leaf of $\F(A)$ in at most one point. We will call
  such an $S_a$, a simple transversal through $a$.  

  Let $a$ be an arbitrary $A$-path. We denote by $\gamma$ its base
  path and by $x\in M$ its initial point, and we let
  $a(t)=\xi(t,\gamma(t))$ for some compactly supported, time
  dependent, section $\xi$ of $A$. Then we define
  \[ \sigma_{\xi}: M\to P(A),\ 
  \sigma_{\xi}(y)(t)= \xi(t, \phi_{\rho\xi}^{t, 0}(y) .\]
  Composing on the left with $\sigma_{\xi}$, we obtain a smooth
  injective map
  \[ T_{\xi}: P(A)\to P(A),\ T_{\xi}(b)= \sigma_{\xi}(\t(b))b .\]
  By the hypothesis and the previous proposition, there exists some open
  set $U\subset A$ containing $0_x$ such that 
  \[ \Exp_{U}=\Exp_{\nabla}: U\to P(A)\] 
  induces a simple transversal $S_x\subset P(A)$ through
  $0_x$. Applying $T_{\xi}$, we obtain a simple transversal
  $T_{\xi}(S_x)$ through $T_{\xi}(0_x)$. But $T_{\xi}(0_x)=0_{\xi}a$
  is homotopic to (hence in the same leaf) as $a$. Using the holonomy
  of the foliation $\F(A)$ along any path from $T(0_x)$ to $a$, we
  obtain a simple transversal through $a$.
  
  Next, we have to prove that the Lie algebroid of $\G(A)$ is
  isomorphic to $A$.  It is clear from the definitions that $A$ can be
  identified with $T^{\s}_M \G(A)$, as vector bundles, and that under
  this identification $\rho$ coincides with the differential of the
  target $\t$. So we need only to check that the bracket of
  right-invariant vector fields on $\G(A)$ is identified with the
  bracket of sections of $A$. For this we note that, on one hand, the
  bracket is completely determined by the infinitesimal flow of
  sections through the basic formula (\ref{Lie-flows}). On the other
  hand, we also know that the exponential
  $\exp:\Gamma(A)\to\Gamma(\G(A))$ is injective in a neighborhood of
  the zero section, and so Exercise \ref{exer:exponential} shows that
  the infinitesimal flow of a section $\al$ is the infinitesimal flow
  of the right-invariant vector field on $\G(A)$ determined by
  $\al$. Hence, we must have $A(\G(A))=A$.
\end{proof}

\section{End of the proof: injectivity of the exponential}

In this section we complete the proof of the main theorem. 
According to Theorem \ref{criteria-exp} in the last section, 
it suffices to show: 

\begin{proposition}
  Any $x\in M$ admits a neighborhood $V$ such that the exponential
  $\overline{\Exp}_{\nabla}: V\to \G(A)$ is injective.
\end{proposition}

In other words, we look for $V$ such that $\Exp_{\nabla}: V\to P(A)$
intersects each leaf of $\F(A)$ in at most one point. This will be
proven in several steps by a sequence of reductions and careful
choices. 

Fix $x\in M$, choose local coordinates around $x$, and let $\nabla$ be
the canonical flat connection on the coordinate neighborhood. We also
choose a small neighborhood $U$ of $0_x$ in $A$ so that the
exponential map $\Exp_{\nabla}: U\to P(A)$ is defined and is
transverse to $\F(A)$. 

\begin{clm} 
  We may choose $U$ such that, for any $v\in U\cap \mathfrak{g}_y$ ($y\in M$)
  with the property that $\Exp_{\nabla}(v)$ is homotopic to $0_{y}$, we must 
  have $v\in Z(\mathfrak{g}_y)$.
\end{clm}

Given a norm $|\cdot |$ on $A$, the set $\set{|[v, w]|: v, w\in
\mathfrak{g}_{y}\text{ with } |v|=|w|=1}$, where $y\in M$ varies in a
neighborhood of $x$, is bounded. Rescaling $|\cdot |$ if necessary, we
find a neighborhood $D$ of $x$ in $M$, and a norm $|\cdot |$ on
$A_{D}= \{ v:\pi(v)\in D\}$, such that $|[v, w]|\leq |v| |w|$ for all
$v,w\in\mathfrak{g}_y$ with $y\in D$. We now choose $U$ so that
$U\subset A_D$ and $|v|< 2\pi $ for all $v\in U$.  If $v$ is as in the
claim, it follows (see Exercise \ref{exer:hol:hom}) that parallel
transport $T_{v}: \mathfrak{g}_{y}\to \mathfrak{g}_y$ along the
constant $A$-path $v$ is the identity. But $T_v$ is precisely the
exponential of the linear map $\ad_{v}: \mathfrak{g}_{y}\to
\mathfrak{g}_y$.  Since $||\ad_v||\leq |v|< 2\pi$, it follows that all
the eigenvalues of $\ad_v$ are of norm less then $2\pi$, and then we
deduce that all these eigenvalues must be zero. In conclusion
$\ad_v=0$, and the claim follows.

\begin{clm}
  We may choose $U$ such that, if $v\in U\cap \mathfrak{g}_y$ ($y\in M$)
  has the property that $\Exp_{\nabla}(v)$ is homotopic to $0_y$ then
  $v= 0_y$.
\end{clm}

Obviously this is just a restatement of the obstruction assumptions,
combined with the previous claim.

\begin{clm}
  We may choose $U$ such that, if $v\in U$ has the property that the
  base path of $\Exp_{\nabla}(v)$ is closed, then $v\in \mathfrak{g}_y$.
\end{clm}

To see this, we write the equations for the geodesics:
\[
\dot{x^i}=b^i_p(x(t))a^p,
\quad
\dot{a^p}=0,
\]
and we apply the following Period Bounding Lemma:

\begin{lemma} 
  Let $\Omega\subset \mathbb{R}^n$ be an open set, let $F: \Omega\to
  \mathbb{R}^n$ be a smooth map, and assume that
  \[ \sup_{x\in \Omega}
  \left|\frac{\partial F^j}{\partial x_k}(x)\right|\leq L, \quad (1\le
  j,k\le n).\]
  Then, any periodic solution of the equation
  $\dot{x}(t)=F(x(t))$ must have period $T$ such that
  \[ T\geq \frac{2\pi}{L}. \]
\end{lemma}

Hence, it suffices to choose $U\subset A_{D}$,
where $D$ is chosen small enough so that 
\[ \sup_{\substack{1\le j,k\le m\\x\in D}}\left\Vert
  \frac{\partial b^j_p}{\partial x^k}(x)a^p\right\Vert< 2\pi,\]
and the claim follows.
\vskip 10 pt

Now, for any open set $O\subset P(A)$, we consider the plaques in $O$
of $\F(A)$, or, equivalently, the leaves of $\F(A)|_{O}$. For $a, b\in
O$, we write $a\sim_{O} b$ if $a$ and $b$ lie in the same plaque.
From now on, we fix $U$ satisfying all the conditions above, and we
choose an open set $O$ such that $\Exp_{\nabla}: U\to P(A)$ intersects
each plaque inside $O$ exactly in one point. This is possible since
$\Exp_{\nabla}$ is transversal to $\F(A)$. Apart from the pair $(O,
U)$, we also choose similar pairs $(O_{i}, U_{i})$, $i= 1, 2$, such
that $O_1O_1\subset O$, $O_{2}O_{2}\subset O_{1}$ and $O_{i}^{-1}=
O_{i}$.

\begin{clm}
  It is possible to choose a neighborhood $V$ of $x$ in $U_{2}$ so that, 
  for all $v\in V$,
  \[ 0_{y}\cdot \Exp_{\nabla}(v) \sim_{O} \Exp_{\nabla}(v) .\] 
\end{clm}
We know that for any $v$ there is a natural homotopy between the two
elements above. This homotopy can be viewed as a smooth map $h:
I\times U \to P(A)$ with $h(0,v)= 0\cdot \Exp_{\nabla}(v)$, $h(1,
v)= \Exp_{\nabla}(v)$, $h(t,0_{y})= 0_{y}$. Since $I$ is compact and
$O$ is open, we can find $V$ around $x$ such that $h(I\times V)\subset
O$. Obviously, $V$ has the desired property. 

\begin{clm}
  It is possible to choose $V$ so that, for all $v, w\in V$,
  \[
  (\Exp_{\nabla}(v)\cdot\overline{\Exp_{\nabla}(w)}) \cdot \Exp_{\nabla}(w)
  \sim_{O} \Exp_{\nabla}(v)
  \]
\end{clm}
This is proved exactly as the previous claim. Ou final claim is:

\begin{clm}
  $\Exp_{\nabla}:V\to P(A)$ intersects each leaf of $\F(A)$ in at most
  one point.
\end{clm}
To see this, let us assume that $v, w\in V$ have $\Exp_{\nabla}(v)\sim
\Exp_{\nabla}(w)$.  Then $a_1:= \Exp_{\nabla}(v)\cdot
\overline{\Exp_{\nabla}(w)}\in O_{1}$ will be homotopic to the trivial
$A$-path $0_y$. On the other hand, by the choice of the pair $(O_1,
U_1)$, $a_1\sim_{O_1} \Exp_{\nabla}(u)$ for an unique $u\in
U_1$. Since $\Exp_{\nabla}(u)$ is equivalent to $0_y$, its base path
must be closed, hence, by claim $4$ above, $u\in \mathfrak{g}_y$.
Using Claim $3$, it follows that $u= 0$, hence $a_{1}\sim_{O_1}
0_{y}$.  Since $O_1O_1\subset O$, this obviously implies that
\[ a_{1}\cdot \Exp_{\nabla}(w) \sim_{O} 0_{y}\cdot\Exp_{\nabla}(w) .\] 
Since $V$ satisfies Claim $5$ and Claim $6$, we get
$\Exp_{\nabla}(v)\sim_{O}\Exp_{\nabla}(w)$. Hence, by the construction
of $O$, $v= w$. This concludes the proof of the main theorem.
\qed

\begin{exercise}
How can one modify this proof, to show that, in the main theorem, it
suffices to require that for each leaf $L$, there exists $x\in L$
satisfying the two obstructions.
\end{exercise}

\begin{exercise}
Define the notion of a \emph{local Lie groupoid} (so that the
structure maps are defined only on appropriate small neighborhoods of
the identity section) and adapt this proof to show that every Lie
algebroid integrates to a local Lie groupoid. 
\end{exercise}

\section{Notes}

Theorem \ref{thm:integrability:main} is the main theorem of these 
lectures and appears first in \cite{CrFe1}. The proof is the
same that we describe here, and makes use of the key Period Bounding
Lemma which is due to Yorke \cite{Yorke}. 

\lecture{An example: integrability and Poisson geometry}%
\label{examples}                   %

In this last lecture we will illustrate some of the previous results
on integrability of Lie bracket with an application to Poisson geometry.

\section{Integrability of Poisson brackets} %
\label{sec:Poisson:geometry}                %

As we saw in Example \ref{ex:Poisson}, to a Poisson
manifold $(M,\{~,~\})$ there is associated a cotangent Lie algebroid
$(T^*M,\rho,[~,~])$ where the anchor $\rho:\Omega^1(M)\to\X(M)$ and
the bracket $[~,~]:\Omega^1(M)\times\Omega^1(M)\to\Omega^1(M)$, on
exact 1-forms is given by:
\begin{equation}
  \label{kozul:bracket:exact}
  \rho(\d f)=\pi^\sharp\d f=X_f,\quad [\d f_1,\d f_2]=\d\{f_1,f_2\},
\end{equation}
for all $f,g\in C^\infty(M)$. Moreover, these formulas determine
uniquely the Lie algebroid structure. Exercise \ref{ex:Koszul} shows that the
Koszul bracket defined by
\begin{equation}
  \label{kozul:bracket}
  [\alpha,\beta]= \Lie_{\pi^{\sharp}(\alpha)}(\beta)- 
  \Lie_{\pi^{\sharp}(\beta)}(\alpha)-\d(\pi(\alpha, \beta)),
\end{equation}
is the unique Lie algebroid bracket on $\Omega^1(M)$ that satisfies 
(\ref{kozul:bracket:exact}).

\begin{exercise} 
  Show that a Lie algebroid structure on $T^*M$ is induced from a
  Poisson bracket on $M$ if and only if (i) the anchor $\rho:T^*M\to
  TM$ is skew-symmetric and (ii) the bracket of closed 1-forms is a
  closed 1-form.
 \end{exercise}

Therefore, Poisson manifolds form a nice class of Lie algebroid
structures. Moreover, all the symplectic and foliated geometry that
combines into a Poisson manifold, often leads to beautiful geometric
interpretations of many of the results and constructions in abstract
Lie algebroid theory. Then, a basic questions is: 
\begin{itemize}
\item What are the Lie groupoids integrating Poisson manifolds?
\end{itemize}

Before we initiate our study of the integrability of Poisson brackets,
it is worth pointing out that there exists yet another connection between
Poisson geometry and Lie groupoid/algebroid theory. In fact, one can
think of Lie algebroids as a special class of Poisson
manifolds, a generalization of the well-known equivalence between finite
dimensional Lie algebras and linear Poisson brackets. 

To explain this, let $p:A\to M$ be a Lie algebroid over a
manifold $M$, with anchor $\rho:A\to TM$, and Lie bracket
$[~,~]:\Gamma(A)\times\Gamma(A)\to\Gamma(A)$. We define a Poisson
bracket on the total space of the dual bundle $A^*$ as follows. The
algebra $C^\infty(A^*)$ of smooth functions on $A^*$ is generated 
by the following two types of functions:
\begin{itemize}
\item the \emph{basic functions} $f\circ p$, where $f\in
 C^\infty(M)$, and
\item the \emph{evaluation by a section} $\al\in\Gamma(A)$, denoted
  $F_\al$, and defined by
  \[ F_\al(\xi):=\langle\al,\xi\rangle\quad (\xi\in A^*).\]
\end{itemize}
In order to define the Poisson bracket on $A^*$ it is enough to define
it on these two kinds of functions, in a way that is compatible with
the Jacobi identity and the Leibniz identity. We do this as follows:
\begin{enumerate}[(i)]
\item The bracket of basic functions is zero:
  \[ \{f\circ \pi,g\circ \pi\}_{A^*}=0,\quad (f,g\in C^\infty(M));\]
\item The bracket of two evaluation functions is the evaluation function of
  their Lie brackets:
  \[ \{F_\al,F_\be\}_{A^*}=F_{[\al,\be]}\quad (\al,\be\in\Gamma(A));\]
\item The bracket of a basic function and an evaluation function is the
  basic function given by applying the anchor:
  \[\{F_\al,f\circ \pi\}_{A^*}=\rho(\al)(f)\circ \pi\quad
  (\al\in\Gamma(A), f\in C^\infty(M)).\]
\end{enumerate}

\begin{exercise}
  Show that these definitions are compatible with the Jacobi and
  Leibniz identities, so they define a Poisson bracket on $A^*$.
\end{exercise}

The following exercise gives some basic properties of this Poisson
bracket.

\begin{exercise}
  Show that the Poisson bracket just defined satisfies the following
  properties: 
  \begin{enumerate}[(a)] 
  \item it is fiberwise linear, i.e., the bracket of functions linear on the
    fibers is a function linear on the fibers;
  \item the hamiltonian vector fields $X_{F_\al}$ project to the vector
    fields $\rho(\al)$.
  \item the flow of a hamiltonian vector field $X_{F_\al}$ is the
    fiberwise transpose of the flow of the section $\al\in\Gamma(A)$:
    \[
    \langle \phi^t_\al(a),\xi\rangle=
    \langle a,\phi^t_{X_{F_\al}}(\xi)\rangle,\quad \forall\ a\in A,\
    \xi\in A^*.
    \]
  \end{enumerate}
  Conversely, show that a fiberwise linear Poisson structure on a
  vector bundle $A^*\to M$ induces a Lie algebroid structure on the dual
  bundle $A$, whose associated Poisson bracket is the original one.
\end{exercise}

Therefore, we see that Lie algebroids and Poisson geometry come
together hand on hand. This lecture is an ellaboration on this theme. 

\section{Contravariant geometry and topology}

For a Poisson manifold $(M,\{~,~\})$ there is a contravariant version
of geometry, dual in certain sense to the usual covariant
geometry.

Let us start by recalling a few basic notions of Poisson
geometry. First of all, we recall that a \emph{Poisson map} is a map
$\phi:M\to N$ between Poisson manifolds that preserves the Poisson
brackets:
\[ \{f\circ\phi,g\circ\phi\}_M=\{f,g\}_N\circ\phi, \]
for all $f,g\in C^\infty(N)$. We recall also that a \emph{Poisson
vector field} is an infinitesimal automorphisms of $\pi$, i.e., a
vector field $X\in\X(M)$ such that $\Lie_X\pi=0$. We denote by
$\X_{\pi}(M)$ the space of Poisson vector fields and by
$\X_{\Ham}(M)\subset \X_{\pi}(M)$ the subspace of Hamiltonian vector
fields.

Some of the concepts we have introduced before for a general Lie
algebroid $A$, are known in Poisson geometry (when $A=T^*M$) under
different names. For example, $T^*M$-paths, $T^*M$-loops and
$T^*M$-homotopies are known as \textbf{cotangent paths},
\textbf{cotangent loops} and \textbf{cotangent homotopies},
respectively. We will denote by $P_\pi(M)$ the space of cotangent
paths (warning: the notation $P(T^*M)$ will have a different meaning
in this lecture). On the other hand, a $T^*M$-connection on a vector
bundle $E$ over a Poisson manifold is known as a \textbf{contravariant
connection} on $E$. Of special interest are the contravariant
connections on $E=T^*M$. By the usual procedure such a connection
induces contravariant connections on any associated bundle ($TM$,
$\wedge TM$, $\wedge T^*M$, $\End(TM)$, etc.)

\begin{exercise}
  Let $(M,\pi)$ be a Poisson manifold. Show that there exists a
  contravariant connection $\nabla$ on $T^*M$ such that $\nabla
  \pi=0$. On the other hand, show that there exists an ordinary
  connection $\bar{\nabla}$ on $T^*M$ such that $\bar{\nabla}\pi=0$
  iff $\pi$ is a regular Poisson structure.
\end{exercise}

The cohomology of the cotangent Lie algebroid of a Poisson manifold
$(M,\pi)$ is known as the \textbf{Poisson cohomology} of the Poisson
manifold $(M,\pi)$ and will be denoted $H_\pi^\bullet(M)$ (instead of
$H^\bullet(T^*M)$).

\begin{exercise}
Show that Poisson cohomology is just the cohomology of the complex of
multivector fields:
\[ \X^\bullet(M):=\Gamma(\wedge^\bullet TM),\]
with differential $\d_{\pi}:\X^{\bullet}(M)\to\X^{\bullet+1}(M)$ given
by taking the Schouten bracket with $\pi$:
\[ \d_{\pi} \theta=[\pi,\theta]. \]
Moreover, verify the following interpretations of the first few
cohomology groups: 
\begin{itemize}
\item $H^0_\pi(M)$ is formed by the functions in the center of the
  Poisson algebra $(C^\infty(M),\{~,~\})$ (often called \emph{Casimirs}).
\item $H^1_\pi(M)$ is the quotient of the Poisson vector fields
  by the hamiltonian vector fields.
\item $H^2_\pi(M)$ is the space of classes of non-trivial infinitesimal
 deformations of the Poisson structure $\pi$. 
\end{itemize}
\end{exercise}

Given a vector field $X\in\X(M)$ and a cotangent path $a:I\to T^*M$ we
define the \emph{integral} of $X$ along $a$ by:
\[ \int_a X=\int_0^1 \langle a(t),X|_{\gamma(t)}\rangle\, \d t,\]
where $\gamma:I\to M$ is the base path of $a$. Notice that if $X_h\in
\X_{\Ham}(M)$ is a hamiltonian vector field, then:
\[ \int_a X_h=h(\gamma(1))-h(\gamma(0)),\]
so the integral in this case depends only on the end-points.

Though the definition of the integral along cotangent paths makes
sense for \emph{any} vector field, we will only be interested in the
integral along \emph{Poisson vector fields}. The reason is that
homotopy invariance holds only for Poisson vector fields, so that one
should think of Poisson vector fields as the analogue of closed
1-forms in this contravariant calculus:

\begin{proposition}
  \label{prop:ctg:homotopy}
  Let $a_{\epsilon}\in P_\pi(M)$ be a family of cotangent paths, whose
  base paths have fixed end-points. Then $a_{\epsilon}$ is a cotangent
  homotopy iff for all $X\in\X_{\pi}(M)$
  \[ \frac{\d}{\d\epsilon} \int_{a_{\epsilon}} X=0.\]
\end{proposition}

\begin{proof} 
  Let $a=a(\epsilon,t)$ be a variation of $A$-paths, fix a connection
  $\nabla$ on $T^*M$ and let $b=b(\epsilon,t)$ be a solution of
  equation (\ref{diffeq}). If we set $I=\langle a,X\rangle$ and
  $J=\langle b,X\rangle$, then a straightforward computation
  using the defining equation (\ref{diffeq}) and the expression
  (\ref{kozul:bracket}) for the Lie bracket, shows that
  \[ \frac{\d I}{\d\epsilon}- \frac{\d J}{\d t}= \Lie_X{\pi}(a,b) .\]
  Integrating first with respect to $t$, using $b(\epsilon, 0)= 0$, and
  then integrating with respect to $\epsilon$, we find that
  \[ \int_{a_{1}}X- \int_{a_{0}}X=
  \int_{b(\cdot,1)} X + \int_{0}^{1}\int_{0}^{1}
  \Lie_X \pi (a, b)\, \d t\,\d\epsilon,\] 
  for any vector field $X$. When $X\in\X_\pi(M)$, the second term on
  the right-hand side drops out. 

  Now if $a_{\epsilon}$ is a cotangent homotopy, then
  $b(\epsilon,1)=0$ and we find that $\int_{a_{0}}X=\int_{a_{1}}X$. A 
  scaling argument shows that $\int_{a_{0}}X=\int_{a_{\epsilon}}X$ for
  any epsilon, so that:
   \[ \frac{\d}{\d\epsilon} \int_{a_{\epsilon}} X=0.\]
  Similarly, if this holds for any $X\in\X_\pi(M)$, we conclude that
  $b(\epsilon,1)=0$ so that $a_{\epsilon}$ is a cotangent homotopy.
\end{proof}

We will denote by $\Sigma(M,x)$ the isotropy group at $x$, which is
the topological group formed by all cotangent homotopy classes of loops
based at $x\in M$ (in the notation of the previous lectures, this is
just $\G_x(T^*M)$). Also, we will denote by $\Sigma(M,x)^0$ the restricted
isotropy group, formed by classes of cotangent loops covering loops that
are contractible in the symplectic leaf. Then we have the short exact sequence:
\[ 1\rmap \Sigma(M,x)^0\rmap \Sigma(M,x)\rmap \pi_1(S,x)\rmap 1,\]
where $S$ is the symplectic leaf through $x$.

\begin{exercise}
  Show that integration gives a homomorphism:
  \[ \int:\Sigma(M,x)\to H^1_{\pi}(M)^*.\]
\end{exercise}

All these sugest to think of $\Sigma(M,x)$ as the Poisson
fundamental group of $M$ at $x$. We remind that these groups will be,
in general, non-isomorphic as $x$ varies (but, of course, $\Sigma(M,x)$ is
isomorphic to $\Sigma(M,y)$ if $x$ and $y$ belong to the same
symplectic leaf). We should also keep in mind that these groups,
in general, are not discrete. In fact, we know that we have an
isomorphism:
\[ \Sigma(M,x)^0\simeq \G(\gg_x)/\tilde{\NN}_x,\]
where $\tilde{\NN}_x$ is the monodromy group of our Poisson manifold at
$x$ and $\gg_x$ is the isotropy Lie algebra at $x$. 

As one could expect, for a Poisson manifold $M$ the isotropy Lie
algebra and the monodromy groups have nice geometric interpretations. 

First, recall that if the Poisson structure vanishes at a point $x_0$
(i.e., $\pi|_{x_0}=0$), then there is a well-defined linear
approximation to $\pi$ at $x_0$, which is a certain linear Poisson
structure on the tangent space $T_{x_0}M$. Equivalently, we have the
Lie algebra structure on $T^*_{x_0}M$ defined by:
\[ [\d_{x_0}f,\d_{x_0}g]=\d_{x_0}\{f,g\}.\]
At a point $x\in M$ of higher rank, the Weinstein splitting theorem
shows that any small transverse $N$ to the symplectic leaf $S$ through
$x$, inherits a Poisson structure which vanishes at $x$. This gives a
linear Poisson structure on $T_x N$, or equivalently, a Lie algebra
structure on $T_x^* N$. A diferent choice of transversal produces
isomorphic linear approximations. It is easy to check that on the
normal space $\nu_x(S)$ we obtain an intrinsically defined Poisson
structure, and so also an intrinsically defined Lie algebra structure
on the dual $\nu^*_x(S)$. This just the isotropy Lie algebra at $x$:
\[ \gg_x=\nu^*_x(S).\]

Now, let us turn to a geometric interpretation of the monodromy. We
will assume that our Poisson manifold is \emph{regular}(\footnote{By
regular, we mean that the rank of the Poisson structure $\pi$ is
constant, so the symplectic leaves form a regular foliation $\F$. If
you prefer, you can assume that we are restricting to a neighborhood
of a regular leaf.}) so both $T\F=\Im\pi^{\sharp}$ and
$\nu^*(\F)=\ker\#$ are vector bundles over $M$, and the isotropy Lie
algebras $\nu^*(\F)_{x}$ are all abelian. This amounts to several
simplifications. For instance, the long exact sequence of of the
monodromy reduces to
\[ 
\xymatrix{\cdots\ar[r]&\pi_{2}(S,x)\ar[r]^{\partial}& \nu^*(S)_x
\ar[r]& \Sigma(M, x)\ar[r]& \pi_1(S,x) }.
\] 
The monodromy groups can be described (or defined) as the image of
$\partial$, which, in turn, is given by integration of a canonical
cohomology 2-class $[\Omega_{\sigma}]\in H^2(S,\nu_{x}^{*}(S))$. This
class can be computed explicitly by using a section $\sigma$ of
$\pi^\sharp:T^{*}_{S}M\to TS$.

Fix a point $x$ in a Poisson manifold $M$, let $S$ be the symplectic
leaf through $x$, and consider a 2-sphere $\gamma: S^2\to S$, which
maps the north pole $p_{N}$ to $x$. The symplectic area of $\gamma$ is
given, as usual, by
\[ A_{\omega}(\gamma)= \int_{S^2} \gamma^{*}\omega ,\]
where $\omega$ is the symplectic 2-form on the leaf $S$. By a
\textbf{deformation} of $\gamma$ we mean a family $\gamma_{t}: S^2\to M$
of 2-spheres parameterized by $t\in (-\eps,\eps)$, starting at
$\gamma_{0}= \gamma$, and such that for each fixed $t$ the sphere
$\gamma_t$ has image lying entirely in a symplectic leaf. The
\textbf{transversal variation} of $\gamma_t$ (at $t=0$) is the class of
the tangent vector
\[
\var_{\nu}(\gamma_{t})\equiv
\left[\left.\frac{\d}{\d t}\gamma_{t}(p_{N})\right|_{t=0}\right]\in
\nu_{x}(\F).
\]
We shall see below that the quantity
\[ \left.\frac{\d}{\d t}A_{\omega}(\gamma_{t})\right|_{t=0} \]
only depends on the homotopy class of $\gamma$ and on
$\var_{\nu}(\gamma_{t})$. Finally, the formula
\[ \langle A^{'}_{\omega}(\gamma), var_{\nu}(\gamma_{t})\rangle=
\left.\frac{\d}{\d t}A_{\omega}(\gamma_{t})\right|_{t=0},\] 
applied to different deformations of $\gamma$, gives a well defined
element
\[ A^{'}_{\omega}(\gamma)\in \nu_{x}^{*}(S) .\]

\vskip 15 pt
\begin{center}
\includegraphics[width=.8\textwidth]{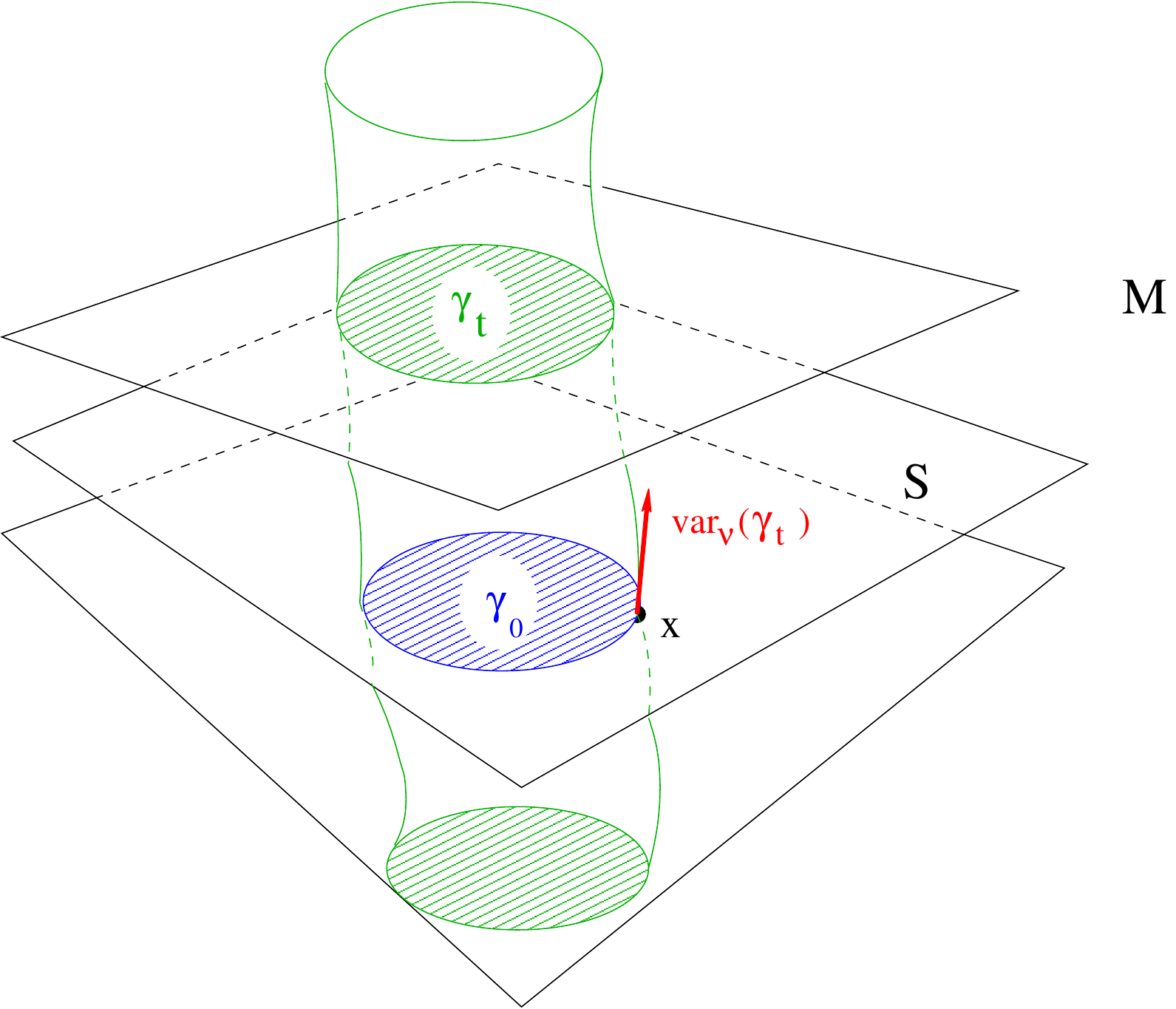}
\end{center}

Now, we have:

\begin{proposition}
  \label{var-area}
  For any regular manifold $M$, 
  \[ \NN_{x}= \{A^{'}_{\omega}(\sigma): \sigma\in \pi_{2}(S,x)\},\]
  where $S$ is the symplectic leaf trough $x$.
\end{proposition}

\begin{proof}
  Recall first that the normal bundle $\nu$ (hence also any associated
  tensor bundle) has a natural flat $\F$-connection
  $\nabla:\Gamma(T\F)\times \Gamma(\nu)\to\Gamma(\nu)$, given by
  \[ \nabla_{X} \overline{Y}= \overline{[X, Y]}.\]
  In terms of the Lie algebroid $T\F$, $\nabla$ is a flat Lie
  algebroid connection giving a Lie algebroid representation of $T\F$
  on the vector bundle $\nu$ (and also, on any associated tensor
  bundle). Therefore, we have the foliated cohomology with
  coefficients in $\nu$, denoted $H^*(\F;\nu)$ (see Exercise
  \ref{ex:cohomology}). Similarly one can talk about cohomology with
  coefficients in any tensorial bundle associated to $\nu$. In the
  special case of trivial coefficients, we recover $H^*(\F)$.
  
  Now, the Poisson tensor in $M$ determines a foliated 2-form
  $\omega\in\Omega^2(\F)$, which is just another way of looking at the
  symplectic forms on the leaves. Therefore, we have a foliated
  cohomology class in the second foliated cohomology group:
  \[ [\omega]\in H^2(\F).\]
  On the other hand, as we saw in Lecture \ref{integrability:I}, a
  splitting $\sigma$ determines a foliated 2-form
  $\Omega_\sigma\in\Omega^2(\F;\nu^*)$, with coefficients in the
  co-normal bundle, and the corresponding foliated cohomology
  class
  \[ [\Omega]\in H^2(\F; \nu^*) \]
  does not depend on the choice of splitting.
  
  These two classes are related in a very simple way: there is
  a map
  \[ \d_{\nu}: H^2(\F)\rmap H^2(\F; \nu^*),\]
  which can be described as follows. We start with a class $[\theta]\in
  H^2(\F)$, represented by a foliated 2-form $\theta$.  As with any
  foliated form, we have $\theta= \widetilde{\theta}|_{T\F}$ for some
  2-form $\widetilde{\theta}\in\Omega^2(M)$. Since
  $\d\widetilde{\theta}|_{\F}= 0$, it follows that the map
  $\Gamma(\wedge^2\F)\to\Gamma(\nu^*)$ defined by
  \[ (X,Y)\mapsto d\widetilde{\theta}(X, Y, -), \]
  gives a closed foliated 2-form with coefficients in $\nu^*$. It is
  easily seen that its cohomology class does not depend on the choice of
  $\widetilde{\theta}$, and this defines $\d_{\nu}$. Now the formula
  $\omega(X_{f},X)= X(f)$, immediately implies that
  \[ \d_{\nu}([\omega])= [\Omega].\]
  Notice that this construction of $\d_{\nu}$ is functorial with respect
  to foliated maps (i.e., maps between foliated spaces which map leaves
  into leaves).

  {From} this perspective, a deformation $\gamma_{t}$ of 2-spheres is
  a foliated map $S^2\times I\to M$, where in $S^2\times I$ we
  consider the foliation $\F_{0}$ whose leaves are the spheres
  $S^2\times\{t\}$. {From} $H^2(S^2)\simeq \mathbb{R}$, we get
  \begin{align*}
    H^{2}(\F_{0})&\simeq C^{\infty}(I)=\Omega^0(I),\\
    H^{2}(\F_{0};\nu)&\simeq C^{\infty}(I) dt=\Omega^1(I),
  \end{align*}
  where the isomorphisms are obtained by integrating over $S^2$.
  Hence, $d_{\nu}$ for $\F_{0}$ becomes the de Rham differential
  $d:\Omega^0(I)\to\Omega^1(I)$. Now, the functoriality of $d_{\nu}$
  with respect to $\gamma_{t}$ when applied to $\omega$ gives
  \[ \frac{\d}{\d t} \int_{S^2}\gamma_{t}^{*}\omega=
  \langle \int_{\gamma_{t}}\Omega,
  \frac{\d}{\d t}\gamma_{t}(p_{N})\rangle.\] 
  This proves the last part of the proposition (and also the
  properties of the variation of the symplectic area stated before). 
\end{proof}

\section{Symplectic groupoids}

Let us now turn to the study of the groupoid
\[ \Sigma(M)=P_\pi(M)/\sim.\]
It is probably a good idea to recall now a few basic examples.

\begin{example}[symplectic manifolds]
  Let $M=S$ be a symplectic manifold with symplectic form
  $\omega=\pi^{-1}$. The symplectic form gives an inverse
  $\omega^{\sharp}:TS\to T^*S$ to the anchor $\pi^{\sharp}:T^*S\to
  TS$. Hence, the cotangent Lie algebroid $T^*S$ is isomorphic to the
  tangent Lie algebroid $TS$ and a cotangent path is determined by its base
  path.
  
  It follows that, if we assume $S$ to be simply connected, then
  $\Sigma(S)$ is just the pair groupoid $S\times S$. Note that
  $\Sigma(S)$ becomes a symplectic manifold with symplectic form
  $\omega\oplus(-\omega)$ (the reason for choosing a minus sign in the
  second factor will be apparent later).
\end{example}

\begin{exercise}
  Extend this example to a non-simply connected symplectic manifold,
  showing that $\Sigma(S)$ still is a symplectic manifold in this case.
\end{exercise}

\begin{example}[trivial Poisson manifold]
  At the other extreme, let $M$ be any manifold with the trivial bracket
  $\{~,~\}\equiv 0$. Then $T^*M$ becomes a bundle of abelian Lie
  algebras. Cotangent paths are just paths in the fibers of $T^*M$ and
  any such path is cotangent homotopic to its average. It follows that
  $\Sigma(M)=T^*M$, with multiplication being addition on the
  fibers. Note that $\Sigma(M)$ is a symplectic manifold with the
  usual canonical symplectic form on the cotangent bundle.
\end{example}

\begin{example}[linear Poisson structures]
  Let $\gg$ be any finite dimensional Lie algebra, and take $M=\gg^*$
  with the Kostant-Kirillov-Souriau bracket:
  \[ 
  \{f_1,f_2\}(\xi)=\langle \xi,[\d_{\xi}f_1,\d_{\xi}f_2]\rangle,
  \quad f_1,f_2\in C^\infty(M),\ \xi\in\gg^*.
  \]
  If $G$ is the 1-connected Lie group integrating $\gg$, we have seen
  before that the cotangent Lie algebroid of $\gg^*$ integrates to the
  1-connected Lie groupoid $\G=T^*G$. By uniqueness, we conclude that
  $\Sigma(\gg^*)=T^*G$. Note again that $\Sigma(\gg^*)$ carries a
  symplectic structure.
\end{example}

In all these examples, the groupoid $\Sigma(M)$ that integrates the
Poisson manifold $M$ carries a symplectic form. It turns out that the
symplectic form is compatible with the groupoid structure in the
following sense:

\begin{definition}
  Let $\G$ be a Lie groupoid and $\omega\in\Omega^\bullet(\G)$. We call
  $\omega$ a \textbf{multiplicative form} if
    \[m^*\omega=p^*_1\omega+p^*_2\omega,\]
  where we denote by $m:\G_2\to\G$ the multiplication in $\G$ defined
  on the submanifold $\G_2\subset \G\times \G$ of composable arrows,
  and by $p_1,p_2:\G_2\to\G$ the projections to the first and
  second factors.
  A \textbf{symplectic groupoid} is a pair $(\G,\omega)$ where $\G$ is
  a Lie groupoid and $\omega\in\Omega^2(\G)$ is a multiplicative
  symplectic form.
\end{definition}

\begin{exercise}
  Check that in all examples above the symplectic forms are
  multiplicative, so all those Poisson manifolds integrate to
  $s$-simply connected symplectic groupoids.
\end{exercise}

In the sequel, if $(S,\omega)$ is a symplectic manifold, we will
denote by $\bar{S}$ the manifold $S$ furnished with the
symplectic form $-\omega$.

\begin{exercise}
  Recall that a map $\phi:(S_1,\omega_1)\to (S_2,\omega_2)$ is
  symplectic if $\phi^*\omega_2=\omega_1$. Show that:
  \begin{enumerate}[(i)]
  \item A map $\phi$ is symplectic iff its graph is a Lagrangian
    submanifold of $S\times\bar{S}$.
  \item A symplectic form $\omega$ in a groupoid $\G$ is
    multiplicative iff the graph of the multiplication 
    \[ \gamma_m=\{(g,h,g\dot h)\in \G\times\G\times\G:(g,h)\in\G^{(2)}\}\] 
    is a Lagrangian submanifold of $\G\times\G\times\bar{\G}$.
  \end{enumerate}
\end{exercise}

Here are some basic properties of a symplectic groupoid, whose proof is
a very instructive exercise.

\begin{exercise}
  Let $(\G,\omega)$ be a symplectic groupoid. Show that:
  \begin{enumerate}[(i)]
  \item the $\s$-fibers and the $\t$-fibers are symplectic orthogonal;
  \item the inverse map $i:\G\to\G$ is an anti-symplectic map;
  \item $M$, viewed as the unit section, is a Lagrangian submanifold.
  \end{enumerate}
\end{exercise}

Finally, we can state and prove the main property of a symplectic
groupoid: the base manifold of a symplectic groupoid has a canonical
Poisson bracket.

\begin{proposition}
  Let $(\G,\omega)$ be a symplectic groupoid over $M$. Then, there
  exists a unique Poisson structure on $M$ such that:
  \begin{enumerate}[(i)]
  \item $\s$ is Poisson and $\t$ is anti-Poisson;
  \item the Lie algebroid of $\G$ is canonically isomorphic to $T^*M$.
  \end{enumerate}
\end{proposition}

\begin{proof}
  Let us call a function on $\G$ of the form $f\circ\s$ a \emph{basic
  function}.  We can assume that $\G$ is $s$-connected, so a function
  $f:\G\to\Rr$ is basic iff it is constant on $s$-fibers. Then, using
  the exercise above, we have the chain of equivalences:
  \begin{align*} 
  f \text{ is basic } &\Longleftrightarrow \d f\in(\Ker \d\s)^0\\
  &\Longleftrightarrow X_f\in(\Ker \d\s)^\perp
  \Longleftrightarrow X_f\in\Ker \d\t.
  \end{align*}
  From this last characterization and the fact that
  \[ X_{\{f_1,f_2\}_\G}=[X_{f_1},X_{f_2}],\]
  we see that the Poisson bracket of any two basic functions is
  a basic function. 

  We conclude that $M$ carries a unique Poisson bracket such that $\s$
  is a Poisson map. Since $\t=\s\circ i$, and the inversion is
  anti-symplectic, we see that $\t$ is anti-Poisson. Now we use the
  following: 

  \begin{exercise}
    Show that if $f\in C^\infty(M)$, then the Hamiltonian vector field
    $-X_{f\circ\t}\in\X_{\Ham}(\G)$ is the right invariant vector field
    in $\G$ that corresponds to the section $\d f\in\Omega^1(M)$.
  \end{exercise}
  
  Hence, to compute the anchor in $A(\G)$ we start with an exact 1-form
  $\d f\in \Omega^1(M)$, we extend it to a right invariant
  vector field $-X_{f\circ\t}$, and apply the differential of the
  target map, obtaining:
  \[ \rho(\d f)=\d\t (-X_{f\circ t})=X_f.\]
  This shows that the anchor for $A(\G)$ coincides with the anchor for $T^*M$. 

  On the other hand, to compute the $A(\G)$-bracket of two exact 1-forms
  $\d f_i\in \Omega^1(M)$, we extend them to the right invariant vector
  fields $-X_{f_i\circ\t}$ and we compute their Lie bracket:
  \begin{align*} 
    [-X_{f_1\circ\t},-X_{f_2\circ\t}]&=X_{\{f_1\circ\t,f_2\circ\t\}_\G}\\
    &=-X_{\{f_1,f_2\}_M\circ\t}.
  \end{align*}
  Therefore, the $A(\G)$-bracket is given by $\d\{f_1,f_2\}_M$,
  which coincides with the Koszul bracket $[\d f_1,\d f_2]$. Since we
  already know that the anchors coincide, the Leibniz identity implies
  that the Lie brackets on $A(\G)$ and $T^*M$ coincide on \emph{any}
  1-forms.
\end{proof}

In the examples of Poisson manifolds that we saw before, $\Sigma(M)$
was a symplectic groupoid. This situation is by no means
excepcional: we will now show that whenever $\Sigma(M)$ is smooth,
it is a symplectic groupoid!

In order to prove this, we will need to look closer at the description
of cotangent homotopies through a Lie algebra action (we saw this in
the previous lecture for general $A$-homotopies). We denote now by
$P(T^*M)$ the space of $C^1$-paths $a:I\to T^*M$ with base path
$\gamma:I\to M$ of class $C^2$. This is a Banach manifold in the
obvious way, and the space of cotangent paths $P_\pi(M)\subset
P(T^*M)$ is a Banach submanifold. Also, we let $P(M)$ denote the space
of $C^2$-paths $\gamma:I\to M$. A basic fact, which in fact explains the
existence of a symplectic structure on $\Sigma(M)$, is that 
\[ P(T^*M)=T^*P(M), \]
so that $P(T^*M)$ carries a natural (weak) symplectic structure
$\omega_{\text{can}}$. 

\begin{exercise}
Identify the tangent space $T_{a}P(T^*M)$ with the space of vector fields
along $a$:
\[ T_{a}P(T^*M)=\set{U:I\to TT^*M: U(t)\in T_{a(t)}T^*M}.\]
and denote by $\omega_0$ the canonical symplectic form on $T^*M$. Show that
the 2-form $\omega_{\text{can}}$ on $P(T^*M)=T^*P(M)$ is given by:
\begin{equation}
  \label{eq:lifted:symp:form}
  (\omega_{\text{can}})_a(U_1, U_2)=\int_0^1 \omega_0(U_1(t), U_2(t))dt,
\end{equation}
for all $U_1, U_2\in T_a P(T^*M)$. Moreover, check that
$\d\omega_{\text{can}}=0$ and that
$\omega_{\text{can}}^\sharp:T_{a}P(T^*M)\to T_{a}^*P(T^*M)$ is 
injective, so that $\omega_{\text{can}}$ is a weak symplectic form. 
\end{exercise}

On the other hand, the Lie algebra
\[
P_{0}\Omega^1(M):=\set{\eta_{t}\in \Omega^1(M), t\in I:\eta_{0}
                =\eta_{1}=0,\ \eta_t\text{ of class }C^1\text{ in }t}
\]
with the pointwise Lie bracket, acts on $P(T^*M)$ in such a way that:
\begin{enumerate}[(a)]
\item the action is tangent to $P_\pi(M)$;
\item two cotangent paths are homotopic if and only if they belong to
the same orbit.
\end{enumerate}
Now we have the following remarkable fact:

\begin{theorem}
  \label{thm:sigma:model}
  The Lie algebra action of $P_{0}\Omega^1(M)$ on
  $(P(T^*M),\omega_{\text{can}})$ is Hamiltonian, with equivariant
  moment map $J:P(T^*M)\to P_{0}\Omega^1(M)^*$ given by
  \begin{equation}
  \label{eq:moment:map}
  \langle J(a),\eta\rangle=
  \int_0^1\langle\frac{d}{dt}\pi(a(t))-\#a(t),\eta(t,\gamma(t))\rangle\,dt.
\end{equation}
\end{theorem}

\begin{proof}
  We need to check that if $\eta_t\in P_{0}\Omega^1(M)$ is a time
  dependent 1-form vanishing at $t=0,1$, the the infitnitesimal
  generator $X_\eta$ coincides with the Hamiltonian vector field
  with Hamiltonian function $a\mapsto \langle J(a),\eta\rangle$. We
  leave this for the reader to check. (Hint: It is enough to prove
  this at $A$-paths $a$ fitting in a coordinate domain, so one can use
  local coordinates.)
\end{proof}

Note that set of cotangent paths $P_{\pi}(M)$ is just the level set
$J^{-1}(0)$. Hence, our groupoid $\Sigma(M)$ can be described
alternatively as a symplectic quotient:
\begin{equation}
  \label{eq:sigma:model}
  \Sigma(M)=\widetilde{P}(T^*M)//P_{0}\Omega(M).
\end{equation}
and we have:

\begin{corollary}
  Whenever $\Sigma(M)$ is smooth, it admits a canonical multiplicative
  symplectic form $\omega$, so that $(\Sigma(M),\omega)$ is a
  symplectic groupoid. Moreover, the induced Poisson bracket on the
  base is the original one.
\end{corollary}

\begin{proof}[Proof of Theorem \ref{thm:sigma:model}]
  We only need to check that the reduced symplectic form $\omega$ is
  multiplicative, i.e., satisfies
  \[ m^*\omega=\pi^*_1\omega+\pi^*_2\omega,\]
  where $m:\Sigma(M)_2\to\Sigma(M)$ is the multiplication in
  $\Sigma(M)$ and $\pi_1,\pi_2:\Sigma(M)_2\to\Sigma(M)$ are the
  projections to the first and second factors. But the additivity of
  the integral and expression (\ref{eq:lifted:symp:form}) shows that
  that this condition holds already at the level of $P(T^*M)$, hence
  it must hold also on the reduced space $\Sigma(M)$.
\end{proof}

Putting all these together, and usings the general integrability
results of Lecture \ref{integrability:II}, we arrive at the following
esult:

\begin{theorem}
  \label{thm:main-Poisson} 
  For a Poisson manifold $M$, the following are equivalent:
  \begin{enumerate}[(i)]
  \item $M$ is integrable by a symplectic groupoid.
  \item The algebroid $T^*M$ is integrable.
  \item The groupoid $\Sigma(M)$ is a smooth manifold.
  \item The monodromy groups $\NN_{x}$, with $x\in M$, are locally
    uniformly discrete. 
  \end{enumerate}
  In this case, $\Sigma(M)$ is the unique $s$-simply connected,
  symplectic groupoid which integrates $M$.
\end{theorem}

\begin{proof}
  If $\Sigma$ is a symplectic groupoid integrating $M$, then its
  associated algebroid is isomorphic to $T^*M$, and this shows that (i)
  implies (ii). We have seen in Lecture 2 that (ii) implies (iii), by
  taking any Lie groupoid $\G$ integrating $T^*M$, and then taking the
  groupoid formed by the universal covers of the $\s$-fibers of $\G$
  together. Obviously (iii) implies (i), since we just saw that
  $\Sigma(M)$ is a symplectic groupois. Finaly, the equivalence of (iii)
  and (iv) is just a special case of our main integrability theorem for
  Lie algebroids.
\end{proof}

Note that not every Lie groupoid integrating the cotangent bundle
$T^*M$ of a Poisson manifold is a symplectic groupoid, as the
following example shows:

\begin{example}
  Take $M=\Ss^3$ with the trivial Poisson bracket, so that the $s$-simply
  connected groupoid integrating $M$ is $\Sigma(\Ss^3)=T^*\Ss^3$,
  multiplication being addition on the fibers. Another
  integrating groupoid is, for example, $\G=\Ss^3\times \Tt^3$ where $\s=\t$
  is projection in the first factor and multiplication is
  multiplication on the torus $\Tt^3$. This groupoid cannot be a
  symplectic groupoid because of the following exercise:

\begin{exercise}
  Show that the manifold $\Ss^3\times\Tt^3$ does not admit any
  symplectic structure.
\end{exercise}

\end{example}

At this point it is natural to ask, for a general Lie algebroid $A$, what is
the relationship between the integrability of $A$ and the
integrability of the Poisson manifold $A^*$.  The following exercise
shows that these two problems are actually equivalent.

\begin{exercise}
  Let $A$ be a Lie algebroid and consider the natural fiberwise linear
  Poisson structure on $A^*$. Check that the natural injection
  $A\hookrightarrow T^*A^*$ is a Lie algebroid morphism, so that:
  \begin{enumerate}[(a)]
  \item If $A$ is integrable, then the Poisson manifold $A^*$
    is also integrable and $\Sigma(A^*)=T^*\G(A)$.
  \item If $A^*$ is integrable, then $A$ is also integrable and
    $\G(A)\hookrightarrow\Sigma(A^*)$ is a subgroupoid.
  \end{enumerate}
\end{exercise}

This exercise is one instance of the intrinsic connection between Poisson
and symplectic geometry on one hand, and Lie algebroid and groupoid theory
on the other hand. In the next section, we will go deeper into this connection
which is sometimes not so obvious.

It is now time to look at some examples.

\begin{example}[Poisson manifolds of dimension 2]
\label{2-dim} 
The lowest dimension one can have non-trivial Poisson manifolds is
2. However, such a Poisson manifold will have it follows immediately 
from Theorem \ref{thm:main-Poisson} and the description of the
monodromy groups that in dimension 2 all Poisson manifolds are integrable:

\begin{corollary}
\label{cor:2-dim} Any 2-dimensional Poisson manifold is
integrable.
\end{corollary}

Corollary \ref{cor:2-dim} can be partially generalized to higher dimensions
in the following sense: any $2n$-dimensional Poisson manifold whose
Poisson tensor has rank $2n$ on a dense, open set, is
integrable. The proof of this fact is more involved (see the notes at the
end of this chapter).
\end{example}

\begin{example}[A non-integrable Poisson manifold]
\label{example:basic} 
Already in dimension 3 there are examples of non-integrable 
Poisson manifolds. Let us consider $M=\mathbb{R}^3$ with the 
Poisson bracket:
\begin{equation}
\label{eq:PB:su2}
 \{f,g\}= \det
\begin{pmatrix}
x & y & z \\
\frac{\partial f}{\partial x}&\frac{\partial f}{\partial y}&
\frac{\partial f}{\partial z}\\
\frac{\partial g}{\partial x}&\frac{\partial g}{\partial y}&
\frac{\partial g}{\partial z}
\end{pmatrix}.
\end{equation}
Notice that this is just the linear Poisson structure on $\mathfrak{su}(2)^*$
when we identify $\mathfrak{su}(2)$ and $\mathbb{R}^3$ with the exterior
product. 

Let us choose any smooth function $a=a(R)$ on $M$, 
which depends only on the radius $R$, and which is strictly positive for $R>0$. 
We multiply the previous brackets by $a$, and we denote by $M_{a}$
the resulting Poisson manifold. The bracket on $\Omega^1(M_a)$ is
computed using the Leibniz identity and we get
\begin{align*}
[dx^2,dx^3]&=a dx^1+ b x^1 R \bar{n},\\
[dx^3,dx^1]&=a dx^2+ b x^2 R \bar{n},\\
[dx^1,dx^2]&=a dx^3+ b x^3 R \bar{n},
\end{align*}
where $\bar{n}= \frac{1}{R}\sum_i x^idx^i$ and $b(R)= a'(R)/R$.
The bundle map $\#: T^*M_{a}\to TM_{a}$ is just
\[ \#(dx^i)=a \bar{v}^{i}, \qquad i=1,2,3 \]
where $\bar{v}^i$ is the infinitesimal generator of a rotation about the
$i$-axis:
\[
\bar{v}^1=x^3\frac{\partial}{\partial x^2}-x^2\frac{\partial}{\partial
x^3}, \quad \bar{v}^2=x^1\frac{\partial}{\partial
x^3}-x^3\frac{\partial}{\partial x^1}, \quad
\bar{v}^3=x^2\frac{\partial}{\partial x^1}-x^1\frac{\partial}{\partial
x^2}.
\]
The leaves of the symplectic foliation of $M_{a}$ are the spheres
$S_{R}^{2}\subset\mathbb{R}^3$ centered at the origin, and the
origin is the only singular point. 

We will compute the monodromy function $r_{\NN}$ in two distict fashions.
First, using the obvious metric on $T^*M_{a}$, we restrict to
a leaf $S_{R}^{2}$ with $R>0$, and we choose as a splitting of 
$\#$ the map defined by
\[ \sigma(\bar{v}^i)= \frac{1}{a}(dx^i-\frac{x^i}{R}\bar{n}),\]
with curvature the center-valued 2-form
\[ \Omega_{\sigma}= \frac{Ra'- a}{a^2R^3}\omega \bar{n},\]
where $\omega= x^1 dx^2\wedge dx^3+ x^2dx^3\wedge dx^1+
x^3dx^1\wedge dx^2$. Since $\int_{S^2_R}\omega=4\pi R^3$ it
follows that
\[
\NN_{(x, y, z)}\simeq 4\pi \frac{Ra'- a}{a^2}\Zz\bar{n}\subset
\Rr\bar{n}.
\]

On the other hand, the canonical generator of $\pi_{2}(S^{2}_{R})$ 
defines \emph{the} symplectic area of $S^2_R$, which is easily computable:
\[ A_{a}(R)= 4\pi \frac{R}{a(R)}.\]
We recover in this way the relationship between the monodromy and
the variation of the symplectic area (Proposition \ref{var-area}).

Also, observe that
\[
r_\NN(x,y,z)= \left\{
\begin{array}{ll}
+\infty \qquad&\  \text{if}\ R=0\ \text{ or }  A'_{a}(R)= 0,\\
\\
A'_{a}(R) \qquad&\  \text{otherwise,}
\end{array}
\right.
\]
so the monodromy might vary in a non-trivial fashion, even for nearby
regular leaves. Our computation also gives the isotropy groups
\[
\Sigma(M_{a}, (x, y, z))\cong \left\{
\begin{array}{llll}
\mathbb{R}^3 \qquad&\  \text{if}\ R=0,\ a(0)= 0,\\
SU(2) \qquad&\ \text{if}\ R=0,\ a(0)\neq 0,\\
\mathbb{R} \qquad&\  \text{if}\ R\neq 0,\  A'_{a}(R)= 0,\\
S^1 \qquad&\  \text{if}\ R\neq 0,\ A'_{a}(R)\neq 0.
\end{array}
\right.
\]
\end{example}

\begin{example}[Heisenberg-Poisson manifolds]
\label{Heisenberg} 
The \textbf{Heisenberg-Poisson manifold} $M(S)$
associated to a symplectic manifold $S$, is the manifold 
$S\times\mathbb{R}$ with the Poisson structure given by $\{f, g\}= t\{f_{t},
g_{t}\}_S$, where $t$ stands for the real parameter, and $f_{t}$
denotes the function on $S$ obtained from $f$ by fixing the value of
$t$. We have the following result:

\begin{corollary} 
\label{cor:heisenberg}
For a symplectic manifold $S$, the following are equivalent:
\begin{enumerate}[(i)]
\item The Poisson-Heisenberg manifold $M(S)$ is integrable;
\item $\widetilde{S}$ is pre-quantizable.
\end{enumerate}
\end{corollary}

We recall that condition (ii) is usually stated as follows: when
we  pull back the symplectic form $\omega$ on $S$ to a 2-form
$\widetilde{\omega}$ on the covering space $\widetilde{S}$, the
group of periods
\[ \set{\int_\gamma
  \widetilde{\omega}:\gamma\in\H_2(\widetilde{S},\Zz)}\subset \Rr\]
is a multiple of $\mathbb{Z}$. Note that this group coincides with the
\textbf{group of spherical periods} of $\omega$
\[ \P(\omega)=\set{\int_\gamma \omega:\gamma\in\pi_2(S)},\]
so that (ii) says that $\P(\omega)\subset\Rr$ is a multiple of
$\Zz$.

\begin{proof}[Proof of Corollary \ref{cor:heisenberg}]
We have to compute the monodromy groups. The singular symplectic
leaves are the points in $S\times \{0\}$ and they clearly have
vanishing monodromy groups. The regular symplectic leaves are the
submanifolds $S\times \{t\}$, where $t\neq 0$, with symplectic
form $\omega/t$. To compute their monodromy groups we invoke 
Proposition \ref{var-area}, so that we immediately get
\[ \NN_{(x, t)}= \frac{1}{t}\P(\omega)\subset \mathbb{R},\]
and the result follows.
\end{proof}

\begin{exercise}
Show that the arrows in $\Sigma(M(S))$ are of two types:
\begin{enumerate}[(a)]
\item arrows which start and end at $(x, 0)$, which form a group
isomorphic to the additive subgroup of $T_{x, 0}^{*}M(S)$;
\item arrows inside the symplectic leaves $S\times \{t\}$, $t\neq
0$, which consist of equivalence classes of pairs $(\gamma, v)$, where
$\gamma$ is a path in $S$ and $v\in \mathbb{R}$. Two such pairs
$(\gamma_i,v_i)$ are equivalent if and only if there is a homotopy
$\gamma(\eps,s)$ (with fixed end points) between the $\gamma_{i}$'s,
such that $v_1-v_0=\frac{1}{t}\int\gamma^*\omega$.
\end{enumerate}
\end{exercise}
\end{example}

\begin{remark}
The two necessary and sufficient conditions that must hold 
in order for $\Sigma(M)$ to be smooth, namely:
\begin{enumerate}[(i)]
\item $r_\NN(x)>0$, and
\item $\liminf_{y\to x}r_\NN(y)>0$,
\end{enumerate}
are independent. In fact, if we choose a symplectic manifold 
which does not satisfy the pre-quantization condition, the 
Poisson-Heisenberg manifold of Example \ref{Heisenberg} gives a non-integrable Poisson manifold 
violating condition (i). On the other hand, Example \ref{example:basic} 
gives examples of non-integrable Poisson manifolds in which
condition (i) is satisfied, but condition (ii) is not.
\end{remark}

\section{The symplectization functor}

So far, to every Poisson manifold $(M,\pi)$ we have associated a groupoid 
$\Sigma(M)$. In the integrable case, this is a symplectic groupoid 
integrating the cotangent algebroid $T^*M$. If  
$\phi:(M_1,\pi_1)\to (M_2,\pi_2)$ is a Poisson map, in general, it 
\emph{does not} induce a Lie algebroid morphism of the associated cotangent
Lie algebroids, and hence it will not yield a Lie groupoid homomorphism. 

We will see now that we can still define $\Sigma$ on Poisson morphims
such that we obtain a functor from the \emph{Poisson category} to a 
certain \emph{symplectic groupoid ``category''}. This functor should not be 
confused with \emph{integration functor} $\G$ which we have studied before and
which goes from the category of Lie algebroids to the category of Lie groupoids. 
We will call $\Sigma$ the \textbf{symplectization functor} and will see that 
it entails a rich geometry, which is not present in the integration functor.

Henceforth, we will denote by \textbf{Poiss} the Poisson category, in which the objects are
the Poisson manifolds and the morphisms are the Poisson maps. We already know 
what the effect of $\Sigma$ on objects is, so let us look
at its effect on a Poisson morphism $\phi:(M_1,\pi_1)\to (M_2,\pi_2)$. 
In order to find out what the answer should be, let us
recall different ways of expressing the condition for a map to be Poisson:

\begin{exercise}
Let $\phi:(M_1,\pi_1)\to (M_2,\pi_2)$ be a smooth map between two Poisson manifolds. 
Show that the following conditions are equivalent:
\begin{enumerate}[(a)]
	\item The map $\phi$ preserves Poisson brackets: $\{f\circ\phi,g\circ\phi\}_1=\{f,g\}_2\circ\phi$.
	\item The Poisson bivectors are $\phi$-related: $\phi_*\pi_1=\pi_2$.
	\item $\Graph(\phi)\subset M_1\times\overline{M}_2$ is a coisotropic submanifold
	(\footnote{For a Poisson manifold $M$, the notation $\overline{M}$ means the same manifold with 
	the symmetric Poisson structure.}).
\end{enumerate}
(Hint: A submanifold $N$ of a Poisson manifold $(M,\pi)$ is called \emph{coisotropic} if
	$\pi^\sharp(TN)^0\subset TN$, where $(TN)^0\subset T^*M$ denotes the annihilator of $TN$.) 
\end{exercise}

Hence, to integrate a Poisson morphism, we just need to know what objects integrate
coisotropic submanifolds of a Poisson manifold. This is solved as follows:

\begin{theorem}
If $\G\tto C$ is a Lagrangian subgroupoid of a symplectic groupoid $\Sigma\tto M$, 
then $C\subset M$ is a coisotropic submanifold. Conversely, if $C$ is a coisotropic
submanifold of an integrable Poisson manifold $(M,\pi)$, there exists a Lagrangian
subgroupoid $\G\tto C$ of $\Sigma(M)$ that integrates $C$. 
\end{theorem}

\begin{proof}
Let us explain why coisotropic submanifolds integrate to Lagrangian
subgroupoids. Note that $C$ is a coisotropic submanifold of $(M,\pi)$ iff
its conormal bundle $\nu^*(C):=(TC)^0\subset T^*M$ is a Lie subalgebroid
of the cotangent Lie algebroid $T^*M$. Therefore, if $(M,\pi)$ is 
integrable, then there exists a source connected Lie subgroupoid $\G\tto C$
of the groupoid $\Sigma(M)\tto M$ that integrates $\nu^*(C)$. Now we claim
that the restriction of the symplectic form $\omega$ to $\G$ vanishes which,
combined with $\dim\G=1/2\dim\Sigma(M)$, implies that $\G$ is Lagrangian.

To prove our claim, we observe that for the canonical symplectic form $\omega_0$ 
on $T^*M$ the submanifold $\nu^*(C)\subset T^*M$ is Lagrangian. Using the explicit
expression (\ref{eq:lifted:symp:form}) for the symplectic form $\omega_{\text{can}}$
on $P(T^*M)$, we see immediately that space of paths $P(\nu^*(C))\subset P(T^*M)$ is 
isotropic. It follows that the symplectic form $\omega$ on the symplectic quotient 
$\Sigma(M)$ restricts to zero on the submanifold $\G=P(\nu^*(C))//P_0\Omega^1(M)$, 
as we claimed.
\end{proof}

We have now found what Poisson morphisms integrate to:

\begin{corollary}
\label{cor:mor}
Let $\phi:(M_1,\pi_1)\to (M_2,\pi_2)$ be a Poisson map between two integrable
Poisson manifolds. Then $\phi$ integrates to a Lagrangian subgroupoid
$\Sigma(\phi)\subset \Sigma(M_1)\times\overline{\Sigma(M_2)}$.
\end{corollary}

Let us introduce now the symplectic ``category'' of Alan Weinstein, which will be 
denoted by \textbf{Symp}. In this ``category'' the objects are the symplectic 
manifolds and the morphisms are the \emph{canonical relations}, which are defined 
as follows. 

\begin{definition}
If $(S_1,\omega_1)$ and $(S_2,\omega_2)$ are two symplectic manifolds, 
then a \textbf{canonical relation} from $S_1$ to $S_2$ is a Lagrangian 
submanifold $L\subset S_1\times\overline{S_2}$.
\end{definition}

If $L_1\in \Mor(S_1,S_2)$ and $L_2\in\Mor(S_2,S_3)$ are canonical relations, their composition
is the usual composition of relations:
\[ L_1\circ L_2:=\{(x,z)\in S_1\times S_3~|~\exists y\in S_2,\text{ with }(x,y)\in L_1
\text{ and }(y,z)\in L_2\}.\]

\begin{exercise}
Let $\phi:S_1\to S_2$ be a diffeomorphims between two symplectic manifolds. Show
that $\phi$ is a symplectomorphism iff its graph is a canonical relation.
\end{exercise}

Hence, the canonical relations enlarge the space of symplectomorphisms.
There is, however, a problem: $L_1\circ L_2$ may not be a smooth submanifold of 
$S_1\times \overline{S_3}$. One needs the two canonical relations to intersect cleanly, as explained
in the following exercise.

\begin{exercise}
Let $L_1\in \Mor(S_1,S_2)$ and $L_2\in\Mor(S_2,S_3)$ be canonical relations, and
denote by $L_1\star L_2$ their fiber product over $S_2$:
\[
\xymatrix{
L_1\star L_2\ar[r]\ar[d]& L_2\ar[d]\\
L_1\ar[r]& S_2
}
\]
Finally, let $p:S_1\times S_2\times S_3\to S_1\times S_3$ be the projection into
the first and last factor.
\begin{enumerate}
\item[(a)] Verify that $L_1\circ L_2=p(L_1\star L_2)$;
\end{enumerate}
Assume the \textbf{clean intersection property}: $L_1\star L_2$ is a 
manifold and its tangent spaces are fiber products of the tangent spaces of $L_1$ and $L_2$:
\[
\xymatrix{
T_{(s_1,s_2,s_3)}(L_1\star L_2)\ar[r]\ar[d]& T_{(s_2,s_3)}L_2\ar[d]\\
T_{(s_1,s_2)}L_1\ar[r]&T_{s_2}S_2
}
\]
\begin{enumerate}
\item[(b)] Show that $L_1\circ L_2$ is an (immersed) Lagrangian submanifold of $S_1\times S_3$.
\end{enumerate}
\end{exercise}

Therefore, in \textbf{Symp} composition is not always defined and Weinstein proposed 
to name it a ``category'', with quotation marks. One can make an analogy between this
``category'' and the category of hermitian vector spaces, as in the following table:
\vskip 10 pt

\begin{center}
  \begin{tabular}{|l||l|}
      \hline
      Symplectic ``category''  & hermitian category \\
      \hline
      Symplectic manifold $(S,\omega)$  & hermitian vector space $H$ \\
      $\overline{S}$ 	            & dual vector space $H^*$ \\
      $S_1\times S_2$             & tensor product $H_1\otimes H_2$\\
      $S_1\times \overline{S_2}$  & homomorphisms $\Hom(S_1,S_2)$\\
      a point $0$	                & $\Cc$\\
      \hline
  \end{tabular}
\end{center}
\vskip 10 pt

Just as an element in the hermitian vector space $H$ can be though of as morphism
from $\Cc\to H$, we can think of an element in the symplectic manifold $S$ as
a morphism in $\Mor(0,S)$, i.e., a Lagrangian submanifold of $S$. As a special
case, an element in the $\Hom$-object $S_1\times \overline{S_2}$ is a Lagrangian
submanifold, so we recover the canonical relations.

\begin{exercise}
In the symplectic ``category'' there exists an involution
which takes an object $S$ to itself $S^\dag:=S$ and a morphism 
$L\in\Mor(S_1,S_2)$ to the \emph{dual morphism} 
\[ L^\dag:=\{(y,x)\in S_2\times\overline{S_1}~|~(x,y)\in L\}\in\Mor(S_2,S_1).\]
Call a morphism \emph{unitary} if its dual equals its inverse: $L^\dag=L^{-1}$. 
Check that the unitary morphisms are just the symplectomorphisms.
\end{exercise}

The symplectic groupoids form a ``subcategory'' \textbf{SympGrp} of the 
symplectic ``category'', where the morphisms $\Mor(\G_1,\G_2)$ are the Lagrangian 
subgroupoids of $\G_1\times\overline{\G_2}$.  We can summarize the 
previous discussion by saying that $\Sigma$ is a covariant functor from 
\textbf{Poiss} to \textbf{SympGrp}. The reminder of this section is dedicated 
to the study of some basic properties of the functor $\Sigma$, namely, we look at its
effect on various geometric constructions in Poisson geometry.

\subsection{Sub-objects}
A sub-object in the Poisson category is just a Poisson submanifold $N$ of 
a Poisson manifold $(M,\pi)$. As the following example shows, a Poisson 
submanifold of an integrable Poisson manifold may fail to be integrable.

\begin{example}
Let $M=\mathfrak{so}(4)^*$ with its linear Poisson structure. Note that we have an
isomorphism 
\[ \mathfrak{so}(4)\simeq \mathfrak{so}(3)\oplus\mathfrak{so}(3)=\mathfrak{su}(2)\oplus\mathfrak{su}(2),\] 
so that we can think of $M$ as $\Rr^6=\Rr^3\times\Rr^3$ with the product Poisson structure,
where each factor has the Poisson structure (\ref{eq:PB:su2}). For the usual 
euclidean structure on $\Rr^6$, the function $f(x)=||x||^2$ is a Casimir, so that
the unit sphere
\[ \Ss^5=\{x\in\Rr^6~:~||x||=1\}, \]
is a Poisson submanifold of $M$. 
\begin{exercise}
Determine which points of $\Ss^5$ are regular points. Using Proposition \ref{var-area},
show that if $x=(x_1,x_2)\in\Ss^5$ is a regular point then:
\[ \NN_x\simeq \{n_1||x_1||+n_2||x_2||:n_1,n_2\in\Zz\}\subset\Rr.\]
\end{exercise}
Therefore, $\Ss^5$ is not an integrable Poisson manifold.
\end{example}

Henceforth, we will assume that both $(M,\pi)$ and the Poisson submanifold
$N$ are integrable. Then the inclusion $i:N\hookrightarrow M$ is a Poisson 
morphism which, according to Corollary \ref{cor:mor}, integrates to a 
Lagrangian subgroupoid $\Sigma(i)\subset\Sigma(N)\times\overline{\Sigma(M)}$. 
However, this is not the end of the story.

For a Poisson submanifold $N\subset M$, let us consider the set of equivalence classes of
cotangent paths that take their values in the restricted subbundle $T^*_NM$:
\[ \Sigma_N(M):=\{[a]\in\Sigma(M)~|~a:I\to T^*_NM\}.\]
This is a Lie subgroupoid of $\Sigma(M)$: it is the Lie subgroupoid that integrates 
the Lie subalgebroid $T^*_NM\subset T^*M$. Moreover, this subalgebroid is a coisotropic 
submanifold of the symplectic manifold $T^*M$, and it follows that $\Sigma_N(M)\subset \Sigma(M)$
is a coisotropic Lie subgroupoid. The fact that the closed 2-form $\omega$ on $\Sigma_N(M)$ is multiplicative
implies that if we factor by its kernel foliation, we still obtain a symplectic groupoid, and in fact
(see \cite{CrFe1}):
\[ \Sigma(N)\simeq \Sigma_N(M)/\Ker\omega.\]

These two constructions are related as follows:

\begin{theorem}
Let $i:N\hookrightarrow M$ be an integrable Poisson submanifold of an 
integrable Poisson manifold, and $\Sigma(i)\subset\Sigma(N)\times\overline{\Sigma(M)}$
the corresponding Lagrangian subgroupoid. For the restriction of the projections on each factor:
\[
\xymatrix{
&\Sigma(i)\ar[dl]_{\pi_1}\ar[dr]^{\pi_2}\\
\Sigma(N)&&\overline{\Sigma(M)}}\]
$\pi_2$ is a diffeomorphism onto the coisotropic subgroupoid
$\Sigma_N(M)\subset \Sigma(M)$ above, and the groupoid morphism:
\[ (\pi_1)\circ(\pi_2)^{-1}:\Sigma_N(M)\to \Sigma(N),\]
corresponds to the quotient map $\Sigma_N(M)\to\Sigma_N(M)/\Ker\omega\simeq\Sigma(N)$.
\end{theorem}

A very special situation happens when the exact sequence of Lie algebroids
\begin{equation}
\label{eq:split:Poisson}
\xymatrix{0\ar[r]& \nu^*(N)\ar[r]& T^*_NM\ar[r]& T^*N\ar[r]& 0}
\end{equation}
splits: in this case, the splitting $\phi:T^*M\to T^*_NM\subset T^*M$ 
integrates to a symplectic Lie groupoid homomorphism 
$\Phi:\Sigma(N)\to\Sigma(M)$, which realizes $\Sigma(N)$ as a symplectic subgroupoid
of $\Sigma(M)$. Poisson submanifolds of this sort maybe called \textbf{Lie-Poisson submanifolds},
and there are topological obstructions on a Poisson submanifold for this to happen.
The following exercise discusses the case of a regular symplectic leaf, so that
$\nu^*(N)$ is bundle of abelian Lie algebras.

\begin{exercise}\footnote{to clarify}
Let $N$ be an embedded Poisson submanifold of a Poisson manifold $(M,\pi)$, and choose
any splitting $\phi$ of the short exact sequence (\ref{eq:split:Poisson}). Denote
by $\nabla$ the connection defined by this splitting:
\[ \nabla_\al \gamma=[\phi(\al),\gamma],\quad \al\in\Omega^1(N),\gamma\in \Gamma(\nu^*(N)),\]
and by $\Omega$ the curvature of this splitting:
\[ \Omega(\al,\be):=[\phi(\al),\phi(\be)]-\phi([\al,\be]).\]
Show that:
\begin{enumerate}[(a)]
\item The connection $\nabla$ is independent of the choice of the 
splitting so that $\nu^*(N)$ is canonically a flat Poisson vector
bundle over $N$.
\item The curvature defines a Poisson cohomology 2-class $[\Omega]\in
H^2_{\Pi}(N,\nu^*(N))$ which is the obstruction for $N$ to be a Lie-Poisson
submanifold.
\item If $N$ is a Lie-Poisson submanifold, then $N$ is integrable.
\item If $N$ is simply connected, then $N$ is a
Lie-Poisson submanifold if and only if its monodromy group vanishes.
\end{enumerate}
\end{exercise}

\subsection{Quotients}

What we have just seen for sub-objects is typical: though the functor $\Sigma$ gives us
some indication of what the integration of a certain geometric construction is, there is
often extra geometry hidden in the symplectization. Another instance of this happens when 
one looks at quotients.

Let $G$ be a Lie group that acts smoothly by Poisson diffeomorphisms on a Poisson manifold 
$(M,\pi)$. We will denote the action by $\Psi:G\times M\to M$ and we will also write 
$\Psi(g,x)=g\cdot x$. For each $g\in G$, we set:
\[ \Psi_g:M\to M,\ x\mapsto g\cdot x,\]
so that each $\Psi_g$ is a Poisson diffeomorphism. 

Now we apply the functor $\Sigma$. For each $g\in G$, we obtain a Lagrangian
subgroupoid $\Sigma(\Psi_g)\subset \Sigma(M)\times\overline{\Sigma(M)}$. This 
Lagrangian subgroupoid is, in fact, the graph of a symplectic Lie groupoid
automorphism, which we denote by the same symbol $\Sigma(\Psi_g):\Sigma(M)\to\Sigma(M)$.
Also, it is not hard to check that 
\[ \Sigma(\Psi):G\times \Sigma(M)\to\Sigma(M), (g,[a])\mapsto g\cdot [a]:=\Sigma(\Psi_g)([a]),\]
defines a symplectic smooth action of $G$ on $\Sigma(M)$. Briefly, $\Sigma$
lifts a Poisson action $\Psi:G\times M\to M$ to a symplectic action 
$\Sigma(\Psi):G\times \Sigma(M)\to\Sigma(M)$ by groupoid automorphisms. However, this is not
the end of the story.

Let us look closer at how one lifts the action from $M$ to $\Sigma(M)$. First of all, recall 
that any smooth action $G\times M\to M$ has a lifted cotangent action $G\times T^*M\to T^*M$. 
This yields, by composition, an action of $G$ on cotangent paths: if $a:I\to T^*M$ is a 
cotangent path we just move it around
\[ (g\cdot a)(t):=g\cdot a(t).\]
The fact that the original action $G\times M\to M$ is Poisson yields that
(i) $g\cdot a$ is a cotangent path whenever $a$ is a cotangent path, and (ii) if
$a_0$ and $a_1$ are cotangent homotopic then so are the translated paths $g\cdot a_0$ and $g\cdot a_1$.
Therefore, we have a well-defined action of $G$ on cotangent homotopy classes and this is just the
lifted action: $\Sigma(\Psi_g)([a])=[g\cdot a]$.

Now we invoke a simple (but important) fact from symplectic geometry: for any 
action $G\times M\to M$ the lifted cotangent action $G\times T^*M\to T^*M$ is 
a hamiltonian action with equivariant momentum map $j:T^*M\to\gg^*$ given by:
\[ j:T^*M\to\gg^*,\ \langle j(\al_x),\xi\rangle:=\langle \al_x,X_\xi(x)\rangle,\]
where $X_\xi\in\X(M)$ is the infinitesimal generator associated with $\xi\in\gg$.
This yields immediately the fact that the lifted action $\Sigma(\Psi):G\times \Sigma(M)\to\Sigma(M)$
is also hamiltonian (\footnote{Recall again that, after all, the symplectic structure on
$\Sigma(M)$ comes from canonical symplectic structures on cotangent bundles.}). The 
equivariant momentum map $J:\Sigma(M)\to\gg^*$ for the lifted action is given by:
\begin{equation}
\label{eq:lifted:momentum:map}
\langle J([a]),\xi\rangle:=\int_0^1 j(a(t))\d t=\int_a X_\xi.
\end{equation}
Since each $X_\xi$ is a Poisson vector field, the last expression shows that only
the cotangent homotopy class of $a$ matters, and $J$ is indeed well-defined. 
Expression \ref{eq:lifted:momentum:map} means that we can see the momentum map of 
the $\Sigma(\Psi)$-action in two ways:
\begin{itemize}
\item It is the integration of the momentum map of the lifted cotangent action;
\item It is the integration of the infinitesimal generators along cotangent paths.
\end{itemize}
In any case, expression (\ref{eq:lifted:momentum:map}) for the momentum map 
shows that it satisfies the following additive property:
\[ J([a_0]\cdot [a_1])=J([a_0])+J([a_1]).\]
Hence $J$ is a groupoid homomorphism from $\Sigma(M)$ to the additive group $(\gg^*,+)$
or, which is the same, $J$ is differentiable \emph{groupoid 1-cocycle}. Moreover, 
this cocycle is \emph{exact} iff there exists a map $\mu:M\to\gg^*$ such that
\[ J=\mu\circ\t-\mu\circ\s,\]
and this happens precisely iff the original Poisson action $\Psi:G\times M\to M$ is hamiltonian
with equivariant momentum map $\mu:M\to\gg^*$. We summarize all this in the following
theorem:

\begin{theorem}
\label{thm:lifted:action}
Let $\Psi:G\times M\to M$ be a smooth action of a Lie group $G$ on a
Poisson manifold $M$ by Poisson diffeomorphisms. There exists a
lifted action $\Sigma(\Psi):G\time\Sigma(M)\to\Sigma(M)$ by symplectic groupoid
automorphisms. This lifted $G$-action is Hamiltonian and admits the
momentum map $J:\Sigma(M)\to\gg^*$ given by (\ref{eq:lifted:momentum:map}). 
Furthermore:
\begin{enumerate}[(i)]
\item The momentum map $J$ is $G$-equivariant and is a groupoid 1-cocycle.
\item The $G$-action on $M$ is hamiltonian with momentum map $\mu:M\to\gg^*$
if and only if $J$ is an exact cocycle.
\end{enumerate}
\end{theorem}

Let us now assume that the Poisson action $\Psi:G\times M\to M$ is 
proper and free. These assumptions guarantee that $M/G$ is a smooth manifold. The
space $C^\infty(M/G)$ of smooth functions on the quotient is naturally identified 
with the space $C^\infty(M)^G$ of $G$-invariant functions on $M$. Since the Poisson 
bracket of $G$-invariant functions is a $G$-invariant function, we have a quotient 
Poisson structure on $M/G$ such that the natural projection $M\to M/G$ is a Poisson map.

\begin{exercise}
Show that if $(M,\pi)$ is an integrable Poisson structure and $\Psi:G\times M\to M$ is 
proper and free Poisson action, then $M/G$ is also an integrable Poisson manifold.
\end{exercise}

So the question arises: what is the relationship between the symplectic groupoids 
$\Sigma(M)$ and $\Sigma(M/G)$?

First notice that if the original Poisson action $\Psi:G\times M\to M$ is 
proper and free, so is the lifted action $\Sigma(\Psi):G\times\Sigma(M)\to\Sigma(M)$. 
Therefore $0\in\gg^*$ is a regular value of the momentum map $J:\Sigma(M)\to\gg^*$.
Let us look at the symplectic quotient:
\[ \Sigma(M)//G:=J^{-1}(0)/G. \]
Since $J$ is a groupoid homomorphism, its kernel $J^{-1}(0)\subset\Sigma(M)$
is a Lie subgroupoid. Since $J$ is $G$-equivariant, the action leaves $J^{-1}(0)$ invariant and 
the restricted action is a free action by groupoid automorphisms. Hence, the groupoid structure
descends to a groupoid structure $\Sigma(M)//G\tto M/G$. 

\begin{exercise}
Show that $\Sigma(M)//G$ is a symplectic groupoid that integrates the 
quotient Poisson manifold $M/G$. 
\end{exercise}

In general, however, it is not true that:
\[ \Sigma(M/G)=\Sigma(M)//G,\]
so, in general, symplectization \emph{does not} commute with reduction. First of all, 
$J^{-1}(0)$ may not be connected, so that $\Sigma(M)//G$ may not have source connected 
fibers. Even if we restrict to $J^{-1}(0)^c$, the connected component of the identity section 
(so that the source fibers of $\Sigma(M)//G$ are connected) these fibers may have a non-trivial
fundamental group, as shown by the followig exercise:

\begin{exercise}
Consider the proper and free action of $\Ss^1$ on $M=\Cc^2-\{0\}$ 
given by $\theta\cdot(z_1,z_2)=(e^{i\theta}z_1,e^{-i\theta}z_2)$. For 
the canonical symplectic structure $\omega=\frac{i}{2}(\d z_1\wedge\d \overline{z}_1+
\d z_2\wedge\d \overline{z}_2)$, this $\Ss^1$-action on $M$ is Poisson. 
Determine both groupoids $\Sigma(M)//G$ and $\Sigma(M/G)$ and show that they are not isomorphic.
\end{exercise}

We will not get here into the subtelties of this problem.

\section{Notes}

The notion of a symplectic groupoid and its relation to Poisson geometry
appears first in Weinstein \cite{Wein1} and is discuss in greater detail
in \cite{CoDaWe}. The integration of Poisson manifolds presented has its 
roots in the work of Cattaneo and Felder \cite{CaFe} on the Poisson sigma model 
and was presented first by us in \cite{CrFe2}.

The main positive integrability result of \cite{CaFe} states that any 
Poisson structure on $\mathbb{R}^2$ is integrable. This, in turn, is a 
consequence of a general result due to Debord \cite{Deb} (see also 
\cite[Corollary 5.9]{CrFe1}) that states that a Lie algebroid with almost 
injective anchor is integrable. The results on the integrability of
regular Poisson manifolds and its relation to the variation of symplectic
areas are due to Alcade Cuesta and Hector \cite{AlHe}. The two examples we 
presented (the modified $\mathfrak{su}^*(2)$-bracket and the Heisenberg-Poisson 
manifolds) are both due to Weinstein \cite{Wein1,Wein4}.

The ``symplectic category'' appears early in the modern history of symplectic
geometry and is due to Weinstein \cite{Wein0}. Also, in his first note on
symplectic groupoids he observes that a Poisson morphism integrates to
a canonical relation. Coisotropic submanifolds are first studied in a systematic
way by Weinstein in \cite{Wein3}, in connection with Poisson groupoids. The problem 
of integrating coisotropic submanifolds was solved by Cattaneo and Felder in \cite{CaFe1}.
The integration of Poisson submanifolds	(and more general submanifolds) was discussed first 
by Xu in \cite{Xu0}, and then completed in \cite{CrFe2}. Actions of groups on groupoids, 
and some of the properties on lifted actions can be found, in one form or another, 
in \cite{CoDaMo,CoDaWe,MiWe,WeXu,Xu1}. This problem, as well as other issues 
such as non-free actions, convexity, etc, is the subject of ongoing research (see \cite{FeOrRa}).

The symplectization functor can be (and should be!) applied to many other constructions
in Poisson geometry. For example, one can look at the symplectization of Poisson
fibrations (in the sense of \cite{Vor}) and this yields groupoids which are 
symplectic fibrations (in the sense of \cite{GLS,MaSa}). Another example, is
provided by the connected sum construction in Poisson
geometry (see \cite{IM}) which out of two Poisson manifolds
$M_1$ and $M_2$ yields, under some conditions, a new Poisson manifold $M_1\# M_2$. A result
of the sort:
\[ \Sigma(M_1\# M_2)=\Sigma(M_1)\#\Sigma(M_2),\]
should be true (here, on the right-hand side one has a symplectic connected sum). A third 
example, is in the theory of Poisson-Nijenhuis manifolds where the application of 
the $\Sigma$ functor leads to a symplectic-Nijenhuis groupoid (\cite{Cr}), and this
is relevant both in the study of the generalized complex structures and of 
integrable hierarchies associated with a PN-manifold. The symplectization functor 
$\Sigma$, which we have just started understanding, should play an important role in 
many other problems in Poisson geometry.

\bibliographystyle{amsplain}

\begin{thebibliography}{11}

\bibitem{AlHe} F.~Alcade Cuesta and G.~Hector, Int\'egration
  symplectique des vari\'et\'es de Poisson r\'eguli\`eres,
  \emph{Israel J.~Math.}~\textbf{90} (1995), 125--165.
  
\bibitem{Mol} R.~Almeida and P.~Molino, Suites d'Atiyah et
  feuilletages transversalement complets, \emph{C.~R.~Acad.~Sci.~Paris
  S\'er.~I Math.~}\textbf{300} (1985), 13--15.

\bibitem{BFFGK} K.~Behrend, L.~Fantechi, W.~Fulton, L.~Goettsche and A.~Kresch,
  An Introduction to Stacks, book in preparation.

\bibitem{BeXu} K.~Behrend and P.~Xu, Differentiable Stacks and Gerbes,
  preprint \emph{math.DG/0605694} (2006).
  
\bibitem{Bra} W.~Brandt, \"Uber eine Verallgemeinerung des
  Gruppenbegriffes, \emph{Math. Ann.~}\textbf{96} (1926), 360--366.

\bibitem{Bry} R.~Bryant, Bochner-K\"ahler metrics, \emph{J.~Amer.~Math.~Soc.~}\textbf{14} 
	(2001), no. 3, 623--715.

\bibitem{CaWe} A.~Cannas da Silva and A.~Weinstein, \emph{Geometric
    Models for Noncommutative Algebras}, Berkeley Mathematics
    Lectures, vol.~\textbf{10}, American Math.~Soc.~, Providence,
    1999.
  
\bibitem{CaFe} A.S.~Cattaneo and G.~Felder, Poisson sigma models and
  symplectic groupoids, in \emph{Quantization of Singular Symplectic
  Quotients}, (ed. N.~P.~Landsman, M.~Pflaum, M.~Schlichenmeier),
  Progress in Mathematics \textbf{198} (2001), 41--73.
  
\bibitem{CaFe1} A.S.~Cattaneo and G.~Felder, Coisotropic submanifolds in Poisson 
	geometry and branes in the Poisson sigma model, 
	\emph{Lett.~Math.~Phys.~}\textbf{69} (2004), 157--175.
	
\bibitem{CoDaMo}
	M.~Condevaux, P.~Dazord, and P.~Molino, G\'eom\'etrie du moment, Travaux
  du S\'eminaire Sud-Rhodanien de G\'eom\'etrie, I, \emph{Publ.~D\'ep. Math.~Nouvelle
  S\'er.~B}, \textbf{88}, Univ.~Claude-Bernard, Lyon, 1988, pp.~131--160.
  
\bibitem{Connes}
	A.~Connes, \emph{Noncommutative Geometry}, Academic Press, 1984.
	
\bibitem{CoDaWe}
	A.~Coste, P.~Dazord, and A.~Weinstein, Groupo\"\i des symplectiques,
  Publications du D\'epartement de Math\'ematiques. Nouvelle S\'erie. A, Vol.\
  2, \emph{Publ.~D\'ep.~Math.~Nouvelle S\'er.~A}, \textbf{87}, Univ. Claude-Bernard, Lyon,
	1987, pp.~i--ii, 1--62.
  
\bibitem{Cr} M.~Crainic, Differentiable and algebroid cohomology, van
  Est isomorphisms, and characteristic classes,
  \emph{Comment.~Math.~Helv.~}\textbf{78} (2003), no. 4, 681--721.

\bibitem{Cr2} M.~Crainic, Generalized complex structures and 
	Lie brackets, preprint math.DG/0412097.

\bibitem{CrFe1} M.~Crainic and R.~L.~Fernandes, Integrability of Lie
  brackets, \emph{Ann.~of Math.~(2)} \textbf{157} (2003), 575--620.

\bibitem{CrFe2} M.~Crainic and R.~L.~Fernandes, Integrability of Poisson
  brackets, \emph{J.~Differential Geom.}~\textbf{66} (2004), 71--137.

\bibitem{DaHe} P.~ Dazord and G.~Hector, Int\'egration symplectique
  des vari\'et\'es de Poisson totalement asph\'eriques, in
  \emph{Symplectic Geometry, Groupoids and Integrable Systems, MSRI
  Publ.}, \textbf{20} (1991), 37--72.
  
\bibitem{Deb} C.~Debord, \emph{Feuilletages singulaires et
	groupo\"{\i}des d'holonomie}, Ph.~D. Thesis, Universit\'e Paul
	Sabatier, Toulouse, 2000.

\bibitem{DoLa} A.~Douady and M.~Lazard, Espaces fibr\'es en alg\`ebres
  de Lie et en groupes, \emph{Invent. Math.}~\textbf{1} (1966),
  133--151.

\bibitem{DuZu} J.-P.~Dufour, N.T.~Zung, \emph{Poisson structures and
    their normal forms}, Progress in Mathematics, \textbf{242},
    Birkh\"auser Verlag, Basel, 2005.  

\bibitem{DuKo} J.~Duistermaat and J.~Kolk, \emph{Lie Groups},
  Springer-Verlag Berlin Heidelberg, 2000.

\bibitem{Ehr} C.~Ehresmann, {\OE}uvres compl\`etes et
  comment\'ees. {I}-1,2. {T}opologie alg\'ebrique et g\'eom\'etrie
  diff\'erentielle, \emph{Cahiers Topologie G\'eom. Differentielle}
  \textbf{24} (1983) suppl.~1.
  
\bibitem{ELW} S.~Evens, J.-H.~Lu and A.~Weinstein, 
	Transverse measures, the modular class and a cohomology 
	pairing for Lie algebroids, \emph{Quart.~J.~Math.~Oxford (2)} 
	\textbf{50} (1999), 417--436.
  
\bibitem{Fer1} R.L.~Fernandes, Lie algebroids, holonomy and
  characteristic classe, \emph{Adv.~in Math.~}\textbf{170} (2002),
  119--179.

\bibitem{Fer2} R.L.~Fernandes, Connections in Poisson Geometry I:
  Holonomy and Invariants, \emph{J.~Differential Geometry}~\textbf{54}
  (2000), 303--366.

\bibitem{FeOrRa}
	R.L.~Fernandes, J.P.~Ortega and T.~Ratiu, Momentum maps in 
	Poisson geometry, in preparation.

\bibitem{Ginz} V.~Ginzburg, Grothendieck groups of Poisson vector
  bundles, \emph{J.~Symplectic Geometry} \textbf{1} (2001), 121--169.

\bibitem{Gro} A.~Grothendieck, S\'eminaire Henri Cartan, 13i\`eme
  ann\'ee: 1960/61.  Fasc. 1, Exp. 7, 9--13; Fasc. 2, Exp. 14--17,
  deuzi\`eme \'edition, Secr\'etariat math\'ematique, Paris, 1962.

\bibitem{GLS} 
	V.~Guillemin, E.~Lerman and S.~Sternberg, \emph{Symplectic 
  fibrations and multiplicity diagrams}, Cambridge University Press, 
  Cambridge, 1996. 

\bibitem{Haef} A.~Haefliger, Homotopy and integrability,
  Lecture Notes in Mathematics, Vol.~197, Springer, Berlin, 1971, 133--163.

\bibitem{HiMa} P.J.~Higgins and K.~Mackenzie, Algebraic constructions
  in the category of Lie algebroids, \emph{J. Algebra}~\textbf{129}
  (1990), 194--230.

\bibitem{Hueb1}
J.~Huebschmann, \emph{Duality for Lie-Rinehart algebras and the
modular class}, J.~reine angew. Math.~\textbf{510} (1999), 103--159.

\bibitem{IM} 
	A.~Ibort and D.~Martinez Torres, A new construction of
  Poisson manifolds, \emph{J. Symplectic Geometry}~\textbf{2} (2003),
  83--107.

\bibitem{Kar} M.~Karas\"ev, Analogues of objects of the theory of Lie
  groups for nonlinear Poisson brackets, \emph{Izv.~Akad.~Nauk
  SSSR Ser.~Mat.~}\textbf{50} (1986), no.3, 508--538, 638.

\bibitem{KuSp} A.~Kumpera, D.~Spencer, \emph{Lie equations. Vol. I:
    General theory}, Annals of Mathematics Studies, no.73, Princeton
    University Press, Princeton, N.J., 1972.  

\bibitem{Lang} S.~Lang, \emph{Differential Manifolds},
  Springer-Verlag, New-York, 1985.

\bibitem{Mack1} K.~Mackenzie, \emph{Lie Groupoids and Lie Algebroids
    in Differential Geometry}, London Math. Soc.~Lecture Notes Series
    \textbf{124}, Cambridge Univ.~Press, Cambridge, 1987.

\bibitem{Mack2} K.~Mackenzie, Lie algebroids and Lie pseudoalgebras,
  \emph{Bull.~London Math.~Soc.~}\textbf{27} (1995), 97--147.

\bibitem{Mack3} K.~Mackenzie, \emph{General theory of Lie groupoids
    and Lie algebroids}, London Mathematical Society Lecture Note
    Series \textbf{213}, Cambridge Univ.~Press, Cambridge, 2005.

\bibitem{MaXu} K.~Mackenzie and P.~Xu, Integration of Lie
  bialgebroids, \emph{Topology} \textbf{39} (2000), 445--467.
  
\bibitem{MaSa} D.~McDuff and D.~Salamon, \emph{Introduction to
  symplectic topology}, $2^\text{nd}$ edition, Oxford Mathematical
  Monographs, Oxford University Press, New York, 1998.

\bibitem{MiWe}
	K.~Mikami and A.~Weinstein, Moments and reduction for symplectic
  groupoids, \emph{Publ.~Res.~Inst.~Math.~Sci.~}\textbf{24} (1988), no.~1, 121--140.

\bibitem{Molino} P.~Molino, \'Etude des feuilletages transversalement
  complets et applications,
  \emph{Ann. Scient. \'Ec. Norm. Sup.}~\textbf{10} (1977), 289--307


\bibitem{MoMr} I.~Moerdijk and J.~Mr\v{c}un, \emph{Introduction to
  foliations and Lie groupoids.} Cambridge University Press,
  Cambridge, 2003.
  
\bibitem{Moer1} I.~Moerdijk, Orbifolds as groupoids: an introduction, in
  \emph{Orbifolds in mathematics and physics}, 205--222, Contemp. Math.,
  \textbf{310}, Amer. Math.~Soc., Providence, RI, 2002. 
   
\bibitem{Moer2} I.~Moerdijk, Introduction to the language of stacks and
  gerbes, University of Utrecht, 2002, preprint \emph{math.AT/0212266}.
  
\bibitem{MoMrc} I.~Moerdijk an  d J.~Mr\v{c}un, On integrability of 
	infinitesimal actions, preprint \emph{math.DG/0406558}.

\bibitem{Ni} V.~Nistor, Groupoids and the integration of Lie
  algebroids, \emph{J.~Math. Soc.~Japan} \textbf{52} (2000), 847--868.

\bibitem{NiXuWe} V.~Nistor, A.~Weinstein and P.~Xu, Pseudodifferential
  operators on differential groupoids,
  \emph{Pacific~J.~Math.}~\textbf{189} (1999), 117--152.
 
\bibitem{Palais} R.~Palais, A global formulation of the Lie theory of
  transformations, \emph{Mem.~Amer.~Math. Soc.}~\textbf{22}, 1957.

\bibitem{Pra1} J.~Pradines, Th\'eorie de Lie pour les groupo\"{\i}des diff\'erentiables. 
	Relations entre propri\'et\'es locales et globales, \emph{C.~R.~Acad.~Sci.~Paris, S\'erie A} 
	\textbf{263} (1966) 907--910.
	
\bibitem{Pra2} J.~Pradines, Th\'eorie de Lie pour les groupo\"{\i}des diff\'erentiables. 
	Calcul diff\'erenetiel dans la cat\'egorie des groupo\"{\i}des infinit\'esimaux, 
	\emph{C.~R.~Acad. Sci.~Paris, S\'erie A} \textbf{264} (1967) 245--248.

\bibitem{Pra3} J.~Pradines, G\'eom\'etrie diff\'erentielle au-dessus d'un groupo\"{\i}de,
	\emph{C.~R. Acad.~Sci.~Paris, S\'erie A} \textbf{266} (1968) 1194--1196.

\bibitem{Pra4} J.~Pradines, Troisi\`eme th\'eor\'eme de Lie pour les groupo\"{\i}des differentiables,
	\emph{C.~R.~Acad. Sci.~Paris, S\'erie A} \textbf{267} (1968) 21--23.
	
\bibitem{Sev} P.~Severa, Some title containing the words "homotopy"
  and "symplectic", e.g.~this one, preprint \emph{math.SG/0105080}.

\bibitem{Sing} I.~M.~~Singer and S.~Sternberg, The infinite groups of
  Lie and Cartan, \emph{J.~Analyse Math.}~\textbf{15} (1965), 1--114.

\bibitem{Suss} H.~Sussmann, Orbits of families of vector fields and
  integrability of distributions,
  \emph{Trans.~Amer. Math.~Soc.}~\textbf{180} (1973), 171--188.
  
\bibitem{Vain} A.~Va\u{\i}ntrob, Lie algebroids and homological vector fields,
	\emph{Uspekhi Mat.~Nauk} \textbf{52} (1997), no.~2(314), 161--162.
 
\bibitem{Vor} 
	Y.~Vorobjev, Coupling tensors and Poisson geometry 
  near a single symplectic leaf, \emph{Banach Center Publ.}~\textbf{54} 
  (2001), 249--274. 

\bibitem{Yorke} J.~Yorke, Periods of periodic solutions and the
  Lipschitz constant, \emph{Proc. Amer.~Math.~Soc.}~\textbf{22}
  (1969), 509--512.
  
\bibitem{Wein0}
	A.~Weinstein, Symplectic geometry, \emph{Bull.~Amer.~Math.~Soc.~(N.S.)} 
	\textbf{5} (1981), no.~1, 1--13. 

\bibitem{Wein1} A.~Weinstein, Symplectic groupoids and Poisson
  manifolds, \emph{Bull.~(New Series) Amer.~Math. Soc.}~\textbf{16}
  (1987), 101--104.
  
\bibitem{Wein3} A.~Weinstein, Coisotropic calculus and Poisson
  groupoids, \emph{J.~Math. Soc.~Japan} \textbf{40} (1988), 705--727. 

\bibitem{Wein4} A.~Weinstein, Blowing up realizations of Heisenberg-Poisson
  manifolds, \emph{Bull.~Sci. ~Math.}~\textbf{113} (1989), 381--406.
  
\bibitem{Wein5} A.~Weinstein, Linearization of regular proper
  groupoids, \emph{J.~Inst. Math.~Jussieu} \textbf{1} (2002), no.3, 493--511.

\bibitem{WeXu}
	A.~Weinstein and P.~Xu, Extensions of symplectic groupoids and
  quantization, \emph{J.~Reine Angew.~Math.~}\textbf{417} (1991), 159--189.

\bibitem{Xu0} P.~Xu, Dirac submanifolds and Poisson involutions, 
  \emph{Ann.~Sci.~\'Ecole Norm.~Sup.~}(4) \textbf{36} (2003), 403--430.

\bibitem{Xu1}
	P.~Xu, \emph{Symplectic groupoids of reduced Poisson spaces}, 
	\emph{C.~R.~Acad. Sci.~Paris, S\'erie I Math.} \textbf{314} (1992) 457--461.

\bibitem{Xu2}
	P.~Xu, \emph{Gerstenhaber algebras and BV-algebras in Poisson
	geometry}, Comm.~Math.~Phys.~\textbf{200} (1999), 545-560.
	
\bibitem{Xu3} P.~Xu, Morita equivalence of Poisson manifolds,
  \emph{Comm.~Math.~Phys.} \textbf{142} (1991), 493--509.

\bibitem{Xu4} P.~Xu, Morita equivalent symplectic groupoids, \emph{in}
  Symplectic geometry, groupoids, and integrable systems (Berkeley, CA, 1989),
  291--311, \emph{Math.~Sci.~Res.~Inst.~Publ.}, \textbf{20}, Springer,
  New York, 1991.

\bibitem{Zak}
	S.~Zakrzewski, Quantum and classical pseudogroups I and II, 
	\emph{Comm. Math.~Phys.~}\textbf{134} (1990), 347--395.

\bibitem{Zun} N.T.~Zung, Proper Groupoids and Momentum Maps: Linearization,
  Affinity and Convexity, preprint \emph{math.SG/0407208} (2004).
\end{thebibliography}
\def\lllll{}

\end{document}